%% file: Lu-2010-5.tex
\documentclass{amsart}

\usepackage{amsthm, amssymb}
\usepackage{amsfonts, amscd}
\usepackage{color}
\usepackage{epsf,epsfig}
\usepackage{graphicx}
\usepackage{epsfig,multicol}
\input xyv2
\input xy

\newtheorem{thm}{Theorem}[section]
\newtheorem{prop}[thm]{Proposition}
\newtheorem{lem}[thm]{Lemma}
\newtheorem{cor}[thm]{Corollary}

\theoremstyle{definition}

\newtheorem{defn}[thm]{Definition}

\newtheorem{example}{Example}

\theoremstyle{remark}

\newtheorem{rem}{Remark}

\def \R {\hbox{\rm I \kern -5.5pt R}}

\def \Hom {{\rm Hom}}

\def \R {{\bf R}\/}

\begin{document}
\title[Graphs and  $({\Bbb Z}_2)^k$-actions]{\large \bf Graphs and  $({\Bbb Z}_2)^k$-actions}
\renewcommand{\thefootnote}{\fnsymbol{footnote}}
\author[Zhi L\"u]{Zhi L\"u}
\thanks{Supported by  grants from NSFC (No. 10931005 and No. 10671034) and NSF of Shanghai (No. 10ZR1403600). }
\address{Institute of Mathematics, School of Mathematical Science, Fudan University, Shanghai,
200433, People's Republic of China.}
 \email{zlu@fudan.edu.cn}
\date{}
\begin{abstract}
Let $\mathcal{A}_n^k$ denote all nonbounding effective smooth
$({\Bbb Z}_2)^k$-actions  on  $n$-dimensional smooth closed
connected manifolds, each of which is cobordant to one with  finite fixed
 set. Motivated by GKM theory, we can associate to each action of $\mathcal{A}_n^k$ a $({\Bbb Z}_2)^k$-colored regular graph of valence $n$.
Together with the combinatorics of colored graphs, equivariant
cobordism and the tom Dieck-Kosniowski-Stong localization theorem,
\begin{enumerate}
\item[$\bullet$] we give a lower bound
for the number of fixed points of an action in $\mathcal{A}_n^k$,
which can become the best possible in some cases;
\item[$\bullet$] we  determine the existence and the equivariant cobordism
classification of all actions in $\mathcal{A}_n^k(h)$ with $h=3,4$,
where $\mathcal{A}_n^k(h)$ is the subset of
  $\mathcal{A}_n^k$, each of which is equivariantly cobordant to an effective $({\Bbb Z}_2)^k$-action
  fixing just $h$ isolated points (note: it is well-known that
  $\mathcal{A}_n^k(h)$ is empty if $h=1,2$);
\item[$\bullet$] we characterize the explicit relationships among  tangent representations at
fixed  points of each action in $\mathcal{A}_n^k(h)$ with $h=3,4$,
which actually give the explicit solution of the Smith problem in
such cases.
\end{enumerate}
As an application, we also study the minimum number of fixed points
of all actions in $\mathcal{A}_n^k$.
 \end{abstract}

\maketitle

\addcontentsline{toc}{chapter}{Bibliography}

\footnote[0]{}

\section{Introduction}

\subsection{Background} \label{s1-1}Using the work of Chang and Skjelbred \cite{cs},  Goresky,
Kottwitz and MacPherson in \cite{gkm} established the GKM theory, saying that  the
 equivariant cohomology of certain algebraic varieties with
  complex torus actions can explicitly be calculated in terms of their associated labeled
 regular  graphs (also called  {\em GKM graphs}). This indicates  that  there is an essential link between torus actions and the combinatorics of labeled regular graphs. Later on,
 Tolman and Weitsman \cite{t}  gave a simple proof of this
result in the symplectic setting, and showed that such GKM graphs can be produced from
a kind of $T^k$-manifolds (called the GKM manifolds).
  In subsequent works (see, e.g.,
  \cite{gz1}-\cite{gz5}),
  Guillemin and Zara developed the GKM theory combinatorially.
 In \cite{bgh}, Biss, Guillemin
and Holm studied the mod 2 version of the GKM theory, and showed
that a mod 2 GKM manifold can  still be associated to a unique
labeled graph (called the mod 2 GKM graph), and its equivariant
cohomology can be read out from this graph, where a mod 2 GKM
manifold $X$ is the real locus of a GKM symplectic manifold $M$ with
a Hamiltonian $T^k$-action, which is the fixed point set of an
anti-symplectic involution (compatible with the $T^k$-action) on
$M$. Note that $X$ naturally admits an action of $({\Bbb Z}_2)^k$
such that $X^{({\Bbb Z}_2)^k}=M^{T^k}$.  In recent years,
the GKM theory has  been further developed  in a variety of different
directions (see, e.g.,  \cite{cl}-\cite{c}, \cite{g}-\cite{gh1},
\cite{gh2}-\cite{gsz}, \cite{hl}, \cite{mmp}, \cite{z}). For
example, Guillemin and Holm in \cite{gh2} established the GKM theory
for torus actions with nonisolated fixed points, and Goldin and Holm
in \cite{gh1} generalized the GKM theory to the case in which the
one-skeleton  has dimension at most 4, so the mod 2 GKM theory may be
generalized to the case in which the one-skeleton has dimension at most
2 (see \cite{bgh}), where the one-skeleton of a $T^k$-manifold (resp. $({\Bbb Z}_2)^k$-manifold) $M$  consists of the points $p\in M$ with its isotropy subgroup $G_p$ having dimension $\geq k-1$.

\subsection{Motivation and Problems}An important feature of the GKM theory is that
equivariant topology can be associated with the combinatorics of
labeled regular graphs. In this paper, we shall further extend this
idea to  general $({\Bbb Z}_2)^k$-actions.

\vskip .2cm Throughout this paper assume that $G=({\Bbb Z}_2)^k$
unless stated otherwise. Suppose that  $(\Phi,
  M^n)$ is an   effective smooth $G$-action on a smooth connected closed
  manifold $M^n$ with a finite fixed point set (i.e.,  $0<\vert M^G\vert<+\infty$). Because
  the $G$-action is effective, the $G$-representation on the tangent space at a fixed point
  must be faithful and that implies  $n\geq k$. Conner and Floyd showed in \cite[Theorem 25.1]{cf}
   that if an involution on a closed manifold fixes a finite set,
  then the number of the fixed points must be even. Based upon this,
  we shall show in Section~\ref{s2} that $(\Phi,
  M^n)$ can always be associated to an $n$-valent regular graph $\Gamma_{(\Phi, M)}$ with $M^G$ as its
  vertex set such that there is a natural map $\alpha$ from the set
  $E_{\Gamma_{(\Phi, M)}}$ of all edges of $\Gamma_{(\Phi, M)}$ to
  $\Hom(G,{\Bbb Z}_2)$ with the following properties:
\begin{enumerate}
\item[(P1)] For each vertex $p$ of $\Gamma_{(\Phi, M)}$, the image
set $\alpha(E_p)$ spans $\Hom(G,{\Bbb Z}_2)$,  where $E_p$ denotes
the set of edges adjacent to $p$.
\item[(P2)] For each edge $e$ of $\Gamma_{(\Phi, M)}$,
$$\alpha(E_p)\equiv \alpha(E_q) \mod \alpha(e)$$
where $p$ and $q$ are two endpoints of $e$.
\end{enumerate}
Here we call the pair $(\Gamma_{(\Phi, M)},\alpha)$ a {\em
  $G$-colored graph} of $(\Phi, M^n)$. Indeed, generally $\Gamma_{(\Phi,
  M)}$ is not determined uniquely, but the set $\{\alpha(E_p)| p\in
  M^G\}$ is independent of the choice of $\Gamma_{(\Phi,
  M)}$. The $G$-colored graph $(\Gamma_{(\Phi, M)},\alpha)$ provides much important topological
  information of $(\Phi, M^n)$. Actually, since all irreducible real representations of $G$ can be identified with all elements
of $\mbox{Hom}(G,{\Bbb Z}_2)$, the set $\{\alpha(E_p)| p\in
  M^G\}$ gives all $G$-representations on
  tangent spaces at fixed points. In particular, $\{\alpha(E_p)| p\in
  M^G\}$ also determines a complete equivariant
  cobordism invariant $\mathcal{P}_{(\Phi,M^n)}$ of
  $(\Phi, M^n)$ where $\mathcal{P}_{(\Phi,M^n)}$ is
obtained by deleting all possible same pairs in $\{\alpha(E_p)| p\in
  M^G\}$ (see Theorem~\ref{ce}).
  This leads us to  consider the following questions:

\begin{enumerate}
 \item[(Q1)] What about the lower bound of the number $\vert M^G\vert$?
 \item[(Q2)] How to determine the existence and the equivariant cobordism classification of such $G$-actions $(\Phi,  M^n)$?
\item[(Q3)] How is the solution of the Smith problem for $(\Phi,
  M^n)$?
\end{enumerate}

\begin{rem}
The original Smith problem \cite{s1} says that if a smooth closed
manifold $X$ homotopic to a sphere admits an action of a finite
group $H$ such that the fixed point set is formed by only two isolated
points $u, v$, then are the tangent representations at $u$ and $v$
isomorphic? The problem was  solved by Atiyah-Bott
\cite{ab} for $H$ to be cyclic $p$-group ($p$ a prime) and by Milnor \cite{m} when $X$ is a homology sphere with a
semi-free $H$-action for $H$ an arbitrary connected compact Lie group. For further development  in this subject and counterexamples, see, e.g., \cite{group}, \cite{pr}, \cite{cs}.    More
generally, if $X$ is not restricted to be a homology sphere but a
closed manifold and if the number $\vert X^H\vert $ of fixed points
is greater than two, as stated in \cite{gz3}, the question of how
the tangent representations of $H$ at distinct fixed points are
related to each other is still open and is known as the  Smith
problem.
\end{rem}

 With respect to the above three questions, we shall pay more attention on   the case in which $G$-actions are
nonbounding. This is based upon the following  result by Stong in
\cite{s3}, which implies that if a $G$-action $(\Phi, M^n)$ with
$0<\vert M^G\vert<+\infty$ is bounding, then the number $\vert
M^G\vert$ may be reduced to be {\em zero} from the viewpoint of
cobordism.

 \begin{thm} [Stong] \label{thm1}
   Suppose $(\Phi, M^n)$ is  a smooth $G$-action
 on an $n$-dimensional smooth closed manifold $M^n$. Then $(\Phi, M^n)$ bounds equivariantly if and
 only if it is cobordant to a $G$-action $(\Psi, N^n)$
 with $N^G$ empty.
 \end{thm}

\begin{rem} Conner and Floyd \cite{cf} showed that when $k=1$, any
$G$-action
  $(\Phi, M^n)$ with  $\vert M^G\vert<+\infty$ always
  bounds equivariantly. This implies that $k\geq 2$ if $(\Phi, M^n)$ with  $\vert M^G\vert<+\infty$
  is nonbounding. Note that for $n=0$, there is a trivial action on a single
point which is nonbounding, and an action is nonbounding if and only
if it has an odd number of fixed points. Implicitly we will have
$n>0$ and $k>1$ throughout.
\end{rem}

 \subsection{Statements of main results} Now let $\mathcal{A}_n^k$ denote the set of all nonbounding
 effective smooth $G$-actions  on  smooth connected closed  $n$-manifolds such that each such $G$-action
 is cobordant to one fixing a finite set, where $n\geq k\geq 2$. Given a $G$-action  $(\Phi,
  M^n)$ in $\mathcal{A}_n^k$, let  $(\Gamma_{(\Phi, M)},\alpha)$ be a
  $G$-colored graph of $(\Phi, M^n)$. With respect to the question
  (Q1), we shall use
 $(\Gamma_{(\Phi, M)},\alpha)$ to find a lower bound of $|M^G|$ with a combinatorial argument, which
   is  stated as follows.

  \begin{thm} \label{thm2} Let $(\Phi,
  M^n)$ be a $G$-action in $\mathcal{A}_n^k$. Then the number $\vert M^G\vert$ is at least
$1+\lceil{n\over{n-k+1}}\rceil$ where $\lceil r\rceil$ denotes the
least integer greater than or equal to $r$.
\end{thm}

\begin{rem}\label{r1} First, we notice that the
bound estimated in Theorem~\ref{thm2} is the best possible
in some special cases. For instance, when $n=k$ or $2^{k-1}$, we
shall see from Examples~\ref{5-1}-\ref{5-2} in Section~\ref{s4} that
the bound is attainable,
   i.e., $\vert M^G\vert=1+\lceil{n\over{n-k+1}}\rceil$.
   However, when $k=2$ and $n$ is not a power of 2, we can find from Lemma~\ref{k=2} that
   $\vert M^G\vert$ is never equal to this bound. Second, since $n\geq k\geq 2$, we have that $\lceil{n\over{n-k+1}}\rceil\geq
   2$, so $|M^G|\geq 3$.

\end{rem}

\begin{cor}[cf.~{\cite[Theorem 31.3]{cf}}]
There is no manifold $M$ equipped with a $G$-action fixing exactly a single point in any equivalence class of
$\mathcal{A}_n^k$.
\end{cor}
\begin{proof}
If there is a $G$-action fixing exactly an isolated point, then this
action must be nonbounding according to the classical Smith Theorem.
\end{proof}

  With respect to the question (Q2), let $\mathcal{A}_n^k(h)$ be the subset of
  $\mathcal{A}_n^k$, each of which is equivariantly cobordant to an effective $G$-action
  fixing exactly $h$ isolated points.
  Then we know from Remark~\ref{r1} that
\begin{enumerate}
\item[$\bullet$]
 $\mathcal{A}_n^k(1)$ is empty.
\item[$\bullet$] $\mathcal{A}_n^k(2)$ is
also empty. Actually, any $G$-action fixing only two isolated points
bounds equivariantly (see \cite{ks}-\cite{l1}).
\end{enumerate}
  For the case $h\geq
  3$, as
  far as the author knows, the existence and
  equivariant cobordism classification of $G$-actions  has not been
  completely determined  yet except for the following cases:
\begin{enumerate}
\item[$\bullet$] The case $k=2$.
 In \cite[Theorem 31.2]{cf}, Conner and Floyd determined the existence of all $({\Bbb
Z}_2)^2$-actions of $\mathcal{A}_n^2$ (i.e., $\mathcal{A}_n^2$ is
nonempty if and only if $n\geq 2$ is even), and  they classified all
possible $({\Bbb Z}_2)^2$-actions of $\mathcal{A}_n^2$ up to
equivariant cobordism.
\item[$\bullet$] The case $k=n=3$.  In \cite{l3}, the author classified all  $({\Bbb Z}_2)^3$-actions
of $\mathcal{A}_3^3$  up to equivariant cobordism.
\end{enumerate}

  Using colored graphs and the tom Dieck-Kosniowski-Stong localization
  theorem (see Theorem~\ref{dks}),  we determine the existence of actions in $\mathcal{A}_n^k(h)$ with $h=3,
  4$,  and furthermore we
   completely classify up to equivariant cobordism  all possible $G$-actions in $\mathcal{A}_n^k(h)$ with $h=3, 4$.
    Let $(\phi_i, {\Bbb R}P^i)$ be the standard linear $({\Bbb Z}_2)^i$-action
    on ${\Bbb R}P^i$
    defined by $$\phi_i\big((g_1,...,g_i), [x_0,x_1,...,x_i]\big)=[x_0,g_1x_1,...,g_ix_i]$$
    where $(g_1,...,g_i)\in ({\Bbb Z}_2)^i$, which fixes $i+1$ isolated points.
        In Definitions~\ref{op1}-\ref{op2}, we will introduce two operations $\Omega$ and $\Delta$ on $G$-spaces, and
    apply $\Omega$-operation $k-i$ times and $\Delta$-operation $2^s$ times to $(\phi_i, {\Bbb R}P^i)$
    to obtain $\Delta^{2^s}\Omega^{k-i}(\phi_i, {\Bbb R}P^i)$. We then show that this is typical example
    when there are three or four isolated fixed points, as follows.

    \begin{thm} \label{thm3}
For $h=3$, we have that
\begin{enumerate}
\item[$(a)$]
$\mathcal{A}_n^k(3)$ is nonempty if and only if $k\geq 2$ and
$n=2^\ell\geq 2^{k-1}$.
\item[$(b)$] Given an integer $\ell\geq k-1$ with $n=2^\ell$,  each of $\mathcal{A}_{2^\ell}^k(3)$
is equivariantly cobordant to one of
$\big\{\sigma\Delta^{2^{\ell-k+1}}\Omega^{k-2}(\phi_2, {\Bbb
R}P^2)\big| \sigma\in\text{\rm GL}(k,{\Bbb Z}_2)\big\}$, where $\sigma\Delta^{2^{\ell-k+1}}\Omega^{k-2}(\phi_i, {\Bbb R}P^i)$ denotes the action obtained by applying an automorphism $\sigma\in \text{\rm GL}(k,{\Bbb Z}_2)$ of $({\Bbb Z}_2)^k$ to
    $\Delta^{2^{\ell-k+1}}\Omega^{k-2}(\phi_i, {\Bbb R}P^i)$.
\end{enumerate}
    \end{thm}

\begin{thm}\label{thm4}
For $h=4$, we have that
\begin{enumerate}
\item[$(a)$] $\mathcal{A}_n^k(4)$ is nonempty if and only if  $k\geq 3$
and $n$ is in the range $$\bigcup_{\ell \geq k-3}[3\cdot 2^{\ell},
5\cdot 2^{\ell}].$$
\item[$(b)$] Given an integer $\ell\geq k-3$ such that $3\cdot 2^{\ell}\leq
n\leq 5\cdot 2^{\ell}$,  each of $\mathcal{A}_n^k(4)$ is
equivariantly cobordant to one of $$\Big\{\sigma\Lambda^{\bf
v}\Delta^{2^{\ell-k+3}}\Omega^{k-3}(\phi_3,{\Bbb R}P^3)\big\vert
{\bf v}\in \mathcal{I}^k(n-3\cdot 2^\ell), \sigma\in \text{\rm
GL}(k,{\Bbb Z}_2)\Big\}$$ where $\mathcal{I}^k(t)$ denotes the set
of all lattices satisfying the equation $x_1+\cdots +x_{2^{k-3}}=t$
with each $0\leq x_i\leq 2^{\ell-k+4}$ in ${\Bbb R}^{2^{k-3}}$,
$\Lambda^{\bf v}$ is  a special operation on
$\Delta^{2^{\ell-k+3}}\Omega^{k-3}(\phi_3,$ ${\Bbb R}P^3)$ (see
Subsection~\ref{cl}), and $\sigma \Lambda^{\bf v}
\Delta^{2^{\ell-k+3}}\Omega^{k-3} (\phi_3,{\Bbb R}P^3)$ denotes the
action obtained by applying an automorphism $\sigma\in \text{\rm
GL}(k,{\Bbb Z}_2)$ to $\Lambda^{\bf
v}\Delta^{2^{\ell-k+3}}\Omega^{k-3}$ $(\phi_3,{\Bbb R}P^3)$.
 \end{enumerate}
\end{thm}

In theory, a colored graph $(\Gamma_{(\Phi,
  M)},\alpha)$ of a $G$-action $(\Phi, M^n)$ fixing a finite set  indicates a possible relationship among
representations on tangent spaces at fixed points, and the tom
Dieck-Kosniowski-Stong localization theorem  gives the algebraic
relationships among them (see Theorem~\ref{dks}). Here combining two
machineries leads us to obtain the {\em explicit} relationships for
the cases $|M^G|=3$ and 4, giving the solution of Smith problem. (Note that  the Smith problem for
$|M^G|=2$ can easily be solved by Theorem~\ref{dks}, i.e., the
tangent representations at two fixed points are isomorphic.) In
particular, this can also be associated with the geometrical
realization of abstract colored graphs.

 \vskip .2cm Let $\Gamma$ be a connected regular graph of valence
$n\geq k$. If there is  a map $\alpha: E_\Gamma\longrightarrow
\text{\rm Hom}(G,{\Bbb Z}_2)$ satisfying (P1) and (P2) as stated
before, then the pair $(\Gamma, \alpha)$ is called {\em an abstract
1-skeleton of type $(n,k)$}, and the set $\{\alpha(E_p)\big| p\in
V_\Gamma\}$ is called the {\em vertex-coloring set} of $(\Gamma,
\alpha)$, where $V_\Gamma$ denotes the set of all vertices of
$\Gamma$. We completely characterize the colored graphs of all actions in $\mathcal{A}_n^k(3)$ as follows.

\begin{thm} \label{thm5}
Let $(\Gamma, \alpha)$ be an abstract 1-skeleton of type $(n,k)$
such that $\Gamma$ contains exactly three vertices $p,q,r$. Then
$(\Gamma, \alpha)$ is realizable as a $G$-colored graph of some
$G$-action in $\mathcal{A}_n^k(3)$ if and only if the following
conditions are satisfied
\begin{enumerate}
\item[$(a)$] $k\geq 2$
and $n=2^\ell$ with $\ell\geq k-1$;
\item[$(b)$]
 there is a basis
$\{\beta_1,...,\beta_{k-1}, \gamma\}$ of $\Hom(G,{\Bbb Z}_2)$ such
that
 $$\alpha(E_p)=\hat{\beta}\cup\hat{\gamma},\ \ \alpha(E_q)=\hat{\beta}\cup\hat{\delta},\ \ \alpha(E_r)=
 \hat{\delta}\cup\hat{\gamma}$$
 where $\hat{\beta}$ is a multiset consisting of all $2^{k-2}$ different sums   with  same multiplicity
$2^{\ell-k+1}$ formed by the odd number
 of elements of $\beta_1,...,\beta_{k-1}$ (i.e., each sum appears exactly $2^{\ell-k+1}$ times
in $\hat{\beta}$ so  $|\hat{\beta}|=2^{\ell-1}$),
 $\hat{\gamma}=\{\gamma+\beta_1+\beta\ \vert
\beta\in\hat{\beta}\}$, and  $\hat{\delta}=\{\gamma+\beta\ \vert
\beta\in\hat{\beta}\}$.
 \end{enumerate}
\end{thm}

We also characterize the tangent representation sets of all actions in $\mathcal{A}_n^k(4)$ as follows.

 \begin{thm} \label{thm6}
 Let $(\Gamma, \alpha)$ be an abstract 1-skeleton of type
$(n,k)$ such that $\Gamma$ contains exactly four vertices $p,q,r,s$.
Then $\{\alpha(E_p), \alpha(E_q), \alpha(E_r), \alpha(E_s)\}$ is the
fixed data of some $G$-action  in $\mathcal{A}_n^k(4)$ if and only
if the following conditions are satisfied
\begin{enumerate}
\item[$(a)$] $k\geq 3$
and $n$ is in the range $3\cdot 2^\ell\leq n\leq 5\cdot 2^\ell$ for
some $\ell\geq k-3$;
\item[$(b)$] there is a basis $\{\beta_1,...,\beta_{k-2},$ $\gamma,
\delta\}$ of $\Hom(({\Bbb Z}_2)^k,{\Bbb Z}_2)$ such that
\begin{eqnarray*}
\alpha(E_p)=\hat{\beta}\cup\hat{\gamma}\cup\hat{\delta}\cup\hat{\omega},
&& \alpha(E_q)=\hat{\beta}\cup\hat{\eta}\cup
\hat{\varepsilon}\cup\hat{\omega},\\\alpha(E_r)=
 \hat{\gamma}\cup\hat{\varepsilon}\cup\hat{\lambda}\cup\hat{\omega}, &&\alpha(E_s)=
 \hat{\delta}\cup\hat{\eta}\cup\hat{\lambda}\cup\hat{\omega}\end{eqnarray*}
 where $\hat{\beta}$ is a multiset consisting of all $2^{k-3}$ different sums with the same multiplicity
$2^{\ell-k+3}$ formed by the odd number
 of elements of
$\beta_1,...,\beta_{k-2}$, and
$$\begin{cases}
\hat{\gamma}=\{\gamma+\beta_1+\beta|\beta\in\hat{\beta}\}\\
\hat{\delta}=\{\delta+\beta_1+\beta|\beta\in\hat{\beta}\}\\
\hat{\varepsilon}=\{\gamma+\beta|\beta\in\hat{\beta}\}\\
\hat{\eta}=\{\delta+\beta|\beta\in\hat{\beta}\}\\
\hat{\lambda}=\{\gamma+\delta+\beta_1+\beta|\beta\in\hat{\beta}\}\\
|\hat{\omega}|=n-3\cdot 2^{\ell}
\end{cases}$$
and  every element in $\hat{\omega}$ is chosen  in the $2^{k-3}$
different elements of  $\{\gamma+\delta+\beta|\beta\in\hat{\beta}\}$
and has multiplicity at most $2^{\ell-k+4}$.
\end{enumerate}
\end{thm}
\begin{rem}
  The reason why  we only characterize the tangent representation sets of all actions in $\mathcal{A}_n^k(4)$
     is because generally each action  $(\Phi, M^n)$ in $\mathcal{A}_n^k(4)$ may not determine a unique  regular graph
  $\Gamma_{(\Phi, M^n)}$. But when $n=3\cdot 2^\ell$, each action   in
  $\mathcal{A}_{3\cdot 2^\ell}^k(4)$ can determine a unique regular graph, so in
  this case, we can characterize the colored graphs of all actions in $\mathcal{A}_{3\cdot 2^\ell}^k(4)$ (see Theorem~\ref{abstr-graph}).
\end{rem}

The  paper is organized as follows. In Section~\ref{s2} we  review
some basic facts for $G$-representations, and then show how to get
colored regular graphs from $G$-actions. In Section~\ref{s3} we
reformulate Stong's result (i.e., Theorem~\ref{thm1}) into a
complete equivariant cobordism invariant in terms of tangent
representations at fixed points, and then review the tom
Dieck-Kosniowski-Stong localization theorem. In Section~\ref{sect4}
we introduce the $\Delta$-operation and the $\Omega$-operation. In
Section~\ref{s4} we  prove Theorem~\ref{thm2} using
combinatorial argument of colored graphs. In Section~\ref{s5} we
consider  $G$-actions fixing three isolated points, and prove Theorems~\ref{thm3} and \ref{thm5}.
 In Section~\ref{s6} by using the colored graph of $(\phi_3, {\Bbb R}P^3)$, we
  construct two examples with four fixed
 points, which will play an essential role in the argument of  the general case.
 Section~\ref{s7} is the most subtle and delicate arguments in this paper. We spend
 much time on determining the essential structure of the tangent
 representations at fixed points of actions of  $\mathcal{A}_n^k(4)$ in Subsection~\ref{ess-str}.
 Then we determine the existence and the equivariant cobordism
 classification of all actions of $\mathcal{A}_n^k(4)$ in Subsections~\ref{existence}-\ref{cl}, and
 completely
 characterize the relationships among tangent representations at fixed points of each
 action of $\mathcal{A}_n^k(4)$ in Subsection~\ref{char-co}.
 The proofs of
Theorems~\ref{thm4} and \ref{thm6} will be completed  in
Subsections~\ref{cl} and~\ref{char-co}, respectively. As an
application, in Section~\ref{s9} we also study the minimum number of
fixed points of all actions in $\mathcal{A}_n^k$. \vskip .2cm

The author expresses his gratitude to Professor R.E. Stong for his
valuable suggestions and comments, and also to Professor M. Masuda
for helpful conversations.

\section{$G$-representations and graphs of actions}\label{s2}

\subsection{$G$-representations}
Let $G=({\Bbb Z}_2)^k$, and let $\mbox{Hom}(G, {\Bbb Z}_2)$ be the
set of all homomorphisms $\rho: G\longrightarrow {\Bbb Z}_2=\{\pm
1\}$, which consists of $2^k$ distinct homomorphisms. Let $\rho_0$ denote the trivial element in $\mbox{Hom}(G,{\Bbb
Z}_2)$, i.e., $\rho_0(g)=1$ for all $g\in G$. Every irreducible real
representation of $G$ is one-dimensional and has the form
$$\lambda_{\rho}: G\times {\Bbb R}\longrightarrow {\Bbb R}$$ with
$\lambda_{\rho}(t, r) = \rho(t)\cdot r$ for some $\rho\in
\mbox{Hom}(G,{\Bbb Z}_2)$. It is well-known that there is a 1-1 correspondence
between all irreducible real representations of $G$ and all elements
of $\mbox{Hom}(G,{\Bbb Z}_2)$. $\mbox{Hom}(G,{\Bbb Z}_2)$ forms an
abelian group with addition given by
$(\rho+\sigma)(g)=\rho(g)\cdot\sigma(g)$, so it is also a
$k$-dimensional vector space over ${\Bbb Z}_2$ with standard basis
$\{\rho_1,...,\rho_k\}$ where $\rho_i$ is defined by mapping
$g=(g_1,...,g_i,...,g_k)$ to $g_i$.

\vskip .3cm Following \cite{co}-\cite{cf}, let $R_n(G)$  denote the
set generated by the isomorphism classes of $G$-representations of dimension $n$, which
naturally forms a vector space over ${\Bbb Z}_2$. Then
$$R_*(G)=\sum_{n\geq 0}R_n(G)$$ is a graded commutative algebra over
${\Bbb Z}_2$ with unit. The multiplication in $R_*(G)$ is given by
$[V_1]\cdot[V_2]=[V_1\oplus V_2]$. We can identify $R_*(G)$ with the
graded polynomial algebra over ${\Bbb Z}_2$ generated by
$\mbox{Hom}(G,{\Bbb Z}_2)$. Namely,  $R_*(G)$ is isomorphic to the
graded polynomial algebra ${\Bbb Z}_2[\rho_1,...,\rho_k]$. Based
upon this, throughout this paper we use the  convention that
 real representations of $G$ will be denoted by elements of ${\Bbb
Z}_2[\rho_1,...,\rho_k]$.

\subsection{Graphs of actions}
Now let $(\Phi, M^n)$ be an effective smooth $G$-action on a
 smooth closed connected manifold with $0<\vert M^G\vert<+\infty$.
 Here $(\Phi, M^n)$ is not necessarily restricted to be nonbounding.
 Then we are going  to construct a regular graph $\Gamma_{(\Phi, M)}$
associated with the action $(\Phi, M^n)$, such that the set of
vertices of the graph is $M^G$ and the valence of the graph is exactly
$n$.

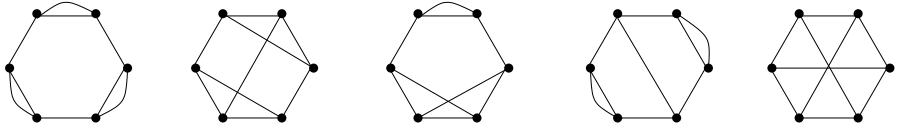
\begin{figure*}[h]
    \input{1.pstex_t}\centering
    \caption[a]{$3$-valent regular connected graphs with 6 vertices }\label{1}
\end{figure*}

\vskip .2cm
We now consider fixed sets of subgroups of $G$.
Given a nontrivial irreducible representation $\rho$, then the
kernel $\ker\rho$ is a subgroup of $G$ isomorphic to $({\Bbb
Z}_2)^{k-1}$. Let $C$ be a component of the fixed set of $\ker\rho$
and let $d=\dim C$. Then the fixed points of the group
$G/\ker\rho\cong {\Bbb Z}_2$ acting on $C$ will be fixed points of
$G$ acting on $M$. If $p$ is a fixed point  of $G/\ker\rho$ acting
on $C$, then $\rho$ occurs as a factor exactly $d$ times in the
tangent $G$-representation (i.e., a  monomial of degree $n$ in
${\Bbb Z}_2[\rho_1,...,\rho_k]$) at $p$ in $M$.
 Assuming $d>0$, Conner and Floyd have shown in \cite{cf}  that the number
of fixed points of  $G/\ker\rho$-action on $C$ must be even. Thus
we may always choose a connected regular graph $\Gamma_{\rho,C}$
with vertices the fixed points of $G/\ker\rho$-action on $C$ with
$d$ edges meeting at each vertex. If $d=1$, then $C$ is a circle
with precisely two fixed points, and we may choose an edge joining
these two fixed points. Clearly, in this case there is only one
choice for $\Gamma_{\rho,C}$. If $d>1$, let $p_1, p_2, ...,
p_{2i-1}, p_{2i}$ be all fixed points of $G/\ker\rho$-action on $C$.
In this case, there can be many choices for $\Gamma_{\rho,C}$ if
$i>1$. For example, when $d=i=3$,  one may gives five kinds of
different choices as shown in Figure~\ref{1}.
\vskip .2cm

Our graph $\Gamma_{(\Phi, M)}$ is now the union of all of these
subgraphs $\Gamma_{\rho, C}$ chosen for each $\rho$ and $C$ along fixed points of $M^G$. Because
the tangent $G$-representation at a fixed point $p$ is a monomial of
degree $n$ in ${\Bbb Z}_2[\rho_1,...,\rho_k]$, exactly $n$ edges
meet at $p$. So $\Gamma_{(\Phi, M)}$ is a regular graph of valence
$n$ with $M^G$ as its vertex set.

\vskip .2cm

\begin{rem}\label{r5}\hskip .1cm
\begin{enumerate}
 \item[a)] In general,  there can be many choices
for our graph $\Gamma_{(\Phi, M)}$.

\item[b)] Unlike GKM graphs, the orientation of associated regular
graphs in our case will not be considered.

\item[c)] Generally, there also may be several edges having the same
endpoints in $\Gamma_{(\Phi, M)}$, i.e., the number $\vert E_e\vert$
for some edge $e\in E_{\Gamma_{(\Phi, M)}}$ may not be 1, where
$E_{\Gamma_{(\Phi, M)}}$ denotes the set of all edges of
$\Gamma_{(\Phi, M)}$, and $E_e$ denotes the set of all edges joining
 two endpoints of $e$.\end{enumerate}
\end{rem}
On the uniqueness of $\Gamma_{(\Phi, M)}$, it is easy to see that
\begin{lem}
$\Gamma_{(\Phi, M)}$ is uniquely determined if  for each
$\rho$ and $C$,
$G/\ker\rho$-action on $C$ fixes only two isolated
points.
\end{lem}

 By the construction of $\Gamma_{(\Phi, M)}$, there is a natural
map
$$\alpha: E_{\Gamma_{(\Phi, M)}}\longrightarrow\Hom(G,{\Bbb Z}_2)$$
such that each edge of $\Gamma_{(\Phi, M)}$ is colored (or labeled)
by a nontrivial element in $\Hom(G,{\Bbb Z}_2)$. Actually, given an
edge $e$ in $E_{\Gamma_{(\Phi, M)}}$, there exists a nontrivial
element $\rho$ in $\Hom(G,{\Bbb Z}_2)$ and a component $C$ of
$M^{\ker\rho}$ such that $e$ is an edge of $\Gamma_{\rho, C}$. Then
$e$ is colored by $\rho$, namely $\alpha(e)=\rho$. The natural map
$\alpha$ has the following two basic properties:
\begin{enumerate}
\item[(P1)] For each vertex $p$ of $\Gamma_{(\Phi, M)}$, the image
set $\alpha(E_p)$ spans $\Hom(G,{\Bbb Z}_2)$.

\item[(P2)] For each edge $e$ of $\Gamma_{(\Phi, M)}$,
$$\alpha(E_p)\equiv \alpha(E_q) \mod \alpha(e)$$
where $p$ and $q$ are two endpoints of $e$.
\end{enumerate}

The property (P1) follows from the following argument. At a fixed
point $p$ in $M^n$, its tangent $G$-representation  is a monomial
$\prod_{e\in E_p}\alpha(e)$ of degree $n$ in ${\Bbb
Z}_2[\rho_1,...,\rho_k]$ where each $\alpha(e)\in \Hom(G,{\Bbb
Z}_2)$ is non-trivial by the definition of $\alpha$. Then the
condition that the action is effective means that $\alpha(e), e\in
E_p$ span $\Hom(G,{\Bbb Z}_2)$.
 The argument of (P2) is as follows: Since $p$ and $q$ are two
endpoints of $e$, there must be a connected component $C$ of the
fixed point set $M^{\ker \alpha(e)}$ such that $p, q\in C$. Thus,
both $p$ and $q$ have the same $\ker \alpha(e)$-representation.
Further, the tangent $G$-representations at $p$ and $q$ in $M$ are
the same when restricted to $\ker \alpha(e)$, and in particular, as
we have noted that $\alpha(e)$ occurs as a factor exactly $\dim C$
times in the tangent $G$-representations at both $p$ and $q$. If
$\sigma$ is another nontrivial irreducible representation occurring
at $\alpha(E_p)$ or $\alpha(E_q)$, then $\sigma$ and
$\sigma+\alpha(e)$ become the same nontrivial representation when
restricted to $\ker\alpha(e)$. Thus, $\alpha(E_p)\equiv \alpha(E_q)
\mod \alpha(e)$.

\vskip .3cm

The pair $(\Gamma_{(\Phi, M^n)}, \alpha)$ is called {\em a
$G$-colored graph} of $(\Phi, M^n)$. By $\mathcal{N}_{(\Phi, M^n)}$
we denote the set $\{\alpha(E_p)| p\in M^G\}$. Note that for $p\in
M^G$, generally $\alpha(E_p)$ may be a multiset. Then
$$\Big\{\prod_{e\in E_p}\alpha(e)\big| p\in M^G\Big\}$$ gives the tangent
representations at all fixed points.  Obviously, both $\alpha(E_p)$
and $\prod_{e\in E_p}\alpha(e)$ are determined by each other. Thus,
throughout this paper {\em $\mathcal{N}_{(\Phi, M^n)}$ will be
understood as all tangent $G$-representations at fixed points}.

\begin{rem}
 Clearly $\mathcal{N}_{(\Phi, M^n)}$ is
independent of the choice of $G$-colored graphs.
\end{rem}

\begin{example}
Let $(\phi_n, {\Bbb R}P^n)$ with $n\geq 2$ be the standard linear
$({\Bbb Z}_2)^n$-action
    on ${\Bbb R}P^n$
    defined by $$\phi_l\big((g_1,...,g_n), [x_0,x_1,...,x_n]\big)=[x_0,g_1x_1,...,g_nx_n], \ (g_1,...,g_n)\in ({\Bbb Z}_2)^n, $$
     which fixes $n+1$ isolated points $[0,...,0,1,0,...,0]$ with 1
in the $i$-th place for $i=0,1,...,n$. It is not difficult to see
that $(\phi_n, {\Bbb R}P^n)$ determines a unique regular graph
$\Gamma_{(\phi_n, {\Bbb R}P^n)}$, which is just the 1-skeleton of an
$n$-simplex, and the ${{n+1}\choose 2}$ edges of $\Gamma_{(\phi_n,
{\Bbb R}P^n)}$ are colored by $\rho_1, ..., \rho_n, \rho_i+\rho_j,
1\leq i<j\leq n$, respectively. Actually, $(\phi_n, {\Bbb R}P^n)$
 is a small  cover over an $n$-simplex $\Delta^n$, and all its possible orbit types can be read out from $\Delta^n$ (see \cite{dj}).  When $n=2, 3$, the colored graph $\Gamma_{(\phi_n,
{\Bbb R}P^n)}$ is shown in Figure~\ref{n2}.
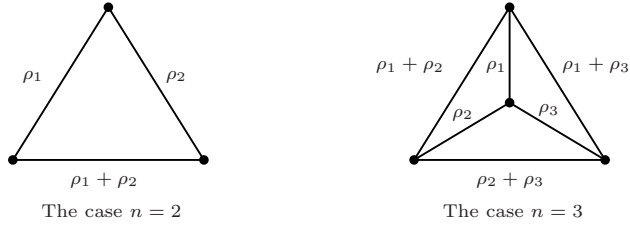
\begin{figure*}[h]
    \input{n2.pstex_t}\centering
    \caption[a]{Colored graphs for the cases $n=2, 3$. }\label{n2}
\end{figure*}
Furthermore, the diagonal action on two copies of $(\phi_2, {\Bbb
R}P^2)$ and the twist involution on the product ${\Bbb R}P^2\times{\Bbb R}P^2$ give a
$({\Bbb Z}_2)^3$-action on ${\Bbb R}P^2\times{\Bbb R}P^2$ fixing
three fixed points, whose colored graph is shown in Figure~\ref{n3}
\begin{figure*}[h]
    \input{n3.pstex_t}\centering
    \caption[a]{The colored graph of the $({\Bbb Z}_2)^3$-action on ${\Bbb R}P^2\times{\Bbb R}P^2$. }\label{n3}
\end{figure*}
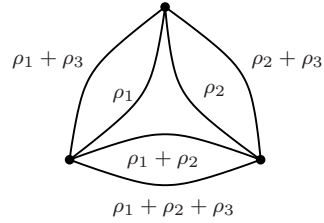
\end{example}

\begin{rem}\label{r6}
In  mod 2 GKM theory, as noted in \cite[Remark 5.9]{bgh}, when $(\Phi, M^n)$ satisfies the
conditions that (1) $(\Phi, M^n)$ is equivariantly formal; (2) $M^G$
is finite; (3) the isotropy weights of the tangent representations
at each fixed point are all distinct and non-zero (i.e., for each
$\rho$ and $C$, $\dim C=1$, so $|E_e|=1$ for each edge
$e\in\Gamma_{(\Phi, M)}$),  its
equivariant cohomology can be explicitly read out from
$(\Gamma_{(\Phi, M)},\alpha)$ as follows:
$$H^*_G(M^n;{\Bbb Z}_2)\cong \{f: M^G\longrightarrow
{\Bbb Z}_2[\rho_1,...,\rho_k]\big| f(p)\equiv f(q) \mod \alpha(e)
\text{ for } e\in E_p\cap E_q\}.$$ In addition, as mentioned in
Subsection~\ref{s1-1},  Biss, Guillemin and Holm in \cite{bgh} also
generalized the mod 2 GKM theory   to the case where the
one-skeleton has dimension at most 2 (i.e., for each $\rho$ and $C$,
$\dim C\leq 2$).

\vskip .2cm  A natural question is that if $(\Phi, M^n)$ with $M^G$ a finite set satisfies the conditions
that (1) $(\Phi, M^n)$ is equivariantly formal; (2) for each $\rho\in \Hom(G,{\Bbb Z}_2)$ and possible components $C$ with $\dim C>2$ of $M^{\ker\rho}$,  the action of $G/\ker\rho$ on  $C$ fixes only two isolated points,
 whether can its equivariant cohomology  be
explicitly read out from its colored graph $(\Gamma_{(\Phi,
M)},\alpha)$? But this is beyond the scope of the current paper.
\end{rem}

\section{A complete equivariant cobordism invariant and the tom Dieck-Kosniowski-Stong localization
theorem}\label{s3}

In this section, we first reformulate Stong's result (i.e.,
Theorem~\ref{thm1}) into a complete equivariant cobordism invariant
in terms of tangent representations at fixed points, and then review
the tom Dieck-Kosniowski-Stong localization theorem (\cite{d},
\cite{ks}).

\vskip .3cm

Throughout the following,  assume that $(\Phi, M^n)$ is an effective
smooth $G$-action on a smooth closed  connected manifold with
$0<\vert M^G\vert<+\infty$, and $(\Gamma_{(\Phi, M^n)}, \alpha)$ is
a colored graph of $(\Phi, M^n)$.

\subsection{A complete equivariant cobordism invariant}

\begin{lem} \label{lm1}
 $(\Phi, M^n)$  is equivariantly cobordant to a $G$-action
$(\Psi, N)$ such that either $\mathcal{N}_{(\Psi, N)}$ is empty or
it is non-empty but all elements of $\mathcal{N}_{(\Psi, N)}$ are
different.
\end{lem}

\begin{proof}
If $p$ and $q$ are two fixed points with $\alpha(E_p)=\alpha(E_q)$,
then one can cut out neighborhoods of $p$ and $q$, each of which
looks like the disc in the associated representation space. One then
glues the resulting boundaries together to obtain a new action,
which is cobordant to the action $(\Phi, M^n)$ by \cite{s3}. This
reduces the number of fixed points by two. One can carry out this
procedure until one obtains the required action $(\Psi, N)$.
\end{proof}

Obviously, $\mathcal{N}_{(\Psi, N)}$ is determined uniquely by
$(\Phi, M^n)$.  Set $$\mathcal{P}_{(\Phi, M^n)}:=
\mathcal{N}_{(\Psi, N)}. $$ Here one calls $\mathcal{P}_{(\Phi,
M^n)}$ the {\em prime tangent $G$-representation set} of $(\Phi,
M^n)$.  By Theorem~\ref{thm1} and Lemma~\ref{lm1}, one has

\begin{thm} \label{ce}
$\mathcal{P}_{(\Phi,
M^n)}$ is a complete equivariant cobordism invariant of $(\Phi,
M^n)$.
\end{thm}

\subsection{The tom Dieck-Kosniowski-Stong localization theorem} In \cite{d}, tom Dieck showed
that the equivariant cobordism class of $(\Phi, M^n)$ is completely
determined by its equivariant Stiefel-Whitney characteristic
numbers, and in particular, the existence of $(\Phi, M^n)$ can be
determined in terms of its fixed data. Later on, Kosniowski and
Stong \cite{ks} gave a more precise formula for the characteristic
numbers of $M^n$ in terms of the fixed data. Combining their works,
one has the following localization theorem.
\begin{thm}  [tom Dieck-Kosniowski-Stong]\label{dks} Let $(\Gamma, \alpha)$ be an abstract 1-skeleton of type $(n,k)$.
 Then the necessary and sufficient condition that
$\{\alpha(E_p)\big| p\in V_\Gamma\}$ is the fixed data of some
$({\Bbb Z}_2)^k$-action $(\Phi, M^n)$ is that for any symmetric
polynomial function $f(x_1,...,x_n)$ over ${\Bbb Z}_2$,
\begin{equation} \label{formula-tks1}\sum_{p\in V_\Gamma}{{f(\alpha(E_p))}\over{\prod_{e\in
E_p}\alpha(e)}}\in {\Bbb Z}_2[\rho_1,...,\rho_k]\end{equation} where
$V_\Gamma$ denotes the set of all vertices of $\Gamma$, and
$f(\alpha(E_p))$ means that $x_1,...,x_n$ in $f(x_1,...,x_n)$ are
replaced by all elements in $\alpha(E_p)$.
\end{thm}

\begin{rem} Originally, as stated in \cite{d} and \cite{ks}, the
formula (\ref{formula-tks1}) should be written as the following form
\begin{equation} \label{formula-tks2}\sum_{p\in V_\Gamma}{{f(\chi^G(p))}\over{\chi^G(p)}}\in H^*(B({\Bbb Z}_2)^k;{\Bbb Z}_2)
={\Bbb Z}_2[t_1,...,t_k]
\end{equation}
 with each  $t_i\in H^1(B({\Bbb Z}_2)^k;{\Bbb Z}_2)$, where
$\chi^G(p)$ denotes the equivariant Euler class of the tangent
representation at $p$, which is a monomial of degree $n$ in
$H^*(B({\Bbb Z}_2)^k;{\Bbb Z}_2)$ and $f(\chi^G(p))$ means that
$x_1,...,x_n$ in $f(x_1,...,x_n)$ are replaced by $n$ factors in
$\chi^G(p)$. If $\{\alpha(E_p)\big| p\in V_\Gamma\}$ is the fixed
data of some $({\Bbb Z}_2)^k$-action $(\Phi, M^n)$, then the
polynomial $\sum_{p\in V_\Gamma}{{f(\chi^G(p))}\over{\chi^G(p)}}\in
H^*(B({\Bbb Z}_2)^k;{\Bbb Z}_2)$  actually means an equivariant
Stiefel-Whitney number of the action. In other words, if we formally
write the the equivariant total Stiefel-Whitney class of the tangent
bundle $\tau(M)$ as $w^G(\tau(M))=\prod_{i=1}^n(1+x_i)$, then the
equivariant Stiefel-Whitney number $f(x_1, ..., x_n)[M]$ can be
calculated by the following formula
\begin{equation}\label{tks3}
f(x_1, ..., x_n)[M]=\sum_{p\in
V_\Gamma}{{f(\chi^G(p))}\over{\chi^G(p)}}\end{equation} where $[M]$
denotes the fundamental homology class of $M$.  For more details,
see \cite{d} and \cite{ks}. The formula~(\ref{tks3}) is an analogue
of the Atiyah-Bott-Berlin-Vergne formula for the case of torus
actions (see \cite{ab1} and \cite{bv}). Since $H^*(B({\Bbb
Z}_2)^k;{\Bbb Z}_2)$ is isomorphic to ${\Bbb
Z}_2[\rho_1,...,\rho_k]$, each class $\chi^G(p)$ uniquely
corresponds to its tangent representation $\alpha(E_p)$, so  we can
give another description of the formula (\ref{formula-tks2}) in
terms of $G$-representations, as stated in (\ref{formula-tks1}).
\end{rem}

\section{Some operations on $G$-actions}\label{sect4}

 Throughout the following, let $(\Phi,
M^n)$ be a $({\Bbb Z}_2)^k$-action on a closed manifold $M$ (note
that here the fixed set of the action is not necessarily restricted
to be finite).

\subsection{Diagonal action}
\begin{defn}\label{op1}
Given an integer $i\geq 1$,   the action
$\underbrace{\Phi\times\cdots\times\Phi}_i$ on
$\underbrace{M^n\times\cdots\times M^n}_i$ defined by $$\big(g,
(x_1, ..., x_i)\big)\longmapsto \big(\Phi(g, x_1), ..., \Phi(g,
x_i)\big)$$ is called an {\em $i$-multi-diagonal action} on $(\Phi,
M^n)$, denoted by $\Delta^i(\Phi, M^n)$.
\end{defn}

For each $i$, $\Delta^i(\Phi, M^n)$ is still a $({\Bbb
Z}_2)^k$-manifold and has dimension $i\cdot n$ such that its fixed
point set is $\underbrace{M^{({\Bbb Z}_2)^k}\times\cdots\times
M^{({\Bbb Z}_2)^k}}_i$. Let $$\xi\longrightarrow M^{({\Bbb
Z}_2)^k}=\bigsqcup_{0\leq j\leq n} \xi^{n-j}\longrightarrow F^j$$ be
the normal bundle of $M^{({\Bbb Z}_2)^k}$ in $M^n$, where $F^j$
denotes the $j$-dimensional fixed part of $M^{({\Bbb Z}_2)^k}$,
which consists of a disjoint union of all $j$-dimensional connected
components, and  if $M^{({\Bbb Z}_2)^k}$ contains no $j$-dimensional
fixed part, then $F^j$ is chosen to be empty. Now if $i$ is a
power of 2, say $2^s$, since ${{2^s}\choose h}\equiv 0\mod 2$ for
any $0<h<2^s$, then it is easy to see that $\Delta^i(\Phi, M^n)$ is
equivariantly cobordant to a $({\Bbb Z}_2)^k$-action $(\Psi, N)$
having fixed point set $\bigsqcup_{0\leq j\leq n}F'^j$, where
$F'^j=\underbrace{F^j\times\cdots\times F^j}_i$. This gives the
following result.

\begin{lem}\label{dia}
If $i$ is a power of 2, then $\Delta^i(\Phi, M^n)$ is equivariantly
cobordant to a $({\Bbb Z}_2)^k$-action $(\Psi, N)$ whose fixed point
set is the diagonal copy of $M^{({\Bbb Z}_2)^k}$ contained in
$$\underbrace{M^{({\Bbb Z}_2)^k}\times\cdots\times M^{({\Bbb
Z}_2)^k}}_i.$$
\end{lem}

\begin{rem}
If $i$ is a power of 2 and $M^{({\Bbb Z}_2)^k}$ is a finite set,
then $(\Psi, N)$ has the same number of fixed points as $(\Phi,
M^n)$.
\end{rem}

\subsection{$\Omega$-operation}
\begin{defn}\label{op2}
 An $\Omega$-{\em operation} on $(\Phi, M^n)$ means that $\Omega$
maps $(\Phi, M^n)$ into a $({\Bbb Z}_2)^{k+1}$-action on $M^n\times
M^n$, denoted by $\Omega{(\Phi, M^n)}$, which is given by  the
diagonal action $\Phi\times \Phi$ together with the involution swapping the factors of $M^n\times M^n$.
The fixed set of $\Omega{(\Phi, M^n)}$ is the copy of the fixed set
of $M^n$ in the diagonal copy of $M^n$ contained in $M^n\times M^n$.
\end{defn}

 $\Omega{(\Phi, M^n)}$ has many same properties as $(\Phi, M^n)$. For example, the fixed point set of $\Omega{(\Phi,
 M^n)}$ has the same dimension as that of $(\Phi, M^n)$. Also, as noted in \cite{l2},
 if the fixed set
of $(\Phi, M^n)$ possesses the linear independence, then this also
is so for $\Omega{(\Phi, M^n)}$.

\vskip .2cm
 Applying the $\Omega$-operation $l$ times to $(\Phi, M^n)$ gives a $({\Bbb Z}_2)^{k+l}$-action denoted by
$\Omega^l{(\Phi, M^n)}$ such that $\dim \Omega^l{(\Phi, M^n)}=2^ln$,
and its fixed point set is
$$\big\{(\underbrace{p, ..., p}_{2^l})\big|p\in M^{({\Bbb Z}_2)^k}\big\}.$$
Clearly, if $(\Phi, M^n)$ has exactly a finite fixed set, then for
each $l\geq 1$,  $\Omega^l{(\Phi, M^n)}$  has a finite fixed set,
and in particular, $\Omega^l{(\Phi, M^n)}$ has the same number of
fixed points as $(\Phi, M^n)$.

\begin{lem} \label{twist}
Suppose that $(\Phi, M^n)$ is a nonbounding $({\Bbb Z}_2)^k$-action
fixing a finite set. Then  for each $l\geq 1$,  $\Omega^l{(\Phi,
M^n)}$ is nonbounding.
\end{lem}
\begin{proof}
It suffices to show that for $l=1$, $\Omega{(\Phi, M^n)}$ is
nonbounding. Let $(\Gamma_{(\Phi, M^n)}, \alpha)$ be a colored graph
of $(\Phi, M^n)$, and let $p$ be a fixed point  of $(\Phi, M^n)$.
Then its tangent representation at $p$ is $\alpha(E_p)$. Let
$\Lambda$ and $\bar{\Lambda}$ be two subsets of $\text{Hom}(({\Bbb
Z}_2)^{k+1}, {\Bbb Z}_2)$ such that both $\Lambda$ and
$\bar{\Lambda}$ are isomorphic to $\text{Hom}(({\Bbb Z}_2)^k, {\Bbb
Z}_2)$ as ${\Bbb Z}_2$ vector spaces, and each $\delta_\rho$ in
$\Lambda$ is $\rho\in \text{Hom}(({\Bbb Z}_2)^k, {\Bbb Z}_2)$ on
$({\Bbb Z}_2)^k$ and 1 on the new ${\Bbb Z}_2$ generator $t_{k+1}$
and each $\bar{\delta}_\rho$ in $\bar{\Lambda}$ is $\rho \in
\text{Hom}(({\Bbb Z}_2)^k, {\Bbb Z}_2)$ on $({\Bbb Z}_2)^k$ and $-1$
on the new ${\Bbb Z}_2$ generator $t_{k+1}$, where $({\Bbb
Z}_2)^{k+1}$ is identified with $({\Bbb Z}_2)^k\oplus
\text{Span}\{t_{k+1}\}$. Now it is easy to see that the tangent
representation at the corresponding fixed point $(p,p)$ in
$\Omega{(\Phi, M^n)}$ is $\{\delta_\rho|\rho\in \alpha(E_p)\}\cup
\{\bar{\delta}_\rho| \rho\in \alpha(E_p)\}$. Since $(\Phi, M^n)$ has
a nonempty  prime tangent representation set, it is also so for
$\Omega{(\Phi, M^n)}$. Then the lemma follows from Theorem~\ref{ce}.
\end{proof}

\subsection{Automorphisms of $({\Bbb Z}_2)^k$}
Given an automorphism $\sigma\in\text{\rm GL}(k,{\Bbb Z}_2)$ of
$({\Bbb Z}_2)^k$, the action $({\Bbb Z}_2)^k\times
M^n\longrightarrow M^n$ defined by
$$(g, x)\longmapsto \Phi(\sigma(g), x)$$
is called the {\em $\sigma$-induced action} of $(\Phi, M^n)$,
denoted by $(\sigma\Phi, M^n)$ or $\sigma(\Phi, M^n)$.

\vskip .2cm

Note that $(\sigma\Phi, M^n)$ is weakly equivariantly homeomorphic
to $(\Phi, M^n)$, but it may not be equivariantly cobordant to
$(\Phi, M^n)$ since generally $\sigma$ will change the tangent
representation set of $(\Phi, M^n)$.

\begin{rem}\label{automorphism}
If $(\Phi, M^n)$ has exactly a finite fixed set, let
$(\Gamma_{(\Phi, M^n)}, \alpha)$ is a colored graph of $(\Phi,
M^n)$, then it is easy to see that $(\Gamma_{(\Phi, M^n)},
\sigma\alpha)$ is a colored graph of $(\sigma\Phi, M^n)$, where
$\sigma\alpha: E_{\Gamma_{(\Phi, M^n)}}\longrightarrow \text{Hom}(G,
{\Bbb Z}_2)$ is defined by
$$\sigma\alpha(e)(g)=\alpha(e)(\sigma(g))$$ for $e\in
E_{\Gamma_{(\Phi, M^n)}}$ and $g\in G$.\end{rem}

\section{The lower bound of $|M^G|$ and examples}\label{s4}

 Let $(\Phi, M^n)$ be a $G$-action in $\mathcal{A}_n^k$, and let
$(\Gamma_{(\Phi, M)}, \alpha)$ be a colored graph associated to
$(\Phi, M^n)$. Since $(\Phi, M^n)$ is assumed to be nonbounding, by
Theorem~\ref{ce} there must be some edges in $\Gamma_{(\Phi, M)}$,
each of which has different colorings at its two endpoints.

\begin{lem} \label{l2}
Let $e\in E_{\Gamma_{(\Phi, M)}}$ be an edge with two endpoints $p$
and $q$ such that $\alpha(E_p)\not=\alpha(E_q)$. Then the number
$\vert E_e\vert$ is at most $n-k+1$, where $E_e$ denotes the set of all edges joining
 two endpoints $p$ and $q$.
\end{lem}
\begin{proof} Without loss of generality, one may assume that
$$\alpha(E_p)=\{\beta_1,\beta_2,...,\beta_t,\gamma_1,...,\gamma_{n-t}\}$$
and $$ \alpha(E_q)=\{\beta_1,\beta_2,...,\beta_t,\gamma'_1,
\gamma'_2,...,\gamma'_{n-t}\}$$ where $\beta_1,\beta_2,...,\beta_t$
are the  representations for the $t$ edges joining $p$ to $q$.

\vskip .2cm
 Suppose that $t\geq n-k+2$.
Choose $\gamma_1, \gamma_2,...,\gamma_s$ to be a maximal linearly
independent set from $\{\gamma_1,...,\gamma_{n-t}\}$ (labeling them
as the first $s$ elements), where $s\leq n-t\leq k-2$. Since
$\alpha(E_p)$ spans $\Hom(G,{\Bbb Z}_2)$ by (P1) of
Section~\ref{s2}, one can choose at least two of the $\beta$'s, say
$\beta_1$ and $\beta_2$, so that the set $\{\beta_1,
\beta_2,\gamma_1,...,\gamma_s\}$ is linearly independent.

\vskip .2cm

Since $\alpha(E_p)\not=\alpha(E_q)$, there must be some nonzero
element $\delta$ in $\Hom(G,{\Bbb Z}_2)$ which occurs more times in
$\alpha(E_p)$ than in  $\alpha(E_q)$. Since each $\beta_i$ occurs
the same number of times in both $\alpha(E_p)$ and $\alpha(E_q)$,
one has that $\delta\in\{\gamma_1,...,\gamma_{n-t}\}$. Since  the
number of times $\delta$ and $\delta+\beta_1$ occurring in
$\alpha(E_p)$ is the same as the number of times $\delta$ and
$\delta+\beta_1$ occurring in  $\alpha(E_q)$ by (P2) of
Section~\ref{s2}, $\delta+\beta_1$ must occur more times in
$\alpha(E_q)$ than in $\alpha(E_p)$. Then $\delta+\beta_1+\beta_2$
must occur more times in $\alpha(E_p)$ than in $\alpha(E_q)$. Thus
one has $\delta+\beta_1+\beta_2\in\{\gamma_1,...,\gamma_{n-t}\}$.
Then $\beta_1+\beta_2$ is in the span of
$\{\gamma_1,\gamma_2,...,\gamma_s\}$ contradicting the linear
independence of $\{\beta_1,\beta_2,\gamma_1,...,\gamma_s\}$.
\end{proof}

Now let us give the proof of Theorem~\ref{thm2}.

\vskip .2cm

 \noindent {\em Proof of Theorem~\ref{thm2}.} By Lemma~\ref{lm1}, we can replace
$(\Phi, M)$ by $(\Psi, N)$ with $\mathcal{N}_{(\Psi,
N)}=\mathcal{P}_{(\Phi, M^n)}$ prime tangent representation set.
This decreases the number of fixed points. Let $p$ be one fixed
point of $N^G$. Then $n$ edges at $p$ of the graph
$\Gamma_{(\Psi,N)}$ would connect $p$ to other fixed points and at
most $n-k+1$ edges can connect $p$ to the same point by
Lemma~\ref{l2}. Thus there must be at least ${n\over {n-k+1}}$ fixed
points not equal to $p$, and so the number of fixed points is at
least $1+\lceil{n\over{n-k+1}}\rceil$. $\hfill\square$

\vskip .2cm

Next let us give two examples to show that the bound established in
Theorem~\ref{thm2} is attainable in some special cases.

\begin{example}\label{5-1}
For $n=k$, one has that ${n\over {n-k+1}}=n=k$,  so this says that
the action has at least $k+1$ fixed points. Consider the standard
linear  $({\Bbb Z}_2)^k$-action $(\phi_k, {\Bbb R}P^k)$. This action
fixes precisely $k+1$ fixed points.
 \end{example}

\begin{example} \label{5-2}
For $k\geq 2$,  applying the $\Omega$-operation $k-2$ times to
$(\phi_2,$ ${\Bbb R}P^2)$ gives a  $({\Bbb Z}_2)^{k}$-action
 $\Omega^{k-2}(\phi_2,$ ${\Bbb R}P^2)$, which
has 3 fixed points and  dimension  $2^{k-1}$. Then
$1+\lceil{{2^{k-1}}\over{2^{k-1}-k+1}}\rceil=3$
 for all $k\geq 2$.
 \end{example}

 The following example illustrates that the bound
established in Theorem~\ref{thm2} can be much smaller than the
actual number of fixed points.

 \vskip .2cm

For $k=2$, Conner and Floyd \cite{cf} found the equivariant
cobordism classes of $({\Bbb Z}_2)^2$-actions with finite fixed set.
Begin with the standard linear $({\Bbb Z}_2)^2$-action
$(\phi_2,{\Bbb R}P^2)$ having three fixed points,  Conner and Floyd
wrote its fixed data as a polynomial
$\rho_1\rho_2+\rho_2\rho_3+\rho_3\rho_1$ in ${\Bbb Z}_2[\rho_1,
\rho_2, \rho_3]$. They then proved that every nonbounding $({\Bbb
Z}_2)^2$-action $(\Phi, M^n)$ having a finite fixed set is
equivariantly cobordant to the $m$-multi-diagonal $({\Bbb
Z}_2)^2$-action $\Delta^m(\phi_2, {\Bbb R}P^2)$ where $n=2m$. Up to
equivariant cobordism, the fixed data of this action can be written
as $(\rho_1\rho_2+\rho_2\rho_3+\rho_3\rho_1)^m.$ The minimum number
of fixed points  is then the number of monomials
$\rho_1^i\rho_2^j\rho_3^h$ in
$(\rho_1\rho_2+\rho_2\rho_3+\rho_3\rho_1)^m$ which have nonzero
coefficient modulo 2. To find this number, let
$m=2^{p_1}+\cdots+2^{p_r}$ be the 2-adic expansion of $m$. Then
\begin{eqnarray*}& & (\rho_1\rho_2+\rho_2\rho_3+\rho_3\rho_1)^m\\
&\equiv&
(\rho^{2^{p_1}}_1\rho^{2^{p_1}}_2+\rho^{2^{p_1}}_2\rho^{2^{p_1}}_3+\rho^{2^{p_1}}_3\rho^{2^{p_1}}_1)\cdots
(\rho^{2^{p_r}}_1\rho^{2^{p_r}}_2+\rho^{2^{p_r}}_2\rho^{2^{p_r}}_3+\rho^{2^{p_r}}_3\rho^{2^{p_r}}_1)\mod
2.
\end{eqnarray*}
This product has $3^r$ monomials and the monomials are all
distinct---for $\rho_1^i\rho_2^j\rho_3^h$ the 2-adic expansions of
$i, j$, and $h$ determine which factors were taken. This gives

\begin{lem}\label{k=2}
Let $(\Phi, M^n)$ be a nonbounding $({\Bbb Z}_2)^2$-action in
$\mathcal{A}_n^2$. Then $n$ is even and $(\Phi, M^n)$ has at least
$3^r$ fixed points where $r$ is the number of terms in the 2-adic
expansion of $n$.
\end{lem}

Of course, $1+\lceil{n\over{n-1}}\rceil=3$, so the bound in
Theorem~\ref{thm2} is now precise exactly for $n=2^s$.

\begin{rem}
The reason why the bound established in Theorem~\ref{thm2} can be
much smaller than the actual number of fixed points is because the
proof of Theorem~\ref{thm2} is only a local analysis for the colored
graph.
\end{rem}

\section{Actions with three fixed points}\label{s5}

\subsection{The existence and the equivariant cobordism classification of actions in $\mathcal{A}_n^k(3)$} Suppose that $\mathcal{A}_n^k(3)$ is nonempty. Take a $G$-action $(\Phi, M^n)$ in $\mathcal{A}_n^k(3)$, without
the loss of generality,
 assume
that $(\Phi,M^n)$ has exactly three fixed points $p,q,r$. Let
$(\Gamma_{(\Phi, M^n)},\alpha)$ be a colored graph of $(\Phi, M^n)$.

\begin{lem}\label{lm3}
$\Gamma_{(\Phi, M^n)}$ is uniquely determined, and
$\mathcal{N}_{(\Phi,M^n)}$ $=\{\alpha(E_p)$, $\alpha(E_q)$,
$\alpha(E_r)\}$ is prime.
\end{lem}
\begin{proof}
 Let $p, q$ be connected by $a$
 edges; $p,r$ by $b$ edges; and $q,r$ by $c$ edges. Then
 $$\begin{cases}a+b=n\\
  b+c=n\\
   a+c=n
   \end{cases}$$ so $n$ must be even and $a=b=c={n\over 2}$. This
   means that $(\Phi, M^n)$ determines a unique graph $\Gamma_{(\Phi,
   M^n)}$.

   \vskip .2cm
    It is obvious that $\mathcal{N}_{(\Phi,M^n)}=\{\alpha(E_p),
\alpha(E_q),\alpha(E_r)\}$ is  prime since any $({\Bbb
Z}_2)^k$-action can not fix a single isolated  point.
\end{proof}
By Lemma~\ref{lm3}
   one can write $$\alpha(E_p)=\hat{\beta}\cup\hat{\gamma},\ \ \
\alpha(E_q)=\hat{\beta}\cup\hat{\delta},\ \ \
\alpha(E_r)=\hat{\delta}\cup\hat{\gamma}$$ where
$\hat{\beta}$, %=\{\beta_1,...,\beta_{n/2}\}$,
$\hat{\gamma}$, %=\{\gamma_1,...,\gamma_{n/2}\}$,
$\hat{\delta}$ %=\{\delta_1,...,\delta_{n/2}\}$
are three multisets formed by elements of $\Hom(({\Bbb Z}_2)^k,{\Bbb
Z}_2)$ with $|\hat{\beta}|=|\hat{\gamma}|=|\hat{\delta}|=n/2$.

\begin{lem} \label{lm4}  For any $\beta\in\hat{\beta}, \gamma\in\hat{\gamma}$, and
$\delta\in\hat{\delta}$, $$
 \beta+\gamma\in
\hat{\delta}, \ \ \gamma+\delta\in \hat{\beta}, \ \ \beta+\delta\in
\hat{\gamma}. $$
\end{lem}

\begin{proof}
Consider an irreducible nontrivial representation $\rho$ and a fixed
component $C$ of $\ker\rho$ acting on $M^n$. Then $({\Bbb
Z}_2)^k/\ker\rho\cong {\Bbb Z}_2$ fixes an even number of points of
$C$. If $C$ has a fixed point, then it has exactly two fixed points,
so $C$ exactly determines a unique $\dim C$-valent regular graph
$\Gamma_{\rho,C}$ with those two fixed points as its vertices. This
means that $\hat{\beta}, \hat{\gamma}, \hat{\delta}$ are all
disjoint.  Further, the lemma follows from the property (P2) in
Section 2.
\end{proof}

\begin{lem}\label{lm5}
There exists an integer $m\geq1$ and a basis
$\{\beta_1,...,\beta_{k-1},\gamma\}$ of $\Hom$ $(({\Bbb Z}_2)^k$,
${\Bbb Z}_2)$ such that
\begin{enumerate}
\item[(i)]
 $n=m\cdot 2^{k-1}$;
\item[(ii)] $\hat{\beta}$ is a multiset consisting of all odd sums with multiplicity $m$ formed by
$\beta_1,...,\beta_{k-1}$;
\item[(iii)] $\hat{\gamma}=\{\gamma+\beta_1+\beta\ \vert
\beta\in\hat{\beta}\}$;
\item[(iv)]$\hat{\delta}=\{\gamma+\beta\ \vert \beta\in\hat{\beta}\}$
\end{enumerate}
where an odd sum formed by $\beta_1,...,\beta_{k-1}$ means the sum
of an odd number of elements in $\beta_1,...,\beta_{k-1}$.
\end{lem}

\begin{proof}
Choose one $\gamma\in\hat{\gamma}$. Then
$\hat{\delta}=\{\gamma+\beta|\beta\in \hat{\beta}\}$ by
Lemma~\ref{lm4}. Choosing any $\beta'\in \hat{\beta}$, one has
$\hat{\gamma}=\{\gamma+\beta+\beta'|\beta\in \hat{\beta}\}$ by
Lemma~\ref{lm4} again. Similarly,
$\hat{\delta}=\{\delta+\beta+\beta'|\beta\in \hat{\beta}\}$ for any
$\delta$ in $\hat{\delta}$.
 Then
$$\hat{\beta}\cup\hat{\gamma}=\{\beta|\beta\in \hat{\beta}\}\cup \{
\gamma+\beta+\beta'|\beta\in \hat{\beta}\}$$ is spanned by elements
$\beta$ of $\hat{\beta}$ and $\gamma=\gamma+\beta'+\beta'$. Since
$\hat{\beta}\cup\hat{\gamma}$ spans $\Hom(({\Bbb Z}_2)^k,{\Bbb
Z}_2)$ by (P1) in Section~\ref{s2}, one has that at least $k-1$
elements of $\hat{\beta}$ must be linearly independent, so $n/2\geq
k-1$.

\vskip .2cm

\noindent {\bf Claim I.} {\em Every sum of an odd number of elements
of $\hat{\beta}$ is again an element of $\hat{\beta}$ and no sum of
an even number of elements of $\hat{\beta}$ is again in
$\hat{\beta}$. The elements $\gamma+$ even sum of $\beta$'s all
belong to $\hat{\gamma}$ and the elements $\gamma+$ odd sum of
$\beta$'s all belong to $\hat{\delta}$.}

\vskip .2cm

In fact, one has from the above argument that
$\hat{\delta}=\{\gamma+\beta|\beta\in \hat{\beta}\}$ and
$\hat{\delta}=\{\delta+\beta+\beta'|\beta\in \hat{\beta}\}$ for any
$\delta$ in $\hat{\delta}$ and any $\beta'\in \hat{\beta}$. Thus,
for any $\beta''\in\hat{\beta}$,  taking $\delta=\gamma+\beta''$
gives
$$\hat{\beta}=\{\beta|\beta\in \hat{\beta}\}=\{\beta+\beta'+\beta''|\beta\in
\hat{\beta}\}.
$$
Replacing $\beta$ by $\beta+\beta'+\beta''$ and repeating this
procedure, one has that every sum of an odd number of elements of
$\hat{\beta}$ is again an element of $\hat{\beta}$. Furthermore, by
Lemma~\ref{lm4}, one has that  the elements $\gamma+$ odd sum of
$\beta$'s all belong to $\hat{\delta}$, the elements $\gamma+$ even
sum of $\beta$'s all belong to $\hat{\gamma}$, and any sum of an
even number of elements of $\hat{\beta}$ is not in $\hat{\beta}$.

\vskip .2cm

\noindent {\bf Claim II.} {\em There are exactly $k-1$ elements of
the $\beta$'s which are linearly independent.}

\vskip .2cm

Actually, one knows from Claim I that an odd sum of the $\beta$'s is
again in $\hat{\beta}$ and if $\gamma\in\hat{\gamma}$ is an even sum
of the $\beta$'s, then any $\delta$ in $\hat{\delta}$ is an odd sum
of the $\beta$'s, so one has that
$\hat{\delta}\cap\hat{\beta}\not=\emptyset$. But this is impossible.
Thus, $\gamma\in\hat{\gamma}$ cannot be a linear combination of
elements of $\hat{\beta}$, so exactly $k-1$ elements of the
$\beta$'s are linearly independent, say $\beta_1,...,\beta_{k-1}$.
In particular,  $\beta_1,...,\beta_{k-1}, \gamma$ form a basis of
$\Hom(({\Bbb Z}_2)^k,{\Bbb Z}_2)$.

\vskip .2cm

Now let $\mathcal{V}$ be the vector space spanned by the elements of
$\hat{\beta}$. Then, by Claims I and II, $\mathcal{V}$ is a
$(k-1)$-dimensional vector space over ${\Bbb Z}_2$ and
$\beta_1,...,\beta_{k-1}$  form a basis of $\mathcal{V}$. The set of
elements of $\mathcal{V}$ belonging to $\hat{\beta}$ must then be
the set of the sums of an odd number of the elements
$\beta_1,...,\beta_{k-1}$ by Claim I. Thus one has

\vskip .2cm

\noindent {\bf Claim III.} {\em $\hat{\beta}$ contains $2^{k-2}$
different elements.}

\vskip .2cm

  Since any two odd sums of the elements of
$\hat{\beta}$ differ by an even sum, and adding an even sum permutes
the elements of $\hat{\beta}$, the elements of $\hat{\beta}$ must
occur with the same multiplicity. Since $\hat{\beta}$ contains
$2^{k-2}$ different elements, one has that $n/2=m\cdot 2^{k-2}$ so
$n=m\cdot 2^{k-1}$, where $m$ is the common multiplicity. This gives

\vskip .2cm

\noindent {\bf Claim IV.} {\em All elements of $\hat{\beta}$ occur
with the same multiplicity $m$ so $n=m\cdot 2^{k-1}$.}

\vskip .2cm

Combining the above arguments completes the proof of the lemma.
\end{proof}

We see from Lemma~\ref{lm5} that  the tangent representation set
$\mathcal{N}_{(\Phi,M^n)}$ $=\{\alpha(E_p)$, $\alpha(E_q),
\alpha(E_r)\}$ is uniquely determined by some basis of $\Hom(({\Bbb
Z}_2)^k,{\Bbb Z}_2)$.  On the other hand, any basis of $\Hom(({\Bbb
Z}_2)^k,{\Bbb Z}_2)$ can be translated into a given basis by an
automorphism of $({\Bbb Z}_2)^k$. Thus, by Theorem~\ref{ce} and
Remark~\ref{automorphism}, it follows that

\begin{prop} \label{p1} Let   $(\Psi, N^n)$ be another  $({\Bbb Z}_2)^k$-action in $\mathcal{A}_n^k(3)$.
Then there is one $\sigma\in \text{\rm GL}(k,{\Bbb Z}_2)$ such that
$(\Phi, M^n)$ is equivariantly cobordant to $(\sigma\Psi, N^n)$.
\end{prop}

Next, let us  look at the case $m=1$.

\vskip .2cm

When $m=1$, one has that $n=2^{k-1}$. In this case,
Example~\ref{5-2} in Section~\ref{s4} provides an example
$\Omega^{k-2}(\phi_2,{\Bbb R}P^2)$ which fixes exactly three
isolated points,  so by Proposition~\ref{p1} any  $({\Bbb
Z}_2)^k$-action  in $\mathcal{A}^k_{2^{k-1}}(3)$ is equivariantly
cobordant to the $({\Bbb Z}_2)^k$-action obtained by applying an
automorphism of $({\Bbb Z}_2)^k$ to $\Omega^{k-2}(\phi_2,{\Bbb
R}P^2)$ to switch representations around. This gives

\begin{cor} \label{cor2}
$\mathcal{A}_{2^{k-1}}^k(3)$ is nonempty. Furthermore, each of
$\mathcal{A}_{2^{k-1}}^k(3)$ is equivariantly cobordant to one of
$\sigma\Omega^{k-2}(\phi_2, {\Bbb R}P^2)$, $\sigma\in \text{\rm
GL}(k,{\Bbb Z}_2)$.
\end{cor}

\begin{rem}\label{r7}
For $k=2$, it is easy to see that up to equivariant cobordism, there
is a unique $({\Bbb Z}_2)^2$-action in $\mathcal{A}_2^2(3)$, which
is exactly $\Omega^0(\phi_2,{\Bbb R}P^2)=(\phi_2,{\Bbb R}P^2)$.
However, for $k>2$, up to equivariant  cobordism there may be  more
$({\Bbb Z}_2)^k$-actions in $\mathcal{A}_{2^{k-1}}^k(3)$. For
example, for $k=3$, up to equivariant  cobordism there are seven
different $({\Bbb Z}_2)^3$-actions in $\mathcal{A}_4^3(3)$, whose
tangent representation sets are listed as follows:
$$\begin{tabular}{|l|l|}
 \hline $i$ & \ \ \ \ \ \ \ \ \ \ \ \ \ \ \ \ \ \ \ \ \ \ \ \ \ \ \ \ \ \ tangent representation set $\mathcal{N}_i$
 \\
\hline 1 &  $\{\rho_1,\rho_2,\rho_3,\rho_1+\rho_2+\rho_3\},
\{\rho_2,\rho_3,\rho_1+\rho_2,\rho_1+\rho_3\}$,\\
&$\{\rho_1,\rho_1+\rho_2,\rho_1+\rho_3,\rho_1+\rho_2+\rho_3\}$\\
\hline 2 & $\{\rho_1,\rho_2,\rho_3,\rho_1+\rho_2+\rho_3\},\{\rho_1,
\rho_3,\rho_1+\rho_2,\rho_2+\rho_3\},$\\
&$\{\rho_2,\rho_1+\rho_2,\rho_2+\rho_3,\rho_1+\rho_2+\rho_3\}$\\
\hline 3 & $\{\rho_1,\rho_2,\rho_3,\rho_1+\rho_2+\rho_3\},\{\rho_1,
\rho_2,\rho_1+\rho_3,\rho_2+\rho_3\},$\\
&$\{\rho_3,\rho_1+\rho_3,\rho_2+\rho_3,\rho_1+\rho_2+\rho_3\}$\\
\hline 4 &
$\{\rho_1,\rho_1+\rho_2,\rho_1+\rho_3,\rho_1+\rho_2+\rho_3\},
\{\rho_1,\rho_2,\rho_1+\rho_3,\rho_2+\rho_3\}$,\\
&$\{\rho_2,\rho_1+\rho_2,\rho_2+\rho_3,\rho_1+\rho_2+\rho_3\}$\\
\hline 5 &
$\{\rho_2,\rho_1+\rho_2,\rho_2+\rho_3,\rho_1+\rho_2+\rho_3\},
\{\rho_2,\rho_3,\rho_1+\rho_2,\rho_1+\rho_3\}$,\\
& $\{\rho_3,\rho_1+\rho_3,\rho_2+\rho_3,\rho_1+\rho_2+\rho_3\}$\\
\hline 6 &
$\{\rho_3,\rho_1+\rho_3,\rho_2+\rho_3,\rho_1+\rho_2+\rho_3\},
\{\rho_1,\rho_3,\rho_1+\rho_2,\rho_2+\rho_3\}$,\\
&$\{\rho_1,\rho_1+\rho_2,\rho_1+\rho_3,\rho_1+\rho_2+\rho_3\}$\\
\hline 7 & $\{\rho_1,\rho_2,\rho_1+\rho_3,\rho_2+\rho_3\},
\{\rho_1,\rho_3,\rho_1+\rho_2,\rho_2+\rho_3\},
\{\rho_2,\rho_3,\rho_1+\rho_2,\rho_1+\rho_3\}$\\
\hline
\end{tabular}$$
 \noindent where  $\{\rho_1,\rho_2,\rho_3\}$ is the standard basis
of $\Hom(({\Bbb Z}_2)^3,{\Bbb Z}_2)$ as stated in Section~\ref{s2}.
\end{rem}

Now let us look at the general case $m\geq 1$.

\begin{lem} \label{lm6}
Let $m\geq 1$. Then $m$ is a power of 2.
\end{lem}
\begin{proof} By Lemma~\ref{lm5}, $\hat{\beta}$ contains $2^{k-2}$ different
elements (say $\bar{\beta}_1,...,\bar{\beta}_{2^{k-2}}$) with same
multiplicity $m$, which consist of all odd sums  formed by
$\beta_1,...,\beta_{k-1}$. Without loss of generality, assume that
$\bar{\beta}_1=\beta_1$. Applying Theorem~\ref{dks}, we use the symmetric function $f=1$ to deduce that
\begin{eqnarray*}
&&{1\over{(\prod_{i=1}^{2^{k-2}}
\bar{\beta}_i)^m\big(\prod_{i=1}^{2^{k-2}}(\gamma+\beta_1+\bar{\beta}_i)\big)^m}}
+{1\over{(\prod_{i=1}^{2^{k-2}} \bar{\beta}_i)^m\big(\prod_{i=1}^{2^{k-2}}(\gamma+\bar{\beta}_i)\big)^m}}\\
&&+ {1\over{\big(\prod_{i=1}^{2^{k-2}}
(\gamma+\bar{\beta}_i)\big)^m\big(\prod_{i=1}^{2^{k-2}}(\gamma+\beta_1+\bar{\beta}_i)\big)^m}}
\end{eqnarray*}
must belong to ${\Bbb Z}_2[\rho_1,...,\rho_{k}]$. Taking the common
denominator, one has that
$${{(\prod_{i=1}^{2^{k-2}}\bar{\beta}_i)^m+\big(\prod_{i=1}^{2^{k-2}}(\gamma+\bar{\beta}_i)\big)^m+
\big(\prod_{i=1}^{2^{k-2}}(\gamma+\beta_1+\bar{\beta}_i)\big)^m}\over{
\big[(\prod_{i=1}^{2^{k-2}}\bar{\beta}_i)\big(\prod_{i=1}^{2^{k-2}}(\gamma+\bar{\beta}_i)\big)\big(\prod_{i=1}^{2^{k-2}}
(\gamma+\beta_1+\bar{\beta}_i)\big)\big]^m}}$$ also belongs to
${\Bbb Z}_2[\rho_1,...,\rho_{k}]$. Since the numerator has smaller
degree than the denominator, this is only possible if the numerator
is zero. Thus, one has
\begin{eqnarray*}(\prod_{i=1}^{2^{k-2}}\bar{\beta}_i)^m+\big(\prod_{i=1}^{2^{k-2}}(\gamma+\bar{\beta}_i)\big)^m+
\big(\prod_{i=1}^{2^{k-2}}(\gamma+\beta_1+\bar{\beta}_i)\big)^m=0.
\end{eqnarray*}
By Lemma~\ref{lm5} and Corollary~\ref{cor2}, when $m=1$, the above
expression still holds, i.e.,
\begin{eqnarray*}\prod_{i=1}^{2^{k-2}}\bar{\beta}_i+\prod_{i=1}^{2^{k-2}}(\gamma+\bar{\beta}_i)+
\prod_{i=1}^{2^{k-2}}(\gamma+\beta_1+\bar{\beta}_i)=0.
\end{eqnarray*}
Furthermore,  one has
\begin{eqnarray*}
 \prod_{i=1}^{2^{k-2}}(\gamma+\bar{\beta}_i) =
\prod_{i=1}^{2^{k-2}}(\gamma+\beta_1+\bar{\beta}_i)+\prod_{i=1}^{2^{k-2}}\bar{\beta}_i
=
\gamma\cdot\prod_{i\not=1}(\gamma+\beta_1+\bar{\beta}_i)+\prod_{i=1}^{2^{k-2}}\bar{\beta}_i
\end{eqnarray*}
and then
\begin{eqnarray*}
& &
(\prod_{i=1}^{2^{k-2}}\bar{\beta}_i)^m+(\prod_{i=1}^{2^{k-2}}(\gamma+\bar{\beta}_i))^m+\big(\prod_{i=1}^{2^{k-2}}(\gamma+\beta_1+\bar{\beta}_i)\big)^m\\
&=&(\prod_{i=1}^{2^{k-2}}\bar{\beta}_i)^m+
\Big\{\gamma\cdot\prod_{i\not=1}(\gamma+\beta_1+\bar{\beta}_i)+\prod_{i=1}^{2^{k-2}}\bar{\beta}_i\Big\}^m
+\Big\{\gamma\cdot\prod_{i\not=1}(\gamma+\beta_1+\bar{\beta}_i)\Big\}^m\\
&=& \sum_{0<j<m}{m\choose
j}\Big\{\gamma\cdot\prod_{i\not=1}(\gamma+\beta_1+\bar{\beta}_i)\Big\}^j(\prod_{i=1}^{2^{k-2}}\bar{\beta}_i)^{m-j}.
\end{eqnarray*}
If $m$ is not a power of 2, then there is a largest $j$, $0<j<m$,
with ${m\choose j}\not=0$ and this sum has a nonzero coefficient for
$$(\gamma^{2^{k-2}})^j(\prod_{i=1}^{2^{k-2}}\bar{\beta}_i)^{m-j}.$$
This is a contradiction.  Thus $m$ must be a power of 2.
\end{proof}

Now, let us complete the proof of Theorems~\ref{thm3}.

\vskip .2cm

\noindent {\em Proof of Theorem~\ref{thm3}.}   If
$\mathcal{A}_n^k(3)$ is nonempty, then we may choose a $G$-action
$(\Phi, M^n)$ in $\mathcal{A}_n^k(3)$. Without loss of generality, assume that
$(\Phi, M^n)$ exactly fixes three isolated points. By
Lemmas~\ref{lm5} and \ref{lm6} one has that $k\geq 2$ and
$n=2^{\ell}$ for some $\ell\geq k-1$. Conversely, let $n=2^{\ell}$
with $\ell\geq k-1\geq 1$. Then the diagonal action on
$2^{\ell-k+1}$ copies of $\Omega^{k-2}(\phi_2,{\Bbb R}P^2)$ gives a
$({\Bbb Z}_2)^k$-action
$\Delta^{2^{\ell-k+1}}\Omega^{k-2}(\phi_2,{\Bbb R}P^2)$, which has
dimension $2^\ell$ and fixes $3^{ 2^{\ell-k+1}}$ isolated points.
However, by Lemma~\ref{dia} this action is actually cobordant to the
action with exactly three fixed points. Thus, $\mathcal{A}_n^k(3)$
is nonempty. This completes the proof of Theorem~\ref{thm3}(a).

\vskip .2cm Now if  $\mathcal{A}_n^k(3)$ is nonempty, then
$n=2^\ell$ with $\ell\geq k-1\geq 1$. Furthermore, by
Proposition~\ref{p1}, each action of $\mathcal{A}_n^k(3)$ is
equivariantly cobordant to one of
$$\sigma\Delta^{2^{\ell-k+1}}\Omega^{k-2}(\phi_2,{\Bbb R}P^2),\ \
\sigma\in \text{\rm GL}(k,{\Bbb Z}_2).$$ $\hfill\Box$

\subsection{The characterization of the colored graphs--Proof of Theorem~\ref{thm5}} \label{thm 1.4} Let $(\Gamma, \alpha)$ be
an abstract 1-skeleton of type $(n, k)$ with exactly
three vertices $p,q,r$.

\vskip .2cm

If $(\Gamma, \alpha)$ is a colored graph of some action $(\Phi,
M^n)$ in $\mathcal{A}_n^k(3)$, then by Lemmas~\ref{lm5} and
\ref{lm6} the necessity of Theorem~\ref{thm5} holds.

\vskip .2cm

Conversely, suppose that $(\Gamma, \alpha)$ satisfies the conditions
(a) and (b) of Theorem~\ref{thm5}. Then it is easy to see that
$\Gamma$ is uniquely determined. Moreover, to prove that  $(\Gamma,
\alpha)$ is a colored graph of some action $(\Phi, M^n)$ in
$\mathcal{A}_n^k(3)$, it suffices to show that $\{\alpha(E_p),
\alpha(E_q), \alpha(E_r)\}$ is the fixed data of some action $(\Phi,
M^n)$ in $\mathcal{A}_n^k(3)$. By the construction of $\hat{\beta},
\hat{\gamma}, \hat{\delta}$ in Theorem~\ref{thm5}(b), we see that
$\hat{\beta}$ exactly contains $2^{k-2}$ different elements (which
consist of all odd sums formed by $\beta_1,...,\beta_{k-1}$), and so
are $\hat{\gamma}$ and $\hat{\delta}$. By $\hat{\beta'}$ (resp.
$\hat{\gamma'}, \hat{\delta'}$) we denote the set formed by those
$2^{k-2}$ different elements in $\hat{\beta}$ (resp. $\hat{\gamma},
\hat{\delta}$). Then we have that
$\hat{\gamma'}=\{\gamma+\beta_1+\beta \ \vert\
\beta\in\hat{\beta'}\}$ and $\hat{\delta'}=\{\gamma+\beta\ \vert\
\beta\in\hat{\beta'}\}$. Obviously, $\hat{\beta'}, \hat{\gamma'},
\hat{\delta'}$ are uniquely determined by the basis $\{\beta_1, ...,
\beta_{k-1}, \gamma\}$.

\vskip .2cm {\bf Claim.} {\em $\{\hat{\beta'}\cup\hat{\gamma'},
\hat{\beta'}\cup\hat{\delta'}, \hat{\delta'}\cup\hat{\gamma'}\}$ is
the fixed data of some action in $\mathcal{A}^k_{2^{k-1}}(3)$}.

\vskip .2cm

We proceed by induction on $k$.  When $k=2$, we have that
$\hat{\beta'}=\{\beta_1\}$, $\hat{\gamma'}=\{\gamma\}$, and
$\hat{\delta'}=\{\gamma+\beta_1\}$. Since $\Hom(({\Bbb Z}_2)^2,{\Bbb
Z}_2)$ contains only three nontrivial elements, without the loss of
generality one may assume that $\beta_1=\rho_1$ and $\gamma=\rho_2$.
Then  $\hat{\beta'}\cup\hat{\gamma'}=\{\rho_1,\rho_2\}$,
$\hat{\beta'}\cup\hat{\delta'}=\{\rho_1,\rho_1+\rho_2\}$, and
$\hat{\delta'}\cup\hat{\gamma'}=\{\rho_1+\rho_2\,\rho_2\}$ are
exactly  the fixed data of the standard linear $({\Bbb
Z}_2)^2$-action on $(\phi_2,{\Bbb R}P^2)$. When $k=l\geq 2$, suppose
inductively that $\hat{\beta'}\cup\hat{\gamma'},
\hat{\beta'}\cup\hat{\delta'}, \hat{\delta'}\cup\hat{\gamma'}$
 are the fixed data of some
$({\Bbb Z}_2)^l$-action  in $\mathcal{A}^l_{2^{l-1}}(3)$.

\vskip .2cm

When $k=l+1$,  let $\hat{\beta'}_1$ denote the set of all odd sums
formed by  $\beta_1,...,\beta_{l-1}$. Then
$\hat{\beta'}_1\subset\hat{\beta'}$ and $\hat{\beta'}_1$ contains
$2^{l-2}$ different elements.  Let $\hat{\beta'}_2$ denote the set
formed by all elements $\beta_1+\beta_{l}+\beta, \beta\in
\hat{\beta'}_1$. Then,  $\hat{\beta'}_2\subset\hat{\beta'}$ and
$\hat{\beta'}_2$ contains $2^{l-2}$ different elements, too. Since
$\beta_1,...,\beta_{l}$ are linearly independent, one has that the
intersection of $\hat{\beta'}_1$ and $\hat{\beta'}_2$ is empty and
$\hat{\beta'}=\hat{\beta'}_1\cup \hat{\beta'}_2$.  Then one has
$$\hat{\gamma'}=\hat{\gamma'}_1\cup\hat{\gamma'}_2\text{ with } \hat{\gamma'}_1\cap\hat{\gamma'}_2=\emptyset$$ where
$\hat{\gamma'}_1=\{\gamma+\beta_1+\beta\vert \beta\in
\hat{\beta'}_1\}$ and $\hat{\gamma'}_2=\{\gamma+\beta_1+\beta\vert
\beta\in \hat{\beta'}_2\}$, and
$$\hat{\delta'}=\hat{\delta'}_1\cup\hat{\delta'}_2\text{ with }\hat{\delta'}_1\cap\hat{\delta'}_2=\emptyset$$
where $\hat{\delta'}_1=\{\gamma+\beta\vert \beta\in
\hat{\beta'}_1\}$ and $\hat{\delta'}_2=\{\gamma+\beta\vert \beta\in
\hat{\beta'}_2\}$. Now let us look at $\hat{\beta'}_1,
\hat{\gamma'}_1$, and $\hat{\delta'}_1$.  Clearly $\hat{\beta'}_1,
\hat{\gamma'}_1$, and $\hat{\delta'}_1$ are exactly formed by
$\beta_1,...,\beta_{l-1},\gamma$. One sees that
$\beta_1,...,\beta_{l-1},\gamma$ span a $l$-dimensional subspace of
$\Hom(({\Bbb Z}_2)^{l+1},{\Bbb Z}_2)$ which is isomorphic to
$\Hom(({\Bbb Z}_2)^{l},{\Bbb Z}_2)$. Now, regarding
$\{\beta_1,...,\beta_{l-1},\gamma\}$ as a basis of $\Hom(({\Bbb
Z}_2)^{l},{\Bbb Z}_2)$,  one has by induction that
$$\hat{\beta'}_1\cup\hat{\gamma'}_1,\ \
\hat{\beta'}_1\cup\hat{\delta'}_1,\ \
\hat{\delta'}_1\cup\hat{\gamma'}_1$$  are the fixed data of some
$({\Bbb Z}_2)^l$-action, denoted by $(\Psi, N^{2^{l-1}})$. Then by
applying $\Omega$-operation to $(\Psi, N^{2^{l-1}})$, as in the
proof of Lemma~\ref{twist} (see also \cite[Lemma 4.1]{l2}),
 the fixed
data of $\Omega(\Psi,N^{2^{l-1}})$ exactly consists of
$$\begin{cases}
\hat{\beta'}\cup\hat{\gamma'}=\hat{\beta'}_1\cup\hat{\beta'}_2\cup
\hat{\gamma'}_1\cup\hat{\gamma'}_2 \\
\hat{\beta'}\cup\hat{\delta'}=\hat{\beta'}_1\cup\hat{\beta'}_2\cup
\hat{\delta'}_1\cup\hat{\delta'}_2\\
\hat{\delta'}\cup\hat{\gamma'}=\hat{\delta'}_1\cup\hat{\delta'}_2\cup
\hat{\gamma'}_1\cup\hat{\gamma'}_2.
\end{cases}$$ This completes the
induction and the proof of the claim.

\vskip .2cm

 Now by Corollary~\ref{cor2}, there is an automorphism $\sigma\in\text{\rm
GL}(k, {\Bbb Z}_2)$ such that $\{\hat{\beta'}\cup\hat{\gamma'},
\hat{\beta'}\cup\hat{\delta'}, \hat{\delta'}\cup\hat{\gamma'}\}$ is
the fixed data of $\sigma\Omega^{k-2}(\phi_2, {\Bbb R}P^2)$.
Moreover, applying $\Delta$-operation $2^{\ell-k+1}$ times to
$\sigma\Omega^{k-2}(\phi_2, {\Bbb R}P^2)$ gives an action
$\Delta^{2^{\ell-k+1}}[\sigma\Omega^{k-2}(\phi_2, {\Bbb R}P^2)]$ in
$\mathcal{A}_n^k(3)$. By Lemma~\ref{dia},
$\Delta^{2^{\ell-k+1}}[\sigma\Omega^{k-2}(\phi_2, {\Bbb R}P^2)]$ is
equivariantly cobordant to a $({\Bbb Z}_2)^k$-action such that its
fixed data is exactly $\{\alpha(E_p), \alpha(E_q), \alpha(E_r)\}$.
This completes the proof of Theorem~\ref{thm5}. $\hfill\Box$

\section{Examples of actions with four fixed points}\label{s6}

 This section is to
show how to obtain new  $({\Bbb Z}_2)^3$-actions from $({\Bbb
Z}_2)^3$-action $(\phi_3, {\Bbb R}P^3)$, which will play an
important role in the study of the general case with four fixed
points.

\vskip .2cm

Begin with  the standard linear $({\Bbb Z}_2)^3$-action
$(\phi_3,{\Bbb R}P^3)$ with four fixed points $p=[1,0,0,0],
q=[0,1,0,0], r=[0,0,1,0]$,  and $s=[0,0,0,1]$. One can easily read
off the tangent representations at four fixed points, and then its
colored graph $(\Gamma_{(\phi_3,{\Bbb R}P^3)}, \alpha)$ can
explicitly be shown in Figure~\ref{2}:
 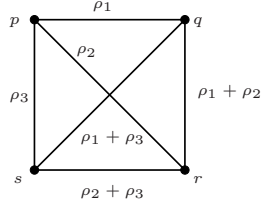
\begin{figure}[h]
   \input{n1.pstex_t}\centering
  \caption[a]{ The colored graph $(\Gamma_{(\phi_3,{\Bbb R}P^3)}, \alpha)$ of $(\phi_3,{\Bbb R}P^3)$}\label{2}
\end{figure}
\vskip .2cm Our  approach to obtain new $({\Bbb Z}_2)^3$-actions is
to modify $(\Gamma_{(\phi_3,{\Bbb R}P^3)}, \alpha)$ into abstract
1-skeleta by adding colored edges on $(\Gamma_{(\phi_3,{\Bbb
R}P^3)}, \alpha)$, and then to use the abstract 1-skeleta to prove
the existence of required new $({\Bbb Z}_2)^3$-actions. We shall see
that only 4- and 5-dimensional new $({\Bbb Z}_2)^3$-actions can be
obtained in such a way.

\subsection{4-dimensional case} 
\begin{figure}[h]
   \input{ld3.pstex_t}\centering
  \caption[a]{Three abstract 1-skeleta of type $(4,3)$ }\label{3}
\end{figure}
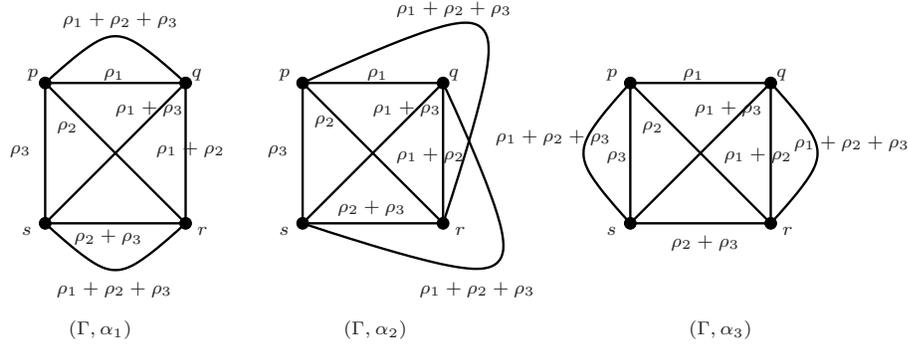

Clearly, by adding two edges to $\Gamma_{(\phi_3,{\Bbb R}P^3)}$ we can produce a unique
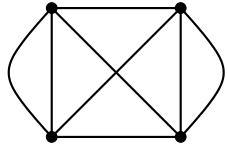
connected regular graph of valence 4 (up to combinatorial
equivalence), denoted by $\Gamma$, shown in the following figure.
 \[   \input{n7.pstex_t}\centering
   \]
Now by adding two colored edges with same coloring
$\rho_1+\rho_2+\rho_3$ into $(\Gamma_{(\phi_3,{\Bbb R}P^3)},
\alpha)$, one can obtain three  abstract 1-skeleta $(\Gamma,
\alpha_1), (\Gamma, \alpha_2), (\Gamma, \alpha_3)$ of type $(4,3)$,
as shown in Figure~\ref{3}.
Obviously, these three  abstract 1-skeleta $(\Gamma, \alpha_1),
(\Gamma, \alpha_2), (\Gamma, \alpha_3)$ have the  same
vertex-coloring set, and in particular, they exactly give all possible abstract
1-skeleta with such a vertex-coloring set.
\begin{lem} \label{4-dim}
There exists a $({\Bbb Z}_2)^3$-action $(\Upsilon, M^4)$
  on a 4-dimensional closed  manifold $M^4$ with four fixed points
  such that its colored graph $(\Gamma_{(\Upsilon, M^4)},
\alpha)$ is one of three  abstract 1-skeleta $(\Gamma, \alpha_1),
(\Gamma, \alpha_2), (\Gamma, \alpha_3)$.
\end{lem}
\begin{proof}
If there is a $({\Bbb Z}_2)^3$-action  on a 4-dimensional closed
manifold such that its colored graph has the same vertex-coloring
set as those three  abstract 1-skeleta $(\Gamma, \alpha_1), (\Gamma,
\alpha_2), (\Gamma, \alpha_3)$, then it must be one
of $(\Gamma, \alpha_1), (\Gamma, \alpha_2), (\Gamma, \alpha_3)$.
Thus, to complete the proof, it suffices to show that the
vertex-coloring set $$\{\alpha_1(E_p), \alpha_1(E_q), \alpha_1(E_r),
\alpha_1(E_s)\}$$ of $(\Gamma, \alpha_1)$ is the fixed data of some
$({\Bbb Z}_2)^3$-action on a closed 4-manifold. By
Theorem~\ref{dks}, this is equivalent to showing that for any
symmetric function $f(x_1,x_2,x_3,x_4)$ of four variables over
${\Bbb Z}_2$
\begin{eqnarray*}
\hat{f}&=& {{f(\alpha_1(E_p))}\over{(\rho_1+\rho_2+\rho_3)
\rho_1\rho_2\rho_3}}+
{{f(\alpha_1(E_q))}\over{(\rho_1+\rho_2+\rho_3)
\rho_1(\rho_1+\rho_2)(\rho_1+\rho_3)}}\\
&& +{{f(\alpha_1(E_r))}\over{(\rho_1+\rho_2+\rho_3)
\rho_2(\rho_1+\rho_2)(\rho_2+\rho_3)}}+
{{f(\alpha_1(E_s))}\over{(\rho_1+\rho_2+\rho_3)
\rho_3(\rho_1+\rho_3)(\rho_2+\rho_3)}}
\end{eqnarray*}
belongs to ${\Bbb Z}_2[\rho_1,\rho_2,\rho_3]$. One may write
$f(x_1,x_2,x_3,x_4)=\sum_{i=0}^jx^i_1f_i(x_2,x_3,x_4)$ to obtain
\begin{eqnarray*}
\hat{f}&=&{1\over{(\rho_1+\rho_2+\rho_3)}}\Bigg\{{{f_0(
\rho_1,\rho_2,\rho_3)}\over{\rho_1\rho_2\rho_3}}+
{{f_0(\rho_1,\rho_1+\rho_2,\rho_1+\rho_3)}\over{
\rho_1(\rho_1+\rho_2)(\rho_1+\rho_3)}}\\
&&+{{f_0(\rho_2, \rho_1+\rho_2,
\rho_2+\rho_3)}\over{\rho_2(\rho_1+\rho_2)(\rho_2+\rho_3)}}+
{{f_0(\rho_3, \rho_1+\rho_3, \rho_2+\rho_3)}\over{\rho_3(\rho_1+\rho_3)(\rho_2+\rho_3)}}\Bigg\}\\
&&+\sum_{i=1}^j(\rho_1+\rho_2+\rho_3)^{i-1}\Bigg\{{{f_i(
\rho_1,\rho_2,\rho_3)}\over{\rho_1\rho_2\rho_3}}+
{{f_i(\rho_1,\rho_1+\rho_2,\rho_1+\rho_3)}\over{
\rho_1(\rho_1+\rho_2)(\rho_1+\rho_3)}}\\
&&+{{f_i(\rho_2, \rho_1+\rho_2,
\rho_2+\rho_3)}\over{\rho_2(\rho_1+\rho_2)(\rho_2+\rho_3)}}+
{{f_i(\rho_3, \rho_1+\rho_3,
\rho_2+\rho_3)}\over{\rho_3(\rho_1+\rho_3)(\rho_2+\rho_3)}}\Bigg\}.
\end{eqnarray*}
Note that for $0\leq i\leq j$, $f_i$ is symmetric.  Then,  from the
colored graph $(\Gamma_{(\phi_3,{\Bbb R}P^3)}, \alpha)$ of the
action $(\phi_3, {\Bbb R}P^3)$, one sees that for $0\leq i\leq j$,
 each term
\begin{eqnarray*}\hat{f}_i&=&{{f_i(
\rho_1,\rho_2,\rho_3)}\over{\rho_1\rho_2\rho_3}}+
{{f_i(\rho_1,\rho_1+\rho_2,\rho_1+\rho_3)}\over{
\rho_1(\rho_1+\rho_2)(\rho_1+\rho_3)}}+{{f_i(\rho_2, \rho_1+\rho_2,
\rho_2+\rho_3)}\over{\rho_2(\rho_1+\rho_2)(\rho_2+\rho_3)}}\\
&&+{{f_i(\rho_3, \rho_1+\rho_3,
\rho_2+\rho_3)}\over{\rho_3(\rho_1+\rho_3)(\rho_2+\rho_3)}}
\end{eqnarray*}
belongs to ${\Bbb Z}_2[\rho_1,\rho_2,\rho_3]$  by Theorem~\ref{dks}.

\vskip .2cm

Clearly, $\hat{f}$ belongs to ${\Bbb Z}_2[\rho_1,\rho_2,\rho_3]$ if
and only if $\hat{f}_0$  is divisible by $\rho_1+\rho_2+\rho_3$. So,
to complete the proof,  it remains to show that $\hat{f}_0$  is
divisible by $\rho_1+\rho_2+\rho_3$. To check that this is true, we
put $\rho_1+\rho_2+\rho_3=0$ and then
\begin{eqnarray*} \hat{f}_0&=&{{f_0(
\rho_1,\rho_2,\rho_3)}\over{\rho_1\rho_2\rho_3}}+
{{f_0(\rho_1,\rho_1+\rho_2,\rho_1+\rho_3)}\over{
\rho_1(\rho_1+\rho_2)(\rho_1+\rho_3)}}
+{{f_0(\rho_2,\rho_1+\rho_2,\rho_2+\rho_3)}\over{
\rho_2(\rho_1+\rho_2)(\rho_2+\rho_3)}}\\ &&+ {{f_0(
\rho_3,\rho_1+\rho_3,\rho_2+\rho_3)}\over{
\rho_3(\rho_1+\rho_3)(\rho_2+\rho_3)}}\\
&=&{{f_0( \rho_1,\rho_2,\rho_3)}\over{\rho_1\rho_2\rho_3}}+{{f_0(
\rho_1,\rho_3,\rho_2)}\over{\rho_1\rho_3\rho_2}}+{{f_0(
\rho_2,\rho_3,\rho_1)}\over{\rho_2\rho_3\rho_1}}+{{f_0(
\rho_3,\rho_2,\rho_1)}\over{\rho_3\rho_2\rho_1}}\\
\end{eqnarray*}
\noindent which is zero since $f_0$ is symmetric. This means that $\hat{f}_0$  is
divisible by $\rho_1+\rho_2+\rho_3$.
\end{proof}
\begin{rem}
 It
is easy to see that any two of $(\Gamma, \alpha_1), (\Gamma,
\alpha_2), (\Gamma, \alpha_3)$ can be translated to each other by
automorphisms of $\text{Hom}(({\Bbb Z}_2)^3,{\Bbb Z}_2)$, i.e., for
any  $\alpha_i$ and $\alpha_j$, there is an automorphism $\theta$
such that the following diagram communts.
 $$
\xymatrix{
                &         \Hom(({\Bbb Z}_2)^3, {\Bbb Z}_2) \ar[dd]^{\theta}     \\
  E_\Gamma \ar[ur]^{\alpha_i} \ar[dr]_{\alpha_j}                 \\
                &          \Hom(({\Bbb Z}_2)^3, {\Bbb Z}_2)                 }
$$
 Thus, each of $(\Gamma,
\alpha_1)$, $(\Gamma, \alpha_2)$, $(\Gamma, \alpha_3)$ can become a
colored graph of some $({\Bbb Z}_2)^3$-action.
\end{rem}

\subsection{5-dimensional case}  
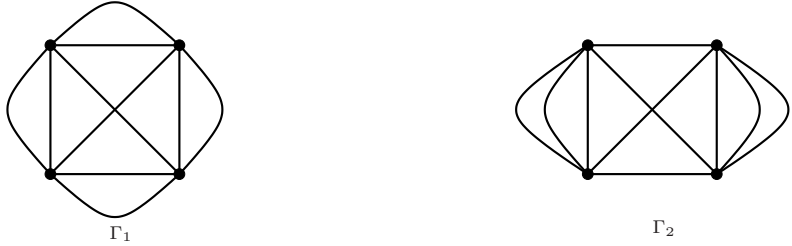
\begin{figure}[h]
   \input{n6.pstex_t}\centering
  \caption[a]{Two regular graphs $\Gamma_1, \Gamma_2$ of
valence 5 }\label{4}
\end{figure}

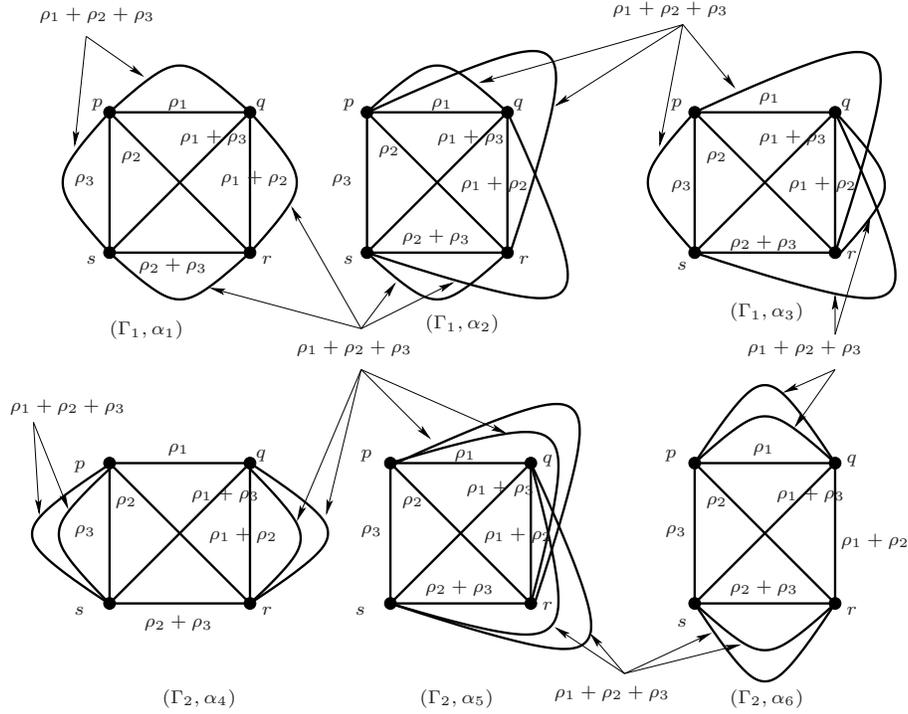
\begin{figure}[h]
   \input{n4.pstex_t}\centering
  \caption[a]{Six abstract 1-skeleta of type $(5,3)$ }\label{ab6}
\end{figure}

An easy observation shows that up
to combinatorial equivalence, by adding four edges to $\Gamma_{(\phi_3,{\Bbb R}P^3)}$
we can merely produce two regular graphs $\Gamma_1, \Gamma_2$ of
valence 5, shown in Figure~\ref{4}.
Furthermore,  by adding four colored edges with same coloring
$\rho_1+\rho_2+\rho_3$ on $(\Gamma_{(\phi_3,{\Bbb R}P^3)}, \alpha)$,
we can obtain six  abstract 1-skeleta of type $(5,3)$, which are
shown in Figure~\ref{ab6}.
These six abstract 1-skeleta have the same vertex-coloring set, and
they give all possible abstract 1-skeleta with such a
vertex-coloring set. Thus, if there is a $({\Bbb Z}_2)^3$-action on
a 5-dimensional closed manifold such that its colored graph has the
same vertex-coloring set as those six abstract 1-skeleta, then its
colored graph must be one of six abstract 1-skeleta.

\begin{lem}\label{5-dim}
There exists a $({\Bbb Z}_2)^3$-action $(\Lambda, M^5)$ on a
5-dimensional closed manifold with four fixed points such that its
colored graph $(\Gamma_{(\Lambda, M^5)}, \alpha)$ is one of six
abstract 1-skeleta in Figure~\ref{4}.\end{lem}

\begin{proof} In a similar way to Lemma~\ref{4-dim}, it suffices to show that the vertex-coloring set
$$\{\alpha_1(E_p), \alpha_1(E_q), \alpha_1(E_r), \alpha_1(E_s)\}$$ of
$(\Gamma_1, \alpha_1)$ is the fixed data of some $({\Bbb
Z}_2)^3$-action on a closed 5-manifold. Let $f(x_1,x_2,x_3,x_4,x_5)$
be a symmetric polynomial function in 5 variables over ${\Bbb Z}_2$
and
\begin{eqnarray*} \hat{f}&=&
{{f(\alpha_1(E_p))}\over{(\rho_1+\rho_2+\rho_3)^2
\rho_1\rho_2\rho_3}}+
{{f(\alpha_1(E_q))}\over{(\rho_1+\rho_2+\rho_3)^2
\rho_1(\rho_1+\rho_2)(\rho_1+\rho_3)}}\\
&& +{{f(\alpha_1(E_r))}\over{(\rho_1+\rho_2+\rho_3)^2
\rho_2(\rho_1+\rho_2)(\rho_2+\rho_3)}}+
{{f(\alpha_1(E_s))}\over{(\rho_1+\rho_2+\rho_3)^2
\rho_3(\rho_1+\rho_3)(\rho_2+\rho_3)}}.
\end{eqnarray*}
By Theorem~\ref{dks} it needs to show that $\hat{f}\in {\Bbb
Z}_2[\rho_1,\rho_2,\rho_3]$ . For this, write
$$f(x_1,x_2,x_3,x_4,x_5)=\sum_{i=0}^jx_1^i\sum_{l=0}^ux_2^lf_{i,l}(x_3,x_4,x_5)=
\sum_{i=0}^j\sum_{l=0}^ux_1^ix_2^lf_{i,l}(x_3,x_4,x_5)$$ with each
$f_{i,l}$ a symmetric function in 3 variables $x_3, x_4, x_5$ and
\begin{eqnarray*}\hat{f}_{i,l}&=&{{f_{i,l}(\rho_1,\rho_2,\rho_3)}\over{\rho_1\rho_2\rho_3}}+
{{f_{i,l}(\rho_1,\rho_1+\rho_2,\rho_1+\rho_3)}\over{\rho_1(\rho_1+\rho_2)(\rho_1+\rho_3)}}+{{f_{i,l}(\rho_2,
\rho_1+\rho_2,
\rho_2+\rho_3)}\over{\rho_2(\rho_1+\rho_2)(\rho_2+\rho_3)}}\\
&&+{{f_{i,l}(\rho_3, \rho_1+\rho_3,
\rho_2+\rho_3)}\over{\rho_3(\rho_1+\rho_3)(\rho_2+\rho_3)}}.
\end{eqnarray*}
Then by direct calculations, one has
\begin{eqnarray*}
\hat{f}&=& {{\hat{f}_{0,0}}\over{(\rho_1+\rho_2+\rho_3)^2}}+
{{\hat{f}_{0,1}+\hat{f}_{1,0}}\over{(\rho_1+\rho_2+\rho_3)}}+\sum_{i+l\geq
2}(\rho_1+\rho_2+\rho_3)^{i+l-2}\hat{f}_{i,l}.
\end{eqnarray*}
One knows from the proof of Lemma~\ref{4-dim} that for any $i$ and
$l$, $\hat{f}_{i,l}$ belongs to ${\Bbb Z}_2[\rho_1,\rho_2,\rho_3]$
and is divisible by $\sigma=\rho_1+\rho_2+\rho_3$. Thus
$\hat{f}\in{\Bbb Z}_2[\rho_1,\rho_2,\rho_3]$ if and only if
$\hat{f}_{0,0}$ is divisible by $\sigma^2$. So
 it remains to prove that $\hat{f}_{0,0}\equiv 0\mod
\sigma^2$. To do this, one will write the symmetric function
$f_{0,0}(x_3,x_4,x_5)$ of 3 variables $x_3,x_4,x_5$ in terms of
elementary symmetric functions $\sigma_1(x_3,x_4,x_5),
\sigma_2(x_3,x_4,x_5)$ and $\sigma_3(x_3,x_4,x_5)$.

\vskip .2cm

Now on the four fixed points, one has that
$\sigma_1(\rho_1,\rho_2,\rho_3)=$
$\sigma_1(\rho_1,\rho_1+\rho_2,\rho_1+\rho_3)= $
$\sigma_1(\rho_2,\rho_1+\rho_2,\rho_2+\rho_3)= $
$\sigma_1(\rho_3,\rho_1+\rho_3,\rho_2+\rho_3)=
\rho_1+\rho_2+\rho_3=\sigma$. Thus, if
$f_{0,0}(x_3,x_4,x_5)=\sigma_1(x_3,x_4,x_5)f'(x_3,x_4,x_5)$ then
$\hat{f}_{0,0}=\sigma\hat{f'}$, and since $\hat{f'}$ is divisible by
$\sigma$ by Lemma~\ref{4-dim}, $\hat{f}_{0,0}$ is divisible by
$\sigma^2$.

\vskip .2cm

So it suffices to only consider
$f_{0,0}(x_3,x_4,x_5)=\sigma_2^v(x_3,x_4,x_5)\sigma_3^w(x_3,x_4,x_5)$.

\vskip .2cm

Write $\rho_3=\rho_1+\rho_2+\sigma$ and consider everything in
${\Bbb Z}_2[\rho_1,\rho_2,\sigma]\cong {\Bbb Z}_2[\rho_1, \rho_2,
\rho_3]$. By direct calculations one has that
\begin{eqnarray*}
& & \sigma_2(\rho_1,\rho_2,\rho_3)= \tau+(\rho_1+\rho_2)\sigma\\ &
&\sigma_2(\rho_1,\rho_1+\rho_2,\rho_1+\rho_3)= \tau+\rho_2\sigma\\ &
& \sigma_2(\rho_2,\rho_1+\rho_2,\rho_2+\rho_3)= \tau+\rho_1\sigma\\
&
&\sigma_2(\rho_3,\rho_1+\rho_3,\rho_2+\rho_3)=\tau+\sigma^2\equiv\tau\mod
\sigma^2
\end{eqnarray*}
where $\tau=\rho_1^2+\rho_1\rho_2+\rho_2^2$, and
\begin{eqnarray*}
&&\sigma_3(\rho_1,\rho_2,\rho_3)=\rho_1\rho_2\rho_3= \varphi+\rho_1\rho_2\sigma\\
&&\sigma_3(\rho_1,\rho_1+\rho_2,\rho_1+\rho_3)=\rho_1(\rho_1+\rho_2)(\rho_1+\rho_3)=\varphi+(\rho_1^2+\rho_1\rho_2)\sigma\\
&&\sigma_3(\rho_2,\rho_1+\rho_2,\rho_2+\rho_3)=\rho_2(\rho_1+\rho_2)(\rho_2+\rho_3)=\varphi+(\rho_2^2+\rho_1\rho_2)\sigma\\
&&\sigma_3(\rho_3,\rho_1+\rho_3,\rho_2+\rho_3)=\rho_3(\rho_1+\rho_3)(\rho_2+\rho_3)\\&&\equiv\varphi+(\rho_1^2+\rho_1\rho_2+\rho_2^2)\sigma\mod
\sigma^2
\end{eqnarray*}
where $\varphi=\rho_1\rho_2(\rho_1+\rho_2)$. Beginning with the
easiest case, take
$f_{0,0}(x_3,x_4,x_5)=\sigma_2^v(x_3,x_4,x_5)\sigma_3^w(x_3,x_4,x_5)$
with $w\geq 1$. This cancels the denominators in the expression of
$\hat{f}_{0,0}$ and
\begin{eqnarray*}
\hat{f}_{0,0}&\equiv&
[\tau+(\rho_1+\rho_2)\sigma]^v[\varphi+\rho_1\rho_2\sigma]^{w-1}+[\tau+\rho_2\sigma]^v
[\varphi+(\rho_1^2+\rho_1\rho_2)\sigma]^{w-1}\\
&&+[\tau+\rho_1\sigma]^v[\varphi+(\rho_2^2+\rho_1\rho_2)\sigma]^{w-1}
+\tau^v[\varphi+(\rho_1^2+\rho_1\rho_2+\rho_2^2]^{w-1}\\
&\equiv&
4\tau^v\varphi^{w-1}+v\tau^{v-1}\varphi^{w-1}[(\rho_1+\rho_2)\sigma+\rho_2\sigma+\rho_1\sigma+0]+(w-1)\tau^v\varphi^{w-2}\\
&&\times[\rho_1\rho_2\sigma+(\rho_1^2+\rho_1\rho_2)\sigma+
(\rho_2^2+\rho_1\rho_2)\sigma+(\rho_1^2+\rho_1\rho_2+\rho_2^2)\sigma]\\
&\equiv& 0\mod \sigma^2
\end{eqnarray*}
For the case $f_{0,0}(x_3,x_4,x_5)=\sigma_2^v(x_3,x_4,x_5)$, it is
convenient to take the common denominator in $\hat{f}_{0,0}$ to be
the product of the four choices for $\sigma_3(x_3,x_4,x_5)$; i.e.,
$$(\varphi+\rho_1\rho_2\sigma)[\varphi+(\rho_1^2+\rho_1\rho_2)\sigma][\varphi+
(\rho_2^2+\rho_1\rho_2)\sigma][\varphi+(\rho_1^2+\rho_1\rho_2+\rho_2^2)\sigma]$$
and then the numerator in $\hat{f}_{0,0}$ becomes (modulo
$\sigma^2$)
\begin{eqnarray*}
&&[\tau+(\rho_1+\rho_2)\sigma]^v[\varphi^3+\varphi^2((\rho_1^2+\rho_1\rho_2)\sigma+
(\rho_2^2+\rho_1\rho_2)\sigma+(\rho_1^2+\rho_1\rho_2+\rho_2^2)\sigma)]\\
&&+[\tau+\rho_2\sigma]^v[\varphi^3+\varphi^2(\rho_1\rho_2\sigma+(\rho_2^2+\rho_1\rho_2)\sigma
+(\alpha_1^2+\alpha_1\alpha_2+\alpha_2^2)\sigma)]\\
&&+[\tau+\rho_1\sigma]^v[\varphi^3+\varphi^2(\rho_1\rho_2\sigma+(\rho_1^2+\rho_1\rho_2)\sigma
+(\rho_1^2+\rho_1\rho_2+\rho_2^2)\sigma)]\\
&&+\tau^v[\varphi^3+\varphi^2(\rho_1\rho_2\sigma+(\rho_1^2+\rho_1\rho_2)\sigma
+(\rho_2^2+\rho_1\rho_2)\sigma)]\\
&\equiv&
4\tau^v\varphi^3+v\tau^{v-1}\varphi^3[(\rho_1+\rho_2)\sigma+\rho_2\sigma+\rho_1\sigma+0]\\
&&+\tau^v\varphi^2[12\rho_1\rho_2\sigma+6\rho_1^2\sigma+6\rho_2^2\sigma]\\
&\equiv& 0\mod \sigma^2.
\end{eqnarray*}
Thus $\hat{f}_{0,0}$ is always divisible by $\sigma^2$.
\end{proof}
\begin{rem}
We  see easily  that any two of $(\Gamma_1, \alpha_1), (\Gamma_1,
\alpha_2), (\Gamma_1, \alpha_3)$ can be translated to each other by
automorphisms of $\text{\rm Hom}(({\Bbb Z}_2)^3, {\Bbb Z}_2)$. Thus,
if one of $(\Gamma_1, \alpha_1)$, $(\Gamma_1, \alpha_2)$,
$(\Gamma_1, \alpha_3)$ can become a colored graph of some $({\Bbb
Z}_2)^3$-action, then so are all three abstract 1-skeleta
$(\Gamma_1, \alpha_1)$, $(\Gamma_1, \alpha_2)$, $(\Gamma_1,
\alpha_3)$. This is also true for $(\Gamma_2, \alpha_4)$,
$(\Gamma_2, \alpha_5)$, $(\Gamma_2, \alpha_6)$. However, we don't
know whether  two kinds of abstract 1-skeleta all can  become the
colored graphs of $({\Bbb Z}_2)^3$-actions, and in addition,  we
cannot determine which of six abstract 1-skeleta is the colored
graph of $(\Lambda, M^5)$.
\end{rem}

\begin{rem} \label{r9}
 In the proof of Lemma~\ref{5-dim}, when $$f_{0,0}(x_3,x_4,x_5)=\sigma_2(x_3,x_4,x_5)\sigma_3(x_3,x_4,x_5)$$ we have that
$$\hat{f}_{0,0}=[\tau+(\rho_1+\rho_2)\sigma]+(\tau+\rho_2\sigma)+(\tau+\rho_1\sigma)+(\tau+\sigma^2)=\sigma^2.$$
Thus, we cannot improve the divisibility of $\hat{f}_{0,0}$. This
also implies that by the above method,  we cannot further modify
$(\phi_3,{\Bbb R}P^3)$ to obtain a 6-dimensional example with four
fixed points.
\end{rem}

Finally, applying the $\Delta$-operation and $\Omega$-operation to
$(\phi_3, {\Bbb R}P^3)$, $(\Upsilon, M^4)$ and $(\Lambda, M^5)$
gives the following result.

\begin{cor}
When $n=3\cdot 2^\ell, 4\cdot 2^\ell, 5\cdot 2^\ell$ with $\ell\geq
k-3\geq 0$, $\mathcal{A}_n^k(4)$ is nonempty.
\end{cor}

\section{Actions with four fixed points}\label{s7}

Suppose that $\mathcal{A}_n^k(4)$ is nonempty. Then $k\geq 3$ since
any nonbounding $({\Bbb Z}_2)^2$-action  cannot exactly fix four
isolated points by the work of Conner and Floyd \cite[Theorem
31.2]{cf} (also see Lemma~\ref{k=2} in Section~\ref{s4}).  Given a
$G$-action $(\Phi, M^n)$ in $\mathcal{A}_n^k(4)$, without the loss
of generality assume that $(\Phi, M^n)$ has exactly 4 fixed points
$p, q,r, s$. Let $(\Gamma_{(\Phi, M^n)},\alpha)$ be a colored graph
of $(\Phi, M^n)$. Then   the tangent representation set of $(\Phi,
M^n)$ is $\mathcal{N}_{(\Phi, M^n)}=\{\alpha(E_p), \alpha(E_q),
\alpha(E_r), \alpha(E_s)\}$.

\subsection{The simple description of $\mathcal{N}_{(\Phi,
M^n)}$}\label{descr}
 Consider a nontrivial irreducible representation $\rho$ in
$\Hom(G,{\Bbb Z}_2)$ and a component $C$ of fixed set of $\ker\rho$.
Then $C$ must contain an even number of fixed points of $(\Phi,
M^n)$, so $C$ contains either all 4 fixed points or 2 fixed points
of $M^G=\{p,q,r,s\}$.

\vskip .2cm

If $C$ contains all 4 fixed points, then $\rho$ appears the same
number of times in the tangent representation at each fixed point.
Actually, $\rho$ appears exactly $\dim C$ times in the tangent
representation at each fixed point. Obviously,  $\dim C$ must be
more than 1, so generally the pair $(\rho, C)$ doesn't determine a
unique connected $\dim C$-valent regular graph $\Gamma_{\rho, C}$
except for $\dim C=2$. Let $\hat{\omega}$ be the multiset of
representations written in $\Hom(G,{\Bbb Z}_2)$ occurring for all
such components $C$ (with multiplicities $\dim C$'s).

\vskip .2cm

If $C$ contains 2 fixed points, then $\rho$ appears the same number
of times in the tangent representations at these two fixed points.
Each pair $(\rho, C)$ determines a unique $\dim C$-valent regular
graph $\Gamma_{\rho, C}$ with two vertices.  Denote the
representations occurring for all such components $C$ with
$\Gamma_{\rho, C}$ containing $p$ and $q$ by $\hat{\beta}$, $p$ and
$r$ by $\hat{\gamma}$, $p$ and $s$ by $\hat{\delta}$, $q$ and $r$ by
$\hat{\varepsilon}$, $q$ and $s$ by $\hat{\eta}$, $r$ and $s$ by
$\hat{\lambda}$, respectively, as shown in Figure~\ref{5}.
\begin{figure}[h]
    \input{n8.pstex_t}\centering
    \caption[a]{\tiny }\label{5}
\end{figure}
Then we can write
$$\alpha(E_p)=\hat{\omega}\cup\hat{\beta}\cup\hat{\gamma}\cup
\hat{\delta}, \
\alpha(E_q)=\hat{\omega}\cup\hat{\beta}\cup\hat{\varepsilon}\cup
\hat{\eta}, \
\alpha(E_r)=\hat{\omega}\cup\hat{\varepsilon}\cup\hat{\gamma}\cup
\hat{\lambda}, \
\alpha(E_s)=\hat{\omega}\cup\hat{\delta}\cup\hat{\eta}\cup
\hat{\lambda}.$$ Obviously, each of the above four expressions
consists of the union of four disjoint sets.
 The
only possible nonempty intersections are
$\hat{\beta}\cap\hat{\lambda}$, $\hat{\gamma}\cap\hat{\eta}$, and
$\hat{\delta}\cap\hat{\varepsilon}$.

\vskip .2cm

Furthermore, we have that $$\begin{cases}
n=\vert\hat{\omega}\vert+\vert\hat{\beta}\vert+\vert\hat{\gamma}\vert+\vert\hat{\delta}\vert\\
n=\vert\hat{\omega}\vert+\vert\hat{\beta}\vert+\vert\hat{\varepsilon}\vert+\vert\hat{\eta}\vert\\
n=\vert\hat{\omega}\vert+\vert\hat{\varepsilon}\vert+\vert\hat{\gamma}\vert+\vert\hat{\lambda}\vert\\
n=\vert\hat{\omega}\vert+\vert\hat{\delta}\vert+\vert\hat{\eta}\vert+\vert\hat{\lambda}\vert\\
\end{cases}$$
and thus
\begin{eqnarray} \label{eq1}
\vert\hat{\beta}\vert=\vert\hat{\lambda}\vert,\ \ \
\vert\hat{\gamma}\vert=\vert\hat{\eta}\vert,\ \ \
\vert\hat{\delta}\vert=\vert\hat{\varepsilon}\vert.\end{eqnarray}
Now let $\hat{\beta}_0=\hat{\lambda}_0$ be the common part of
$\hat{\beta}$ and $\hat{\lambda}$ and let
$\hat{\beta}_1=\hat{\beta}\setminus\hat{\beta}_0$ and
$\hat{\lambda}_1=\hat{\lambda}\setminus\hat{\lambda}_0$.   By
(\ref{eq1}) we then have
$$\vert\hat{\beta}_1\vert=\vert\hat{\beta}\vert-\vert\hat{\beta}_0\vert=\vert\hat{\lambda}\vert-
\vert\hat{\lambda}_0\vert=\vert\hat{\lambda}_1\vert.$$ Similarly we
form $\hat{\gamma}_0, \hat{\gamma}_1, \hat{\eta}_0, \hat{\eta}_1$
with $\hat{\gamma}_0=\hat{\eta}_0$ and
$|\hat{\gamma}_1|=|\hat{\eta}_1|$, and $\hat{\delta}_0,
\hat{\delta}_1$, $\hat{\varepsilon}_0, \hat{\varepsilon}_1$ with
$\hat{\delta}_0=\hat{\varepsilon}_0$ and
$|\hat{\delta}_1|=|\hat{\varepsilon}_1|$. Then four tangent
representations of $\mathcal{N}_{(\Phi, M^n)}$
 become
$$\alpha(E_p)=\hat{\omega}\cup\hat{\beta}_0\cup\hat{\gamma}_0\cup
\hat{\delta}_0\cup\hat{\beta}_1\cup\hat{\gamma}_1\cup
\hat{\delta}_1,$$
$$\alpha(E_q)=\hat{\omega}\cup\hat{\beta}_0\cup\hat{\gamma}_0\cup
\hat{\delta}_0\cup\hat{\beta}_1\cup\hat{\varepsilon}_1\cup
\hat{\eta}_1,$$
$$\alpha(E_r)=\hat{\omega}\cup\hat{\beta}_0\cup\hat{\gamma}_0\cup
\hat{\delta}_0\cup\hat{\varepsilon}_1\cup\hat{\gamma}_1\cup
\hat{\lambda}_1,$$
$$\alpha(E_s)=\hat{\omega}\cup\hat{\beta}_0\cup\hat{\gamma}_0\cup
\hat{\delta}_0\cup\hat{\delta}_1\cup\hat{\eta}_1\cup
\hat{\lambda}_1$$ with
$\vert\hat{\beta}_1\vert=\vert\hat{\lambda}_1\vert$,
$|\hat{\gamma}_1|=|\hat{\eta}_1|$ and
$|\hat{\delta}_1|=|\hat{\varepsilon}_1|$.

\vskip .2cm Note that if $\hat{\omega}$ is empty, then
$\Gamma_{(\Phi, M^n)}$ is uniquely determined.

\subsection{The essential structure of $\mathcal{N}_{(\Phi, M^n)}$}
\label{ess-str} Following the notions of the above subsection, now
let us study the essential structure of $\mathcal{N}_{(\Phi, M^n)}$.

\begin{lem} \hskip.1cm \label{lm7}
\begin{enumerate}
\item[$(a)$]
 $\hat{\beta}_1, \hat{\gamma}_1,
\hat{\delta}_1,\hat{\varepsilon}_1, \hat{\eta}_1, \hat{\lambda}_1$
are all nonempty.
\item[$(b)$] The sets $\hat{\omega}, \hat{\beta}_0=\hat{\lambda}_0,
\hat{\gamma}_0=\hat{\eta}_0, \hat{\delta}_0=\hat{\varepsilon}_0,
\hat{\beta}_1, \hat{\gamma}_1, \hat{\delta}_1,\hat{\varepsilon}_1,
\hat{\eta}_1, \hat{\lambda}_1$ are all disjoint.\end{enumerate}
\end{lem}
\begin{proof}
Since no two of the fixed points have the same representation, one
of the sets $\hat{\beta}_1, \hat{\gamma}_1,
\hat{\delta}_1,\hat{\varepsilon}_1, \hat{\eta}_1, \hat{\lambda}_1$
must be nonempty. Say $\hat{\beta}_1\not=\emptyset$ (then also
$\hat{\lambda}_1\not=\emptyset$). Then by comparing $\alpha(E_p)$
and $\alpha(E_q)$ one of $\hat{\gamma}_1, \hat{\delta}_1,
\hat{\varepsilon}_1, \hat{\eta}_1$ must be nonempty. Say
$\hat{\gamma}_1\not=\emptyset$ (then also
$\hat{\eta}_1\not=\emptyset$).

\vskip .2cm

Let $\beta\in\hat{\beta}$ and $\gamma\in \hat{\gamma}_1$, and
consider $\beta+\gamma$. The representation $\gamma$ occurs more
times in $\alpha(E_p)$ than in $\alpha(E_q)$, but the number of
times that $\gamma$ and $\beta+\gamma$ occur in $\alpha(E_p)$ is the
same as the number of times that $\gamma$ and $\beta+\gamma$ occur
in  $\alpha(E_q)$ by the property (P2) in Section~\ref{s2}, so
$\beta+\gamma$ must occur more times in $\alpha(E_q)$ than in
$\alpha(E_p)$. Thus $\beta+\gamma$ belongs to
$\hat{\eta}_1\cup\hat{\varepsilon}_1$. If $\beta\in\hat{\beta}_0$,
then $\beta\in \hat{\lambda}$ and $\gamma$ occurs more times in
$\alpha(E_r)$ than in $\alpha(E_s)$, so $\beta+\gamma$ occurs more
times in $\alpha(E_s)$ than in $\alpha(E_r)$, and then
$\beta+\gamma$ belongs to $\hat{\delta}_1\cup\hat{\eta}_1$. But
$\beta+\gamma$ is in both $\hat{\eta}_1\cup\hat{\varepsilon}_1$ and
$\hat{\delta}_1\cup\hat{\eta}_1$, so $\beta+\gamma$ belongs to
$\hat{\eta}_1$. If $\beta\in\hat{\beta}_1$, then $\beta$ occurs more
times in $\alpha(E_p)$ than in $\alpha(E_r)$ so $\beta+\gamma$ must
occur more times in $\alpha(E_r)$ than in $\alpha(E_p)$, and
$\beta+\gamma$ belongs to $\hat{\lambda}_1\cup\hat{\varepsilon}_1$.
But  $\beta+\gamma$ is in both $\hat{\eta}_1\cup\hat{\varepsilon}_1$
and $\hat{\lambda}_1\cup\hat{\varepsilon}_1$, so $\beta+\gamma$
belongs to $\hat{\varepsilon}_1$. Since we had supposed that
$\hat{\beta}_1$ and $\hat{\gamma}_1$ are nonempty, we have that
$\hat{\varepsilon}_1$ is nonempty so all of the sets $\hat{\beta}_1,
\hat{\gamma}_1, \hat{\delta}_1,\hat{\varepsilon}_1, \hat{\eta}_1,
\hat{\lambda}_1$ are nonempty.

\vskip .2cm

Further, if $\beta\in \hat{\beta}_0$ and $\beta\in \hat{\beta}_1$,
then $\beta+\gamma$ belongs to both $\hat{\eta}_1$ and
$\hat{\varepsilon}_1$, which is impossible. Thus $\hat{\beta}_0$ and
$\hat{\beta}_1$ must be disjoint. Similarly, $\hat{\gamma}_0$ and
$\hat{\gamma}_1$ are disjoint, and so on. Thus, the sets
$$\hat{\omega}, \hat{\beta}_0=\hat{\lambda}_0,
\hat{\gamma}_0=\hat{\eta}_0, \hat{\delta}_0=\hat{\varepsilon}_0,
\hat{\beta}_1, \hat{\gamma}_1, \hat{\delta}_1,\hat{\varepsilon}_1,
\hat{\eta}_1, \hat{\lambda}_1$$ are all disjoint.
\end{proof}

Now let us analyze the structures of $\hat{\beta}_1, \hat{\gamma}_1,
\hat{\delta}_1,\hat{\varepsilon}_1, \hat{\eta}_1, \hat{\lambda}_1$.

\begin{lem}\hskip .1cm \label{lm8}
\begin{enumerate}
\item[$(a)$] $\hat{\beta}_1, \hat{\gamma}_1,
\hat{\delta}_1,\hat{\varepsilon}_1, \hat{\eta}_1, \hat{\lambda}_1$
have the same number of elements, and all elements of
$\hat{\beta}_1, \hat{\gamma}_1, \hat{\delta}_1,\hat{\varepsilon}_1,
\hat{\eta}_1, \hat{\lambda}_1$  occur with the same multiplicity.
\item[$(b)$] All odd sums of elements of $\hat{\beta}_1$ are again in
$\hat{\beta}_1$ and choose $\gamma\in\hat{\gamma}_1$ and
$\delta\in\hat{\delta}_1$, one has that
\begin{equation} \label{eq2}
\begin{cases} \hat{\gamma}_1=\{\gamma+\beta+\beta'|\beta, \beta'\in \hat{\beta}_1\}\\
\hat{\delta}_1=\{\delta+\beta+\beta'|\beta, \beta'\in \hat{\beta}_1\}\\
\hat{\varepsilon}_1=\{\gamma+\beta|\beta\in\hat{\beta}_1\}\\
\hat{\eta}_1=\{\delta+\beta|\beta\in\hat{\beta}_1\}\\
\hat{\lambda}_1=\{\gamma+\beta+\beta'|\beta, \beta'\in \hat{\beta}_1\}\\
\end{cases}
\end{equation}
and $\gamma,\delta$ are linearly independent of the span of
$\hat{\beta}_1$.
\end{enumerate}
\end{lem}

\begin{proof}
First, by Figure~\ref{5}, Lemma~\ref{lm7} and the property (P2) in
Section~\ref{s2}, we easily obtain that for any
$\beta\in\hat{\beta}_1, \gamma\in \hat{\gamma}_1, \delta\in
\hat{\delta}_1,\varepsilon\in \hat{\varepsilon}_1,
\eta\in\hat{\eta}_1, \lambda\in\hat{\lambda}_1$,
\begin{eqnarray} \label{eq3}
\begin{cases}\beta+\gamma\in\hat{\varepsilon}_1, \beta+\delta\in
\hat{\eta}_1, \beta+\eta\in \hat{\delta}_1, \beta+\varepsilon\in
\hat{\gamma}_1, \gamma+\delta\in\hat{\lambda}_1,
\gamma+\varepsilon\in\hat{\beta}_1,\\
\gamma+\lambda\in\hat{\delta}_1, \delta+\eta\in\hat{\beta}_1,
\delta+\lambda\in\hat{\gamma}_1, \varepsilon+\eta\in\hat{\lambda}_1,
\varepsilon+\lambda\in\hat{\eta}_1,
\eta+\lambda\in\hat{\varepsilon}_1. \end{cases} \end{eqnarray} Now
for $\gamma\in\hat{\gamma}_1$ and $\beta\in\hat{\beta}_1$, the
number of times that $\gamma$ occurs in $\alpha(E_p)$ is the same as
the number of times that $\beta+\gamma$ occurs in $\alpha(E_q)$, so
the number of times that $\gamma$ occurs in $\hat{\gamma}_1$ is the
same as the number of times that $\beta+\gamma$ occurs in
$\hat{\varepsilon}_1$. Also, by (\ref{eq3}) if $\varepsilon\in
\hat{\varepsilon}_1$, then $\beta+\varepsilon\in \hat{\gamma}_1$.
Thus if $\hat{\gamma}_1=\{\gamma_1,...,\gamma_t\}$ then
$\hat{\varepsilon}_1=\{\gamma_1+\beta,...,\gamma_t+\beta\}$. Thus
$\vert\hat{\beta}_1\vert=\vert\hat{\varepsilon}_1\vert$ and
similarly, $\vert\hat{\beta}_1\vert=\vert\hat{\eta}_1\vert$, so all
of the sets $\hat{\beta}_1, \hat{\gamma}_1,
\hat{\delta}_1,\hat{\varepsilon}_1, \hat{\eta}_1, \hat{\lambda}_1$
have the same number of elements.

\vskip .2cm Moreover, the number of times that $\beta+\gamma$ occurs
in $\hat{\varepsilon}_1$ is not only the same as the number of times
that $\gamma$ occurs in $\hat{\gamma}_1$,  and but it also is the
same as the number of times that $\beta$ occurs in $\hat{\beta}_1$.
Thus, the elements of $\hat{\beta}_1, \hat{\gamma}_1,
\hat{\varepsilon}_1$ all occur with the same multiplicity.
Similarly, this also happens for $\hat{\delta}_1, \hat{\gamma}_1,
\hat{\lambda}_1$ (resp. $\hat{\beta}_1,\hat{\delta}_1,
\hat{\eta}_1$). Therefore, the elements of $\hat{\beta}_1,
\hat{\gamma}_1, \hat{\delta}_1,\hat{\varepsilon}_1, \hat{\eta}_1,
\hat{\lambda}_1$ all occur with the same multiplicity.

\vskip .2cm By (\ref{eq3}) and the property (P2) of Section~\ref{s2}
one has that $\hat{\eta}_1=\{\delta+\beta|\beta\in\hat{\beta}_1\}$,
and for any $\eta\in\hat{\eta}_1$ and $\beta'\in\hat{\beta}_1$,
$\eta+\beta'\in\hat{\delta}_1$ and
$\hat{\eta}_1=\{\eta+\beta+\beta'|\beta\in\hat{\beta}_1\}$. Thus,
for any $\beta''\in\hat{\beta}_1$, choosing $\eta=\delta+\beta''$
gives
$$\hat{\beta}_1=\{\beta|\beta\in\hat{\beta}_1\}=\{\beta+\beta'+\beta''|\beta\in\hat{\beta}_1\}.$$
Replacing $\beta$ by $\beta+\beta'+\beta''$ and repeating the above
procedure, we have that every sum of an odd number of elements in
$\hat{\beta}_1$ is again in $\hat{\beta}_1$. Then (\ref{eq2})
follows from this and (\ref{eq3}).

\vskip .2cm

Finally, let us show that $\gamma,\delta$ are linearly independent
of the span of $\hat{\beta}_1$. This is equivalent to proving that
$\gamma+\delta$ cannot be a sum of elements of $\hat{\beta}_1$. If
$\gamma+\delta$ is an odd sum of elements of $\hat{\beta}_1$, then
$\gamma+\delta\in \hat{\lambda}_1\cap\hat{\beta}_1$ so
$\hat{\lambda}_1\cap\hat{\beta}_1$ is nonempty. This is impossible
by Lemma~\ref{lm7}(b). If $\gamma+\delta$ is an even sum of elements
of $\hat{\beta}_1$, then $\gamma\in
\hat{\delta}_1\cap\hat{\gamma}_1$, but this is also impossible.
\end{proof}
\begin{rem}\label{m}
Let $m$ be the multiplicity of  each element of $\hat{\beta}_1$. In
a similar way to the proof of Lemma~\ref{lm5}, it is easy to check
that all elements of $\hat{\beta}_1$ span a $(k-2)$-dimensional
subspace in $\Hom(({\Bbb Z}_2)^k, {\Bbb Z}_2)$, so
$|\hat{\beta}_1|=m\cdot 2^{k-3}$. Thus,
$|\hat{\gamma}_1|=|\hat{\delta}_1|=|\hat{\varepsilon}_1|=|\hat{\eta}_1|=|\hat{\lambda}_1|=m\cdot
2^{k-3}$.
\end{rem}

\begin{defn} Let
$$\underline{\hat{\omega}}:=\hat{\omega}\cup\hat{\beta}_0\cup\hat{\gamma}_0\cup\hat{\delta}_0$$
which is called the {\em changeable part} of $\mathcal{N}_{(\Phi,
M^n)}$. \end{defn} Indeed,
 $(\phi_3,
{\Bbb R}P^3)$ is an example with $\underline{\hat{\omega}}$ being
empty, while Lemmas~\ref{4-dim} and~\ref{5-dim} in Section~\ref{s6}
provide examples with $\underline{\hat{\omega}}$
 nonempty.

\begin{lem}\label{lm9}
 When $|\underline{\hat{\omega}}|>0$,  every element in
$\underline{\hat{\omega}}$ has the form $\gamma+\delta+\beta$,
$\beta\in \hat{\beta}_1$, where $\gamma\in\hat{\gamma}_1$ and
$\delta\in\hat{\delta}_1$.
\end{lem}
\begin{proof}   Let $\xi\in
\underline{\hat{\omega}}$, $\gamma\in\hat{\gamma}_1$ and
$\delta\in\hat{\delta}_1$. If $\xi\in
\hat{\omega}\cup\hat{\beta}_0\cup\hat{\delta}_0$ but
$\xi\not\in\hat{\gamma}_0$, since $\gamma\in \alpha(E_p)$ and
$\gamma\not\in \alpha(E_q)$, one then has that $\gamma+\xi\in
\alpha(E_q)$ so $\gamma+\xi\in \hat{\varepsilon}_1\cup\hat{\eta}_1$
by Lemma~\ref{lm7}. On the other hand, $\gamma\in \alpha(E_p)$ but
$\gamma\not\in \alpha(E_s)$, so $\gamma+\xi\in \alpha(E_s)$ and
$\gamma+\xi\in \hat{\lambda}_1\cup\hat{\eta}_1$. Thus
$\gamma+\xi\in\hat{\eta}_1$. Further, by Lemma~\ref{lm8}, there is
an element $\beta\in \hat{\beta}_1$ such that
$\gamma+\xi=\delta+\beta$, so $\xi=\gamma+\delta+\beta$ as desired.
If $\xi\in \hat{\gamma}_0$, since $\delta\in \alpha(E_p)$ but
$\delta\not\in \alpha(E_q)$ and $\delta\not\in \alpha(E_r)$,
similarly to the above argument, one has that $\delta+\xi\in
\hat{\eta}_1\cup\hat{\varepsilon}_1$ and $\delta+\xi\in
\hat{\gamma}_1\cup\hat{\varepsilon}_1$, so
$\delta+\xi\in\hat{\varepsilon}_1$ and then by Lemma~\ref{lm8},
there is an element $\beta'\in \hat{\beta}_1$ such that
$\delta+\xi=\gamma+\beta'$. Thus $\xi=\gamma+\delta+\beta'$ as
desired.
\end{proof}

\begin{rem}\label{r10}
 The sum of two elements in
$\underline{\hat{\omega}}$ is an even sum of elements in
$\hat{\beta}_1$ and cannot occur as an element in $\alpha(E_p),
\alpha(E_q), \alpha(E_r)$, and $\alpha(E_s)$.
\end{rem}

Together with Lemmas~\ref{lm8}-\ref{lm9} and Remark~\ref{m} one has

\begin{prop}\label{lm10}
There is a basis $\{\beta_1,...,\beta_{k-2}, \gamma,
\delta\}$ of $\Hom({\Bbb Z}_2)^k,{\Bbb Z}_2)$ and an integer $m$
such that
\begin{enumerate}
\item[$(i)$]
 $n=3m\cdot 2^{k-3}+|\underline{\hat{\omega}}|$;
\item[$(ii)$] $\hat{\beta}_1$ is a multiset consisting of all odd sums with  same multiplicity $m$ formed by
$\beta_1,...,\beta_{k-2}$, and
$$\begin{cases}
\hat{\gamma}_1=\{\gamma+\beta_1+\beta|\beta\in\hat{\beta}_1\}\\
\hat{\delta}_1=\{\delta+\beta_1+\beta|\beta\in\hat{\beta}_1\}\\
\hat{\varepsilon}_1=\{\gamma+\beta|\beta\in\hat{\beta}_1\}\\
\hat{\eta}_1=\{\delta+\beta|\beta\in\hat{\beta}_1\}\\
\hat{\lambda}_1=\{\gamma+\delta+\beta_1+\beta|\beta\in\hat{\beta}_1\}\\
\end{cases}$$
\item[$(iii)$] Every element in
$\underline{\hat{\omega}}$ has the form $\gamma+\delta+\beta$,
$\beta\in\hat{\beta}_1$
 if
$\underline{\hat{\omega}}$ is nonempty. \end{enumerate}
\end{prop}

\begin{defn}
If $\underline{\hat{\omega}}=\emptyset$, then $(\Phi, M^n)$ is said
to be {\em $\underline{\hat{\omega}}$-empty}; otherwise, it is said
to be {\em $\underline{\hat{\omega}}$-nonempty}.
\end{defn}

\begin{rem}\label{classify}
 If
$\underline{\hat{\omega}}=\emptyset$, then we see from
Proposition~\ref{lm10} that the tangent representation set
$\mathcal{N}_{(\Phi, M^n)}=\{\alpha(E_p)$, $\alpha(E_q)$,
$\alpha(E_r)$, $\alpha(E_s)\}$ is completely determined by the basis
$\{\beta_1,...,\beta_{k-2}, \gamma, \delta\}$ of $\Hom(({\Bbb
Z}_2)^k,{\Bbb Z}_2)$. Thus, for any two
$\underline{\hat{\omega}}$-empty actions $(\Phi_1, M_1^n)$ and
$(\Phi_2, M_2^n)$ in $\mathcal{A}_n^k(4)$, there is an automorphism
$\sigma\in \text{\rm GL}(k,{\Bbb Z}_2)$ such that $(\Phi_1, M_1^n)$
is equivariantly cobordant to $(\sigma\Phi_2, M_2^n)$.
\end{rem}

 To keep the notation
manageable, $\prod\hat{\beta}_{1}$ means
$\prod_{\beta\in\hat{\beta}_{1}}\beta$, and similarly for
$\prod\hat{\gamma}_{1}$, $\prod\hat{\delta}_{1}$,
$\prod\hat{\varepsilon}_{1}$, $\prod\hat{\eta}_{1}$,
$\prod\hat{\lambda}_{1}$, $\prod\underline{\hat{\omega}}$. Let
$\hat{\beta}_{11}$ (resp. $\hat{\gamma}_{11}$, $\hat{\delta}_{11}$,
$\hat{\varepsilon}_{11}$, $\hat{\eta}_{11}$, $\hat{\lambda}_{11}$)
denote the set of consisting of  all $2^{k-3}$ different elements in
$\hat{\beta}_1$ (resp. $\hat{\gamma}_{1}$, $\hat{\delta}_{1}$,
$\hat{\varepsilon}_{1}$, $\hat{\eta}_{1}$, $\hat{\lambda}_{1}$).
Then, by Proposition~\ref{lm10} one has that
$\prod\hat{\beta}_{1}=(\prod\hat{\beta}_{11})^m$ where
$\prod\hat{\beta}_{11}$ means
$\prod_{\beta\in\hat{\beta}_{11}}\beta$, and similarly
$\prod\hat{\gamma}_{1}=(\prod\hat{\gamma}_{11})^m$,
$\prod\hat{\delta}_{1}=(\prod\hat{\delta}_{11})^m$,
$\prod\hat{\varepsilon}_{1}=(\prod\hat{\varepsilon}_{11})^m$,
$\prod\hat{\eta}_{1}=(\prod\hat{\eta}_{11})^m$,
$\prod\hat{\lambda}_{1}=(\prod\hat{\lambda}_{11})^m$.
\begin{cor} \label{data2}
$\hat{\beta}_{11}\cup\hat{\gamma}_{11}\cup\hat{\delta}_{11},
\hat{\beta}_{11}\cup\hat{\eta}_{11}\cup\hat{\varepsilon}_{11},
\hat{\varepsilon}_{11}\cup\hat{\gamma}_{11}\cup\hat{\lambda}_{11}$,
and $\hat{\lambda}_{11}\cup\hat{\eta}_{11}\cup\hat{\delta}_{11}$ are
the fixed data of some action in $\mathcal{A}^k_{3\cdot
2^{k-3}}(4)$.
\end{cor}
\begin{proof}
 We proceed by induction on $k$. When $k=3$, we have that
$\hat{\beta}_{11}=\{\beta_1\}, \hat{\gamma}_{11}=\{\gamma\},
\hat{\delta}_{11}=\{\delta\},
\hat{\varepsilon}_{11}=\{\gamma+\beta_1\},
\hat{\eta}_{11}=\{\delta+\beta_1\},
\hat{\lambda}_{11}=\{\gamma+\delta\}$.   Then it is easy to see that
$\hat{\beta}_{11}\cup\hat{\gamma}_{11}\cup\hat{\delta}_{11}=\{\beta_1,
\gamma, \delta\},
\hat{\beta}_{11}\cup\hat{\eta}_{11}\cup\hat{\varepsilon}_{11}=\{\beta_1,
\delta+\beta_1, \gamma+\beta_1\},
\hat{\varepsilon}_{11}\cup\hat{\gamma}_{11}\cup\hat{\lambda}_{11}=\{\gamma+\beta_1,
\gamma, \gamma+\delta\}$, and
$\hat{\lambda}_{11}\cup\hat{\eta}_{11}\cup\hat{\delta}_{11}=\{\gamma+\delta,
\delta+\beta_1, \delta\}$ form  the fixed data of
$(\sigma\phi_3,{\Bbb R}P^3)$ for some $\sigma\in \text{\rm
GL}(3,{\Bbb Z}_2)$. When $k=l\geq 3$, suppose inductively that
$\hat{\beta}_{11}\cup\hat{\gamma}_{11}\cup\hat{\delta}_{11},
\hat{\beta}_{11}\cup\hat{\eta}_{11}\cup\hat{\varepsilon}_{11},
\hat{\varepsilon}_{11}\cup\hat{\gamma}_{11}\cup\hat{\lambda}_{11}$,
and $\hat{\lambda}_{11}\cup\hat{\eta}_{11}\cup\hat{\delta}_{11}$ are
the fixed data of some action in $\mathcal{A}^k_{3\cdot
2^{k-3}}(4)$.

\vskip .2cm

When $k=l+1$, in a similar way to the proof of Claim in
Section~\ref{thm 1.4}, let $\hat{\beta'}_{11}$ denote the set of all
odd sums formed by  $\beta_1,...,\beta_{l-2}$. Then
$\hat{\beta'}_{11}\subset\hat{\beta}_{11}$  contains $2^{l-3}$
different elements.  Let $\hat{\beta''}_{11}$ denote the set formed
by all elements $\beta_1+\beta_{l-1}+\beta, \beta\in
\hat{\beta'}_{11}$. Then,
$\hat{\beta''}_{11}\subset\hat{\beta}_{11}$  contains $2^{l-3}$
different elements, too. Moreover,  one has that
$\hat{\beta}_{11}=\hat{\beta'}_{11}\cup \hat{\beta''}_{11}$ and
$\hat{\beta'}_{11}\cap \hat{\beta''}_{11}=\emptyset$. Now let
$$\begin{cases}
\hat{\gamma'}_{11}=\{\gamma+\beta_1+\beta|\beta\in\hat{\beta'}_{11}\}\\
\hat{\delta'}_{11}=\{\delta+\beta_1+\beta|\beta\in\hat{\beta'}_{11}\}\\
\hat{\varepsilon'}_{11}=\{\gamma+\beta|\beta\in\hat{\beta'}_{11}\}\\
\hat{\eta'}_{11}=\{\delta+\beta|\beta\in\hat{\beta'}_{11}\}\\
\hat{\lambda'}_{11}=\{\gamma+\delta+\beta_1+\beta|\beta\in\hat{\beta'}_{11}\}\\
\end{cases}
\text{and }\
\begin{cases}
\hat{\gamma''}_{11}=\{\gamma+\beta_1+\beta|\beta\in\hat{\beta''}_{11}\}\\
\hat{\delta''}_{11}=\{\delta+\beta_1+\beta|\beta\in\hat{\beta''}_{11}\}\\
\hat{\varepsilon''}_{11}=\{\gamma+\beta|\beta\in\hat{\beta''}_{11}\}\\
\hat{\eta''}_{11}=\{\delta+\beta|\beta\in\hat{\beta''}_{11}\}\\
\hat{\lambda''}_{11}=\{\gamma+\delta+\beta_1+\beta|\beta\in\hat{\beta''}_{11}\}.\\
\end{cases}$$
Then
$$\begin{cases}
\hat{\gamma}_{11}=\hat{\gamma'}_{11}\cup \hat{\gamma''}_{11}\text{ with } \hat{\gamma'}_{11}\cap \hat{\gamma''}_{11}=\emptyset\\
\hat{\delta}_{11}=\hat{\delta'}_{11}\cup \hat{\delta''}_{11}\text{ with } \hat{\delta'}_{11}\cap \hat{\delta''}_{11}=\emptyset\\
\hat{\varepsilon}_{11}=\hat{\varepsilon'}_{11}\cup \hat{\varepsilon''}_{11}\text{ with } \hat{\varepsilon'}_{11}\cap \hat{\varepsilon''}_{11}=\emptyset\\
\hat{\eta}_{11}=\hat{\eta'}_{11}\cup \hat{\eta''}_{11}\text{ with } \hat{\eta'}_{11}\cap \hat{\eta''}_{11}=\emptyset\\
\hat{\lambda}_{11}=\hat{\lambda'}_{11}\cup \hat{\lambda''}_{11}\text{ with } \hat{\lambda'}_{11}\cap \hat{\lambda''}_{11}=\emptyset.\\
\end{cases}
$$
 We see that
$\hat{\beta'}_{11}, \hat{\gamma'}_{11}, \hat{\delta'}_{11},
\hat{\varepsilon'}_{11}, \hat{\eta'}_{11}, \hat{\lambda'}_{11}$
  are exactly formed by
$\beta_1,...,\beta_{l-2},\gamma, \delta$.  Now, regarding
$\{\beta_1,...,\beta_{l-2},\gamma, \delta\}$ as a basis of
$\Hom(({\Bbb Z}_2)^{l},{\Bbb Z}_2)$,  one has by induction that
$\hat{\beta'}_{11}\cup\hat{\gamma'}_{11}\cup\hat{\delta'}_{11},
\hat{\beta'}_{11}\cup\hat{\eta'}_{11}\cup\hat{\varepsilon'}_{11},
\hat{\varepsilon'}_{11}\cup\hat{\gamma'}_{11}\cup\hat{\lambda'}_{11}$,
and $\hat{\lambda'}_{11}\cup\hat{\eta'}_{11}\cup\hat{\delta'}_{11}$
are the fixed data of some $({\Bbb Z}_2)^l$-action, denoted by
$(\Psi, N)$. Then by applying $\Omega$-operation to $(\Psi, N)$, as
in the proof of Lemma~\ref{twist},
 the fixed
data of $\Omega(\Psi,N)$ exactly consists of
$$\hat{\beta}_{11}\cup\hat{\gamma}_{11}\cup\hat{\delta}_{11},\ \
\hat{\beta}_{11}\cup\hat{\eta}_{11}\cup\hat{\varepsilon}_{11},\ \
\hat{\varepsilon}_{11}\cup\hat{\gamma}_{11}\cup\hat{\lambda}_{11},\
\ \hat{\lambda}_{11}\cup\hat{\eta}_{11}\cup\hat{\delta}_{11}.$$ This
completes the induction and the proof.
\end{proof}
\vskip .2cm

\begin{cor} \label{lm11}
$m$ is a power of 2.
\end{cor}

\begin{proof}
 Using the formula of Theorem~\ref{dks}, we consider 1 as a symmetric polynomial in
$n=3m\cdot2^{k-3}+\vert\underline{\hat{\omega}}\vert$ variables.
Then we have
\begin{eqnarray*}
&&{1\over{(\prod\underline{\hat{\omega}})(\prod\hat{\beta}_{11})^m(\prod\hat{\gamma}_{11})^m
(\prod\hat{\delta}_{11})^m}}+
{1\over{(\prod\underline{\hat{\omega}})(\prod\hat{\beta}_{11})^m(\prod\hat{\eta}_{11})^m
(\prod\hat{\varepsilon}_{11})^m}}\\
&&+{1\over{(\prod\underline{\hat{\omega}})(\prod\hat{\varepsilon}_{11})^m(\prod\hat{\gamma}_{11})^m
(\prod\hat{\lambda}_{11})^m}}+
{1\over{(\prod\underline{\hat{\omega}})(\prod\hat{\lambda}_{11})^m(\prod\hat{\eta}_{11})^m
(\prod\hat{\delta}_{11})^m}}\\
 &=&
{1\over{(\prod\underline{\hat{\omega}})\big[(\prod\hat{\beta}_{11})(\prod\hat{\gamma}_{11})
(\prod\hat{\delta}_{11})(\prod\hat{\varepsilon}_{11})(\prod\hat{\eta}_{11})
(\prod\hat{\lambda}_{11})\big]^m}}\\
&& \times
\Bigg\{(\prod\hat{\varepsilon}_{11})^m(\prod\hat{\eta}_{11})^m
(\prod\hat{\lambda}_{11})^m+(\prod\hat{\gamma}_{11})^m(\prod\hat{\delta}_{11})^m
(\prod\hat{\lambda}_{11})^m\\
&&+(\prod\hat{\beta}_{11})^m(\prod\hat{\delta}_{11})^m
(\prod\hat{\eta}_{11})^m+(\prod\hat{\beta}_{11})^m(\prod\hat{\gamma}_{11})^m
(\prod\hat{\varepsilon}_{11})^m\Bigg\}
\end{eqnarray*}
must belong to the polynomial algebra ${\Bbb
Z}_2[\rho_1,...,\rho_{k}]$. Because the degree of the numerator is
smaller than the degree of the denominator, this means the numerator
must be zero, i.e.,
\begin{eqnarray*}&&(\prod\hat{\varepsilon}_{11})^m(\prod\hat{\eta}_{11})^m
(\prod\hat{\lambda}_{11})^m
+(\prod\hat{\gamma}_{11})^m(\prod\hat{\delta}_{11})^m
(\prod\hat{\lambda}_{11})^m\\
&&+(\prod\hat{\beta}_{11})^m(\prod\hat{\delta}_{11})^m
(\prod\hat{\eta}_{11})^m
+(\prod\hat{\beta}_{11})^m(\prod\hat{\gamma}_{11})^m
(\prod\hat{\varepsilon}_{11})^m=0\\
\end{eqnarray*}
By Corollary~\ref{data2}, the above equation still holds for $m=1$,
so one has
 \begin{eqnarray*}
 &&(\prod\hat{\varepsilon}_{11})(\prod\hat{\eta}_{11})
(\prod\hat{\lambda}_{11})
+(\prod\hat{\gamma}_{11})(\prod\hat{\delta}_{11})
(\prod\hat{\lambda}_{11})\\
&&+(\prod\hat{\beta}_{11})(\prod\hat{\delta}_{11})
(\prod\hat{\eta}_{11})
+(\prod\hat{\beta}_{11})(\prod\hat{\gamma}_{11})
(\prod\hat{\varepsilon}_{11}) =0\end{eqnarray*}
  and then
 \begin{eqnarray*}
(*)&&\big[(\prod\hat{\gamma}_{11})(\prod\hat{\delta}_{11})
(\prod\hat{\lambda}_{11})
+(\prod\hat{\beta}_{11})(\prod\hat{\delta}_{11})
(\prod\hat{\eta}_{11})\\
&&+(\prod\hat{\beta}_{11})(\prod\hat{\gamma}_{11})
(\prod\hat{\varepsilon}_{11})\big]^m
+(\prod\hat{\gamma}_{11})^m(\prod\hat{\delta}_{11})^m
(\prod\hat{\lambda}_{11})^m\\
&&+(\prod\hat{\beta}_{11})^m(\prod\hat{\delta}_{11})^m
(\prod\hat{\eta}_{11})^m
+(\prod\hat{\beta}_{11})^m(\prod\hat{\gamma}_{11})^m
(\prod\hat{\varepsilon}_{11})^m\\
&&=0.
 \end{eqnarray*}
 Next, we are going to to show that if $m$ is not a power of 2, then
 $(*)$ does not hold. Let
 $m=2^{p_r}+2^{p_{r-1}}+\cdots+2^{p_1}=j+2^{p_1}$ where
 $p_r>p_{r-1}>\cdots>p_1$. Using the 2-adic expansion of $m$ and
 write
 \begin{eqnarray*}
&&\big[(\prod\hat{\gamma}_{11})(\prod\hat{\delta}_{11})
(\prod\hat{\lambda}_{11})
+(\prod\hat{\beta}_{11})(\prod\hat{\delta}_{11})
(\prod\hat{\eta}_{11})\\
&&+(\prod\hat{\beta}_{11})(\prod\hat{\gamma}_{11})
(\prod\hat{\varepsilon}_{11})\big]^m\\
&=&\big[(\prod\hat{\gamma}_{11})(\prod\hat{\delta}_{11})
(\prod\hat{\lambda}_{11})
+(\prod\hat{\beta}_{11})(\prod\hat{\delta}_{11})
(\prod\hat{\eta}_{11})\\
&&+(\prod\hat{\beta}_{11})(\prod\hat{\gamma}_{11})
(\prod\hat{\varepsilon}_{11})\big]^j\times\Big\{
\big[(\prod\hat{\gamma}_{11})(\prod\hat{\delta}_{11})
(\prod\hat{\lambda}_{11})\big]^{2^{p_1}}\\
&&+\big[(\prod\hat{\beta}_{11})(\prod\hat{\delta}_{11})
(\prod\hat{\eta}_{11})\big]^{2^{p_1}}
+\big[(\prod\hat{\beta}_{11})(\prod\hat{\gamma}_{11})
(\prod\hat{\varepsilon}_{11})\big]^{2^{p_1}}\Big\}.
 \end{eqnarray*}

 \noindent In the above equation we seek the terms of largest degree in
 $\gamma$ and $\delta$. For this we have
 $$\begin{cases}
\prod\hat{\gamma}_{11}=\gamma^{2^{k-3}}+\text{ terms of lower
degree }\\
\prod\hat{\delta}_{11}=\delta^{2^{k-3}}+\text{ terms of lower
degree }\\
\prod\hat{\varepsilon}_{11}=\gamma^{2^{k-3}}+\text{ terms of lower
degree }\\
\prod\hat{\eta}_{11}=\delta^{2^{k-3}}+\text{ terms of lower
degree }\\
\prod\hat{\lambda}_{11}=(\gamma+\delta)^{2^{k-3}}+\text{ terms of
lower
degree }\\
 \end{cases}$$
 and $\prod\hat{\beta}_{11}$ has no $\gamma$'s and $\delta$'s. In
 the $j$-th power, the term of largest degree in $\gamma$ and
 $\delta$ is $\big[\gamma\delta(\gamma+\delta)\big]^{j\cdot
 2^{k-3}}$ occurring in the monomial $\big[(\prod\hat{\gamma}_{11})(\prod\hat{\delta}_{11})
(\prod\hat{\lambda}_{11})\big]^j$. In $(*)$ this monomial is
multiplied by
\begin{eqnarray*}
&&\big[(\prod\hat{\beta}_{11})(\prod\hat{\delta}_{11})
(\prod\hat{\eta}_{11})\big]^{2^{p_1}}
+\big[(\prod\hat{\beta}_{11})(\prod\hat{\gamma}_{11})
(\prod\hat{\varepsilon}_{11})\big]^{2^{p_1}}\\
&=&(\prod\hat{\beta}_{11})^{2^{p_1}}\big\{\gamma^{2^{k-2+p_1}}+\delta^{2^{k-2+p_1}}+\text{
terms of lower degree} \big\}.
\end{eqnarray*}
Thus, in $(*)$ the term with largest degree in $\gamma$ and $\delta$
is
$$(\prod\hat{\beta}_{11})^{2^{p_1}}\gamma^{j\cdot2^{k-3}}\delta^{j\cdot2^{k-3}}(\gamma+\delta)^{
j\cdot 2^{k-3}+2^{k-2+p_1}}$$ which is nonzero. But this is
impossible, so $m$ must be a power of 2.
\end{proof}

Throughout the following, assume that $m$ is always a power of 2.

\vskip .2cm

 By  Corollary~\ref{data2}, let $\hat{\beta}_{11}\cup\hat{\gamma}_{11}\cup\hat{\delta}_{11},
\hat{\beta}_{11}\cup\hat{\eta}_{11}\cup\hat{\varepsilon}_{11},
\hat{\varepsilon}_{11}\cup\hat{\gamma}_{11}\cup\hat{\lambda}_{11}$,
and $\hat{\lambda}_{11}\cup\hat{\eta}_{11}\cup\hat{\delta}_{11}$ be
the fixed data of some action $(\Psi, N)$. Then,  applying
$\Delta$-operation $m$ times to $(\Psi, N)$ gives an action
$\Delta^m(\Psi, N)$. Since $m$ is a power of 2, by Lemma~\ref{dia}
$\Delta^m(\Psi, N)$ is equivariantly cobordant to an action whose fixed data exactly consists of
$\hat{\beta}_{1}\cup\hat{\gamma}_{1}\cup\hat{\delta}_{1},
\hat{\beta}_{1}\cup\hat{\eta}_{1}\cup\hat{\varepsilon}_{1},
\hat{\varepsilon}_{1}\cup\hat{\gamma}_{1}\cup\hat{\lambda}_{1}$, and
$\hat{\lambda}_{1}\cup\hat{\eta}_{1}\cup\hat{\delta}_{1}$. This
gives

\begin{cor}\label{data3}
$\hat{\beta}_{1}\cup\hat{\gamma}_{1}\cup\hat{\delta}_{1},
\hat{\beta}_{1}\cup\hat{\eta}_{1}\cup\hat{\varepsilon}_{1},
\hat{\varepsilon}_{1}\cup\hat{\gamma}_{1}\cup\hat{\lambda}_{1}$, and
$\hat{\lambda}_{1}\cup\hat{\eta}_{1}\cup\hat{\delta}_{1}$ are the
fixed data of some action in $\mathcal{A}^k_{3m\cdot 2^{k-3}}(4)$.
\end{cor}

Now let us further analyze the structure of the changeable part
$\underline{\hat{\omega}}$.

\vskip .2cm Let $f(x_1,...,x_{3m\cdot 2^{k-3}})$ be a symmetric
function  in $3m\cdot 2^{k-3}$ variables over ${\Bbb Z}_2$ where $m$
is a power of 2. Write
 $$\hat{f}={{f(\hat{\beta}_1,\hat{\gamma}_1,\hat{\delta}_1)}\over
{\prod
\hat{\beta}_1\prod\hat{\gamma}_1\prod\hat{\delta}_1}}+{{f(\hat{\beta}_1,\hat{\varepsilon}_1,\hat{\eta}_1)}\over
{\prod \hat{\beta}_1\prod\hat{\varepsilon}_1\prod\hat{\eta}_1}}
+{{f(\hat{\varepsilon}_1,\hat{\gamma}_1,\hat{\lambda}_1)}\over
{\prod \hat{\varepsilon}_1\prod\hat{\gamma}_1\prod\hat{\lambda}_1}}
+{{f(\hat{\lambda}_1,\hat{\eta}_1,\hat{\delta}_1)}\over {\prod
\hat{\lambda}_1\prod\hat{\eta}_1\prod\hat{\delta}_1}}$$  Then
$\hat{f}$ has degree $\deg f-3m\cdot 2^{k-3}$. By
Corollary~\ref{data3} and Theorem~\ref{dks}, it follows that
\begin{lem} $\hat{f}\in{\Bbb Z}_2[\rho_1,...,\rho_{k}]$ is a polynomial. In particular, if $\deg f<3m\cdot 2^{k-3}$, then
$\hat{f}=0$.
\end{lem}
\begin{rem}
 $\hat{f}$ will play  an important role in  determining the structure
of
    $\underline{\hat{\omega}}$.
\end{rem}

\begin{lem}\label{omega}
Let $\underline{\hat{\omega}}$ be nonempty. Then  $\hat{f}$ is
divisible by $\prod\underline{\hat{\omega}}$.
\end{lem}

\begin{proof} Since $\underline{\hat{\omega}}$ is nonempty, we have
that $n=3m\cdot 2^{k-3}+|\underline{\hat{\omega}}|$.  Take the following polynomial function
over ${\Bbb Z}_2$ which is symmetric in variables $x_1,...,x_n$
\begin{eqnarray*}
&&g(x_1,...,x_n;\hat{\beta}_{11}, \hat{\gamma}_{11},
\hat{\delta}_{11}, \hat{\varepsilon}_{11}, \hat{\eta}_{11},
\hat{\lambda}_{11})\\
&=&\sum_{i_1,...,i_n}\Big[h(x_{i_1},...,x_{i_{|\underline{\hat{\omega}}|}};\hat{\beta}_{11},
\hat{\gamma}_{11}, \hat{\delta}_{11}, \hat{\varepsilon}_{11},
\hat{\eta}_{11},
\hat{\lambda}_{11})f(x_{i_{|\underline{\hat{\omega}}|+1}},...,x_{i_{|\underline{\hat{\omega}}|+3m\cdot
2^{k-3}}})\Big]
\end{eqnarray*}
where each $\{i_1,...,i_n\}$ is a permutation of $\{1,...,n\}$ and
\begin{eqnarray*}
&&h(x_{i_1},...,x_{i_{|\underline{\hat{\omega}}|}};\hat{\beta}_{11},
\hat{\gamma}_{11}, \hat{\delta}_{11}, \hat{\varepsilon}_{11},
\hat{\eta}_{11},
\hat{\lambda}_{11})\\
&=&
\prod_{j=1}^{|\underline{\hat{\omega}}|}\Big[\prod(x_{i_j}+\hat{\beta}_{11})\prod(x_{i_j}+\hat{\gamma}_{11})\prod(x_{i_j}+\hat{\delta}_{11})
\prod(x_{i_j}+\hat{\varepsilon}_{11})
\prod(x_{i_j}+\hat{\eta}_{11})\\
&&\times \prod(x_{i_j}+\hat{\lambda}_{11})\Big]
\end{eqnarray*}
and $x_{i_j}+\hat{\beta}_{11}=\{x_{i_j}+\beta\vert
\beta\in\hat{\beta}_{11}\}$ (similarly for
$x_{i_j}+\hat{\gamma}_{11}, x_{i_j}+\hat{\delta}_{11},
x_{i_j}+\hat{\varepsilon}_{11}, x_{i_j}+\hat{\eta}_{11},
x_{i_j}+\hat{\lambda}_{11}$). Using the
formula of Theorem~\ref{dks}, we then have that
\begin{eqnarray*}
\hat{g}&=&{{g(\alpha(E_p);\hat{\beta}_{11}, \hat{\gamma}_{11},
\hat{\delta}_{11}, \hat{\varepsilon}_{11}, \hat{\eta}_{11},
\hat{\lambda}_{11})}\over{\prod\underline{\hat{\omega}}
\prod\hat{\beta}_1\prod\hat{\gamma}_1 \prod\hat{\delta}_1}}+
{{g(\alpha(E_q);\hat{\beta}_{11}, \hat{\gamma}_{11},
\hat{\delta}_{11}, \hat{\varepsilon}_{11}, \hat{\eta}_{11},
\hat{\lambda}_{11})}\over{
\prod\underline{\hat{\omega}}\prod\hat{\beta}_1\prod\hat{\eta}_1
\prod\hat{\varepsilon}_1}}\\
&&+{{g(\alpha(E_r);\hat{\beta}_{11}, \hat{\gamma}_{11},
\hat{\delta}_{11}, \hat{\varepsilon}_{11}, \hat{\eta}_{11},
\hat{\lambda}_{11})}\over{\prod\underline{\hat{\omega}}
\prod\hat{\varepsilon}_1\prod\hat{\gamma}_1
\prod\hat{\lambda}_1}}+{{g(\alpha(E_s);\hat{\beta}_{11},
\hat{\gamma}_{11}, \hat{\delta}_{11}, \hat{\varepsilon}_{11},
\hat{\eta}_{11},
\hat{\lambda}_{11})}\over{\prod\underline{\hat{\omega}}\prod\hat{\lambda}_1\prod\hat{\eta}_1
\prod\hat{\delta}_1}}\\
&=& {{h(\underline{\hat{\omega}};\hat{\beta}_{11},
\hat{\gamma}_{11}, \hat{\delta}_{11}, \hat{\varepsilon}_{11},
\hat{\eta}_{11},
\hat{\lambda}_{11})f(\hat{\beta}_1,\hat{\gamma}_1,\hat{\delta}_1)}\over{\prod\underline{\hat{\omega}}
\prod\hat{\beta}_1\prod\hat{\gamma}_1 \prod\hat{\delta}_1}}\\
&&+ {{h(\underline{\hat{\omega}};\hat{\beta}_{11},
\hat{\gamma}_{11}, \hat{\delta}_{11}, \hat{\varepsilon}_{11},
\hat{\eta}_{11},
\hat{\lambda}_{11})f(\hat{\beta}_1,\hat{\eta}_1,\hat{\varepsilon}_1)}\over{
\prod\underline{\hat{\omega}}\prod\hat{\beta}_1\prod\hat{\eta}_1
\prod\hat{\varepsilon}_1}}\\
&&+{{h(\underline{\hat{\omega}};\hat{\beta}_{11}, \hat{\gamma}_{11},
\hat{\delta}_{11}, \hat{\varepsilon}_{11}, \hat{\eta}_{11},
\hat{\lambda}_{11})f(\hat{\varepsilon}_1,\hat{\gamma}_1,\hat{\lambda}_1)}\over{\prod\underline{\hat{\omega}}
\prod\hat{\varepsilon}_1\prod\hat{\gamma}_1
\prod\hat{\lambda}_1}}\\
&&+{{h(\underline{\hat{\omega}};\hat{\beta}_{11}, \hat{\gamma}_{11},
\hat{\delta}_{11}, \hat{\varepsilon}_{11}, \hat{\eta}_{11},
\hat{\lambda}_{11})f(\hat{\lambda}_1,\hat{\eta}_1,\hat{\delta}_1)}\over{\prod\underline{\hat{\omega}}\prod\hat{\lambda}_1\prod\hat{\eta}_1
\prod\hat{\delta}_1}}\\
&=& {{h(\underline{\hat{\omega}};\hat{\beta}_{11},
\hat{\gamma}_{11}, \hat{\delta}_{11}, \hat{\varepsilon}_{11},
\hat{\eta}_{11},
\hat{\lambda}_{11})\hat{f}}\over{\prod\underline{\hat{\omega}}}}
\end{eqnarray*}
which belongs to ${\Bbb Z}_2[\rho_1,...,\rho_{k}]$. We know from
Proposition~\ref{lm10} that each element of
$\underline{\hat{\omega}}$ has the form $\beta+\gamma+\delta$,
$\beta\in \hat{\beta}_{1}$. An easy argument shows that for any
$\beta\in \hat{\beta}_{1}$
\begin{eqnarray*}
\beta+\gamma+\delta+\hat{\beta}_{11}=\hat{\lambda}_{11}, &
\beta+\gamma+\delta+\hat{\gamma}_{11}=\hat{\eta}_{11}, &
\beta+\gamma+\delta+\hat{\delta}_{11}=\hat{\varepsilon}_{11},\\
\beta+\gamma+\delta+\hat{\varepsilon}_{11}=\hat{\delta}_{11}, &
\beta+\gamma+\delta+\hat{\eta}_{11} =\hat{\gamma}_{11}, &
\beta+\gamma+\delta+\hat{\lambda}_{11}=\hat{\beta}_{11}\end{eqnarray*}
so
\begin{eqnarray*} &&h(\underline{\hat{\omega}};\hat{\beta}_{11},
\hat{\gamma}_{11}, \hat{\delta}_{11}, \hat{\varepsilon}_{11},
\hat{\eta}_{11},
\hat{\lambda}_{11})\\
&=&\Big[(\prod\hat{\lambda}_{11})(\prod\hat{\eta}_{11})(\prod\hat{\varepsilon}_{11})
(\prod\hat{\delta}_{11})(\prod\hat{\gamma}_{11})(\prod\hat{\beta}_{11})\Big]^{|\underline{\hat{\omega}}|}.\end{eqnarray*}
Obviously, $h(\underline{\hat{\omega}};\hat{\beta}_{11},
\hat{\gamma}_{11}, \hat{\delta}_{11}, \hat{\varepsilon}_{11},
\hat{\eta}_{11}, \hat{\lambda}_{11})$ is not divisible by any
$\beta+\gamma+\delta, \beta\in \hat{\beta}_{1}$, so it is not
divisible by $\prod\underline{\hat{\omega}}$. Thus we must have
that $\hat{f}$ is divisible by $\prod\underline{\hat{\omega}}$.
\end{proof}

\begin{lem}\label{lm13}
 $\hat{f}$ is divisible by
$\big[\prod(\gamma+\delta+\hat{\beta}_1)\big]^2=\prod_{\beta\in\hat{\beta}_{11}}(\beta+\gamma+\delta)^{2m}$.
\end{lem}

\begin{proof}
Applying $\Omega$-operation $k-3$ times to the $({\Bbb
Z}_2)^3$-action $(\Lambda, M^5)$ of Lemma~\ref{5-dim} in
Section~\ref{s6} gives a $({\Bbb Z}_2)^k$-action
$\Omega^{k-3}(\Lambda, M^5)$. It is easy to see that
$\Omega^{k-3}(\Lambda, M^5)$ is a
$\underline{\hat{\omega}}$-nonempty action  in
$\mathcal{A}^k_{5\cdot 2^{k-3}}(4)$, which has the property that
$\vert\underline{\hat{\omega}}\vert=2^{k-2}$ and  each element of
$\underline{\hat{\omega}}$ occurs exactly two times.  In particular,
by applying an automorphism $\sigma\in \text{\rm GL}(k,{\Bbb Z}_2)$
of $({\Bbb Z}_2)^k$ to $\Omega^{k-3}(\Lambda, M^5)$ (if necessary),
we can choose the same basis as that in Proposition~\ref{lm10}, so
that $\underline{\hat{\omega}}$ is exactly the disjoint union
$$\{\beta+\gamma+\delta\ \vert
\beta\in\hat{\beta}_{11}\}\cup\{\beta+\gamma+\delta\ \vert
\beta\in\hat{\beta}_{11}\}.$$

Now
 by applying $\Delta$-operation $m$ times to  $\sigma\Omega^{k-3}(\Lambda, M^5)$, one may
obtain a $\underline{\hat{\omega}}$-nonempty $({\Bbb Z}_2)^k$-action
$\Delta^m[\sigma\Omega^{k-3}(\Lambda, M^5)]$ in
$\mathcal{A}^k_{5m\cdot 2^{k-3}}(4)$ such that its changeable part
$$\underline{\hat{\omega}}=\{\beta+\gamma+\delta\
\vert \beta\in\hat{\beta}_1\}\cup\{\beta+\gamma+\delta\ \vert
\beta\in\hat{\beta}_1\},$$ i.e., each element of
$\{\beta+\gamma+\delta\ \vert \beta\in\hat{\beta}_{11}\}$ has
multiplicity $2m$ times in $\underline{\hat{\omega}}$.

\vskip .2cm Similarly to the proof of Lemma~\ref{omega},  taking the
following polynomial function over ${\Bbb Z}_2$ which is symmetric
in variables $x_1,...,x_{5m\cdot 2^{k-3}}$
\begin{eqnarray*}
&&g(x_1,...,x_{5m\cdot 2^{k-3}};\hat{\beta}_{11}, \hat{\gamma}_{11},
\hat{\delta}_{11}, \hat{\varepsilon}_{11}, \hat{\eta}_{11},
\hat{\lambda}_{11})\\
&=&\sum_{i_1,...,i_n}\Big[h(x_{i_1},...,x_{i_{m\cdot2^{k-2}}};\hat{\beta}_{11},
\hat{\gamma}_{11}, \hat{\delta}_{11}, \hat{\varepsilon}_{11},
\hat{\eta}_{11},
\hat{\lambda}_{11})f(x_{i_{m\cdot2^{k-2}+1}},...,x_{i_{5m\cdot
2^{k-3}}})\Big]
\end{eqnarray*}
where each $\{i_1,...,i_n\}=\{1,...,n\}$ and
\begin{eqnarray*}
&&h(x_{i_1},...,x_{i_{m\cdot2^{k-2}}};\hat{\beta}_{11},
\hat{\gamma}_{11}, \hat{\delta}_{11}, \hat{\varepsilon}_{11},
\hat{\eta}_{11},
\hat{\lambda}_{11})\\
&=&
\prod_{j=1}^{m\cdot2^{k-2}}\Big[\prod(x_{i_j}+\hat{\beta}_{11})\prod(x_{i_j}+\hat{\gamma}_{11})\prod(x_{i_j}+\hat{\delta}_{11})
\prod(x_{i_j}+\hat{\varepsilon}_{11})
\prod(x_{i_j}+\hat{\eta}_{11})\\
&&\times \prod(x_{i_j}+\hat{\lambda}_{11})\Big].
\end{eqnarray*}
By direct calculations  one  has that
\begin{eqnarray*}
\hat{g} &=&{{h(\underline{\hat{\omega}};\hat{\beta}_{11},
\hat{\gamma}_{11}, \hat{\delta}_{11}, \hat{\varepsilon}_{11},
\hat{\eta}_{11},
\hat{\lambda}_{11})\hat{f}}\over{\prod_{\beta\in\hat{\beta}_{11}}(\beta+\gamma+\delta)^{2m}}}
\end{eqnarray*}
which belongs to ${\Bbb Z}_2[\rho_1,...,\rho_{k}]$ by
Theorem~\ref{dks}, and so  $\hat{f}$ is divisible by
$\prod_{\beta\in\hat{\beta}_{11}}(\beta+\gamma+\delta)^{2m}$.
\end{proof}

\begin{lem}\label{lm14}
Let $f(x_1,...,x_{3m\cdot2^{k-3}})$ be the product
$$\sigma_{m\cdot2^{k-2}}(x_1,...,x_{3m\cdot2^{k-3}})\cdot\sigma_{3m\cdot2^{k-3}}(x_1,...,x_{3m\cdot2^{k-3}})$$
of the  $(m\cdot2^{k-2})$-th elementary symmetric function and the
$(3m\cdot2^{k-3})$-th elementary symmetric function in
$3m\cdot2^{k-3}$ variables. Then
$$\hat{f}=\prod_{\beta\in\hat{\beta}_{11}}(\beta+\gamma+\delta)^{2m}.$$
\end{lem}

\begin{proof}
Since $\deg f=5m\cdot2^{k-3}$, one has that $\deg
\hat{f}=2m\cdot2^{k-3}=m\cdot2^{k-2}$. Thus, in order to prove that
$\hat{f}=\prod_{\beta\in\hat{\beta}_{11}}(\beta+\gamma+\delta)^{2m}$,
by Lemma~\ref{lm13} it suffices to show that $\hat{f}$ is nonzero.
Since $\sigma_{m\cdot2^{k-2}}(x_1,...,x_{3m\cdot2^{k-3}})$ and
$\sigma_{3m\cdot2^{k-3}}(x_1,...,x_{3m\cdot2^{k-3}})$ are elementary
symmetric functions and $m$ is a power of 2, by direct calculations
one has that
\begin{eqnarray*}
\hat{f}&=&{{f(\hat{\beta}_1, \hat{\gamma}_1, \hat{\delta}_1)}\over
{\prod
\hat{\beta}_1\prod\hat{\gamma}_1\prod\hat{\delta}_1}}+{{f(\hat{\beta}_1,
\hat{\eta}_1, \hat{\varepsilon}_1)}\over {\prod
\hat{\beta}_1\prod\hat{\eta}_1\prod\hat{\varepsilon}_1}}+{{f(\hat{\varepsilon}_1,
\hat{\gamma}_1, \hat{\lambda}_1)}\over {\prod
\hat{\varepsilon}_1\prod\hat{\gamma}_1\prod\hat{\lambda}_1}}+{{f(\hat{\lambda}_1,
\hat{\eta}_1, \hat{\delta}_1)}\over {\prod
\hat{\lambda}_1\prod\hat{\eta}_1\prod\hat{\delta}_1}}\\
&=& \sigma_{m\cdot2^{k-2}}(\hat{\beta}_1, \hat{\gamma}_1,
\hat{\delta}_1)+\sigma_{m\cdot2^{k-2}}(\hat{\beta}_1, \hat{\eta}_1,
\hat{\varepsilon}_1)+\sigma_{m\cdot2^{k-2}}(\hat{\varepsilon}_1,
\hat{\gamma}_1,
\hat{\lambda}_1)\\
&&+\sigma_{m\cdot2^{k-2}}(\hat{\lambda}_1,
\hat{\eta}_1, \hat{\delta}_1)\\
&=&\Big[\sigma_{2^{k-2}}(\hat{\beta}_{11}, \hat{\gamma}_{11},
\hat{\delta}_{11})+\sigma_{2^{k-2}}(\hat{\beta}_{11},
\hat{\eta}_{11},
\hat{\varepsilon}_{11})+\sigma_{2^{k-2}}(\hat{\varepsilon}_{11},
\hat{\gamma}_{11}, \hat{\lambda}_{11})\\
&&+\sigma_{2^{k-2}}(\hat{\lambda}_{11}, \hat{\eta}_{11},
\hat{\delta}_{11})\Big]^m.
\end{eqnarray*}
Write $\bar{\sigma}_{2^{k-2}}$ to be
$$\sigma_{2^{k-2}}(\hat{\beta}_{11}, \hat{\gamma}_{11},
\hat{\delta}_{11})+\sigma_{2^{k-2}}(\hat{\beta}_{11},
\hat{\eta}_{11},
\hat{\varepsilon}_{11})+\sigma_{2^{k-2}}(\hat{\varepsilon}_{11},
\hat{\gamma}_{11},
\hat{\lambda}_{11})+\sigma_{2^{k-2}}(\hat{\lambda}_{11},
\hat{\eta}_{11}, \hat{\delta}_{11})$$ so
$\hat{f}=(\bar{\sigma}_{2^{k-2}})^m$. Thus, this may be reduced to
considering the case $m=1$, i.e., it only needs to show that
$\bar{\sigma}_{2^{k-2}}\not=0$ for $k\geq 3$.

\vskip .2cm

We proceed by induction on $k$. When $k=3$, by Remark~\ref{r9} in
Section~\ref{s6} we know that $\bar{\sigma}_{2}$ is nonzero. When
$k=l$, suppose inductively that $\bar{\sigma}_{2^{l-2}}\not=0$. Now
consider the case in which $k=l+1$.

\vskip .2cm

Let $\hat{\beta'}_{11}$ denote the set of all odd sums formed by
$\beta_1,...,\beta_{l-2}$ (note that $l-2=k-3$), so
$\hat{\beta'}_{11}$ exactly contains the half of all elements of
$\hat{\beta}_{11}$ (i.e., $\hat{\beta'}_{11}$ has just $2^{k-4}$
different elements). Write
$\hat{\beta''}_{11}=\{\beta'+\beta_1+\beta_{l-1}\ \vert
\beta'\in\hat{\beta'}_{11}\}$. Then
$\hat{\beta'}_{11}\cap\hat{\beta''}_{11}=\emptyset$ and
$\hat{\beta'}_{11}\cup\hat{\beta''}_{11}=\hat{\beta}_{11}$.
Similarly, let
$$\begin{cases} \hat{\gamma'}_{11}=\{\gamma+\beta_1+\beta'\ \vert
\beta'\in\hat{\beta'}_{11}\}\\
\hat{\delta'}_{11}=\{\delta+\beta_1+\beta'\ \vert
\beta'\in\hat{\beta'}_{11}\}\\
 \hat{\varepsilon'}_{11}=\{\gamma+\beta'\ \vert
\beta'\in\hat{\beta'}_{11}\}\\
\hat{\eta'}_{11}=\{\delta+\beta'\ \vert
\beta'\in\hat{\beta'}_{11}\}\\
\hat{\lambda'}_{11}=\{\gamma+\delta+\beta_1+\beta'\ \vert
\beta'\in\hat{\beta'}_{11}\} \end{cases}$$ so $\hat{\gamma'}_{11}$
(resp. $\hat{\delta'}_{11}$, $\hat{\varepsilon'}_{11}$,
$\hat{\eta'}_{11}$, and $\hat{\lambda'}_{11}$) also contains exactly
the half of all elements of $\hat{\gamma}_{11}$ (resp.
$\hat{\delta}_{11}$, $\hat{\varepsilon}_{11}$, $\hat{\eta}_{11}$,
and $\hat{\lambda}_{11}$). Furthermore,  one has that
$$\begin{cases}
\hat{\gamma''}_{11}=\{\gamma+\beta'+\beta_{l-1}\ \vert
\beta'\in\hat{\beta'}_{11}\} \text{ with }
\hat{\gamma'}_{11}\cap\hat{\gamma''}_{11}=\emptyset \text{ and }
\hat{\gamma'}_{11}\cup\hat{\gamma''}_{11}=\hat{\gamma}_{11}\\
\hat{\delta''}_{11}=\{\delta+\beta'+\beta_{l-1}\ \vert
\beta'\in\hat{\beta'}_{11}\} \text{ with }
\hat{\delta'}_{11}\cap\hat{\delta''}_{11}=\emptyset \text{ and }
\hat{\delta'}_{11}\cup\hat{\delta''}_{11}=\hat{\delta}_{11}\\
\hat{\varepsilon''}_{11}=\{\gamma+\beta'+\beta_1+\beta_{l-1}\ \vert
\beta'\in\hat{\beta'}_{11}\} \text{ with }
\hat{\varepsilon'}_{11}\cap\hat{\varepsilon''}_{11}=\emptyset \text{
and }
\hat{\varepsilon'}_{11}\cup\hat{\varepsilon''}_{11}=\hat{\varepsilon}_{11}\\
\hat{\eta''}_{11}=\{\delta+\beta'+\beta_1+\beta_{l-1}\ \vert
\beta'\in\hat{\beta'}_{11}\} \text{ with }
\hat{\eta'}_{11}\cap\hat{\eta''}_{11}=\emptyset \text{ and }
\hat{\eta'}_{11}\cup\hat{\eta''}_{11}=\hat{\eta}_{11}\\
 \hat{\lambda''}_{11}=\{\gamma+\delta+\beta'+\beta_{l-1}\
\vert \beta'\in\hat{\beta'}_{11}\} \text{ with }
\hat{\lambda'}_{11}\cap\hat{\lambda''}_{11}=\emptyset \text{ and }
\hat{\lambda'}_{11}\cup\hat{\lambda''}_{11}=\hat{\lambda}_{11}.
\end{cases}$$
Obviously, after reducing modulo $\beta_1+\beta_{l-1}$, we see that
$\hat{\beta''}_{11}$ (resp. $\hat{\gamma''}_{11}$,
$\hat{\delta''}_{11}$, $\hat{\varepsilon''}_{11}$,
$\hat{\eta''}_{11}$, and $\hat{\lambda''}_{11}$) becomes
$\hat{\beta'}_{11}$ (resp. $\hat{\gamma'}_{11}$,
$\hat{\delta'}_{11}$, $\hat{\varepsilon'}_{11}$, $\hat{\eta'}_{11}$,
and $\hat{\lambda'}_{11}$). Thus
\begin{eqnarray*} &&\bar{\sigma}_{2^{l-1}}\\&\equiv &
\sigma_{2^{l-1}}(\hat{\beta'}_{11},\hat{\beta'}_{11},
\hat{\gamma'}_{11}, \hat{\gamma'}_{11},\hat{\delta'}_{11},
\hat{\delta'}_{11})+
\sigma_{2^{l-1}}(\hat{\beta'}_{11},\hat{\beta'}_{11},
\hat{\eta'}_{11}, \hat{\eta'}_{11},\hat{\varepsilon'}_{11},
\hat{\varepsilon'}_{11})\\
&&+\sigma_{2^{l-1}}(\hat{\varepsilon'}_{11},\hat{\varepsilon'}_{11},
\hat{\gamma'}_{11}, \hat{\gamma'}_{11},\hat{\lambda'}_{11},
\hat{\lambda'}_{11})+\sigma_{2^{l-1}}(\hat{\lambda'}_{11},\hat{\lambda'}_{11},
\hat{\eta'}_{11}, \hat{\eta'}_{11},\hat{\delta'}_{11},
\hat{\delta'}_{11})\\
&\equiv & \sigma^2_{2^{l-2}}(\hat{\beta'}_{11}, \hat{\gamma'}_{11},
\hat{\delta'}_{11})+ \sigma^2_{2^{l-2}}(\hat{\beta'}_{11},
\hat{\eta'}_{11}, \hat{\varepsilon'}_{11})
+\sigma^2_{2^{l-2}}(\hat{\varepsilon'}_{11},
\hat{\gamma'}_{11}, \hat{\lambda'}_{11})\\
&& +\sigma^2_{2^{l-2}}(\hat{\lambda'}_{11}, \hat{\eta'}_{11},
\hat{\delta'}_{11})\\
&\equiv&\Big\{\sigma_{2^{l-2}}(\hat{\beta'}_{11},
\hat{\gamma'}_{11},  \hat{\delta'}_{11})+
\sigma_{2^{l-2}}(\hat{\beta'}_{11}, \hat{\eta'}_{11},
\hat{\varepsilon'}_{11})
 +\sigma_{2^{l-2}}(\hat{\varepsilon'}_{11},
\hat{\gamma'}_{11}, \hat{\lambda'}_{11})\\
&&+\sigma_{2^{l-2}}(\hat{\lambda'}_{11}, \hat{\eta'}_{11},
\hat{\delta'}_{11})\Big\}^2\\
&\equiv& (\bar{\sigma}_{2^{l-2}})^2\mod (\beta_1+\beta_{l-1})\\
&\not=& 0 \text{ by induction}
 \end{eqnarray*}
 so $\bar{\sigma}_{2^{l-1}}$ must be nonzero. This completes the
induction, and thus $\hat{f}\not=0$.
\end{proof}
Combining  Lemmas~\ref{omega}-\ref{lm14} and Proposition~\ref{lm10},
 it easily follows that

\begin{cor} \label{dim}
Each element of $\{\beta+\gamma+\delta\
\vert\beta\in\hat{\beta}_{11}\}$ occurs at most $2m$ times in
$\underline{\hat{\omega}}$, and so $|\underline{\hat{\omega}}|\leq
m\cdot2^{k-2}$. Furthermore, $n$ is in the range $3m\cdot
2^{k-3}\leq n\leq 5m\cdot 2^{k-3}$.
\end{cor}

Together with all arguments above, now we can give a description of
the essential structure of $\mathcal{N}_{(\Phi, M^n)}$, which is
stated as follows.

\begin{thm}\label{tangent}
Let $(\Phi, M^n)\in \mathcal{A}_n^k(4)$ be an action  with four
fixed points $p,q,r,s$ and let $(\Gamma_{(\Phi, M^n)}, \alpha)$ be a
colored graph of $(\Phi, M^n)$. Then
\begin{enumerate}
\item[$(i)$]
$k\geq 3$ and  $n$ is in the range $3\cdot 2^\ell \leq n\leq 5\cdot
2^\ell$ for some $\ell\geq k-3$;
\item[$(ii)$]  there is a basis $\{\beta_1,...,\beta_{k-2},$ $\gamma,
\delta\}$ of $\Hom(({\Bbb Z}_2)^k,{\Bbb Z}_2)$ such that
\begin{eqnarray*}
\alpha(E_p)=\hat{\beta}\cup\hat{\gamma}\cup\hat{\delta}\cup\hat{\omega},
&& \alpha(E_q)=\hat{\beta}\cup\hat{\eta}\cup
\hat{\varepsilon}\cup\hat{\omega},\\\alpha(E_r)=
 \hat{\gamma}\cup\hat{\varepsilon}\cup\hat{\lambda}\cup\hat{\omega}, &&\alpha(E_s)=
 \hat{\delta}\cup\hat{\eta}\cup\hat{\lambda}\cup\hat{\omega}\end{eqnarray*}
 where $\hat{\beta}$ is a multiset consisting of all sums with  same multiplicity
$2^{\ell-k+3}$ formed by the odd number
 of elements of
$\beta_1,...,\beta_{k-2}$, and
$$\begin{cases}
\hat{\gamma}=\{\gamma+\beta_1+\beta|\beta\in\hat{\beta}\}\\
\hat{\delta}=\{\delta+\beta_1+\beta|\beta\in\hat{\beta}\}\\
\hat{\varepsilon}=\{\gamma+\beta|\beta\in\hat{\beta}\}\\
\hat{\eta}=\{\delta+\beta|\beta\in\hat{\beta}\}\\
\hat{\lambda}=\{\gamma+\delta+\beta_1+\beta|\beta\in\hat{\beta}\}\\
|\hat{\omega}|=n-3\cdot 2^{\ell}
\end{cases}$$
and  every element in $\hat{\omega}$ has the form
$\gamma+\delta+\beta$  and occurs at most $2^{\ell-k+4}$ times if
$\hat{\omega}$ is non-empty, where $\beta\in\hat{\beta}$.
\end{enumerate}
\end{thm}

\subsection{The existence of actions in
$\mathcal{A}_n^k(4)$}\label{existence}

\begin{defn}
Given a basis $B=\{\beta_1,...,\beta_{k-2},\gamma, \delta\}$  of
$\Hom(({\Bbb Z}_2)^k,{\Bbb Z}_2)$ and an integer $\ell\geq k-3$
where $k\geq 3$. We say that $\mathcal{S}_B^{k,\ell}=\{\hat{\beta},
\hat{\gamma}, \hat{\delta},  \hat{\varepsilon}, \hat{\eta},
\hat{\lambda}\}$ is a {\em $\hat{\omega}$-empty
$2^{\ell-k+3}$-multi-structure over $B$} if $\hat{\beta}$ is a
multiset consisting of all odd sums with  same multiplicity
$2^{\ell-k+3}$ formed by $\beta_1,...,\beta_{k-2}$,  and
$$\begin{cases}
\hat{\gamma}=\{\gamma+\beta_1+\beta|\beta\in\hat{\beta}\}\\
\hat{\delta}=\{\delta+\beta_1+\beta|\beta\in\hat{\beta}\}\\
\hat{\varepsilon}=\{\gamma+\beta|\beta\in\hat{\beta}\}\\
\hat{\eta}=\{\delta+\beta|\beta\in\hat{\beta}\}\\
\hat{\lambda}=\{\gamma+\delta+\beta_1+\beta|\beta\in\hat{\beta}\}.
\end{cases}$$
\end{defn}

\begin{rem}
Clearly,  $B$ and $\ell$ determines a unique $\hat{\omega}$-empty
$2^{\ell-k+3}$-multi-structure over $B$. Thus, if $\ell$ is fixed,
then any two $\hat{\omega}$-empty $2^{\ell-k+3}$-multi-structures
$\mathcal{S}_{B_1}^{k,\ell}$ and $\mathcal{S}_{B_2}^{k,\ell}$ can be
translated to each other by automorphisms of $\Hom(({\Bbb
Z}_2)^k,{\Bbb Z}_2)$.
\end{rem}

By Proposition~\ref{lm10}, Corollary~\ref{data3} and
Theorem~\ref{ce}, one has that

\begin{lem}\label{lm}
Given an $\hat{\omega}$-empty $2^{\ell-k+3}$-multi-structure
$\mathcal{S}_B^{k,\ell}=\{\hat{\beta}, \hat{\gamma}, \hat{\delta},
\hat{\varepsilon}, \hat{\eta}, \hat{\lambda}\}$ over
$B=\{\beta_1,...,\beta_{k-2},$ $\gamma, \delta\}$, up to equivariant
cobordism there is a unique $\hat{\omega}$-empty action in
$\mathcal{A}^k_{3\cdot 2^\ell}(4)$ such that its tangent
representation set consisting of
\begin{eqnarray*}
\hat{\beta}\cup\hat{\gamma}\cup\hat{\delta}, \ \
\hat{\beta}\cup\hat{\eta}\cup \hat{\varepsilon}, \ \
 \hat{\gamma}\cup\hat{\varepsilon}\cup\hat{\lambda},
 \ \
 \hat{\delta}\cup\hat{\eta}\cup\hat{\lambda}. \end{eqnarray*}
 Conversely, each action in $\mathcal{A}^k_{3\cdot 2^\ell}(4)$
 determines a unique $\hat{\omega}$-empty
 $2^{\ell-k+3}$-multi-structure.
\end{lem}

\begin{rem}
we easily see from  Proposition~\ref{lm10} that each action of
$\mathcal{A}^k_{3\cdot 2^\ell+t}(4)$ with $t\leq 2^{\ell+1}$ can
still determine a unique $\hat{\omega}$-empty
 $2^{\ell-k+3}$-multi-structure.
\end{rem}

\begin{prop}\label{lm12}
In every dimension $n$ with $3\cdot 2^{\ell}\leq n\leq 5\cdot
2^{\ell}$ for every $\ell\geq k-3$,  $\mathcal{A}_n^k(4)$ is
nonempty.
\end{prop}

\begin{proof}
Given an integer $\ell\geq k-3$, let $n=3\cdot 2^\ell+t$ where
$0\leq t\leq 2^{\ell+1}$. Take an $\hat{\omega}$-empty
$2^{\ell-k+3}$-multi-structure
$\mathcal{S}_B^{k,\ell}=\{\hat{\beta}, \hat{\gamma}, \hat{\delta},
\hat{\varepsilon}, \hat{\eta}, \hat{\lambda}\}$ over
$B=\{\beta_1,...,\beta_{k-2},$ $\gamma, \delta\}$, by Lemma~\ref{lm}
and  Theorem~\ref{dks} one has that for any symmetric function
$f(x_1,...,x_{3\cdot 2^\ell})$ in $3\cdot 2^\ell$ variables over
${\Bbb Z}_2$
  $$\hat{f}={{f(\hat{\beta},\hat{\gamma},\hat{\delta})}\over
{\prod
\hat{\beta}\prod\hat{\gamma}\prod\hat{\delta}}}+{{f(\hat{\beta},\hat{\varepsilon},\hat{\eta})}\over
{\prod \hat{\beta}\prod\hat{\varepsilon}\prod\hat{\eta}}}
+{{f(\hat{\varepsilon},\hat{\gamma},\hat{\lambda})}\over {\prod
\hat{\varepsilon}\prod\hat{\gamma}\prod\hat{\lambda}}}
+{{f(\hat{\lambda},\hat{\eta},\hat{\delta})}\over {\prod
\hat{\lambda}\prod\hat{\eta}\prod\hat{\delta}}}$$ belongs to ${\Bbb
Z}_2[\rho_1, ..., \rho_k]$. Further,  by Lemma~\ref{lm13}, $\hat{f}$
is divisible by $\big[\prod(\gamma+\delta+\hat{\beta})\big]^2$.

\vskip .2cm  Now, for any integer $t\leq 2^{\ell+1}$, one can always
choose a multiset $\Theta$ formed by
$\beta_1,...,\beta_{k-2},\gamma,\delta$ such that $|\Theta|=t$, and
each element of $\Theta$ is chosen in the $2^{k-3}$ different
elements of $\{\beta+\gamma+\delta\ \vert\beta\in\hat{\beta}\}$ and
has multiplicity  at most $2^{\ell-k+4}$. Note that if $t$ is zero
then $\Theta$ is empty. Since $|\hat{\beta}|=2^{\ell}$ and $\ell\geq k-3$, one has that
$\big[\prod(\gamma+\delta+\hat{\beta})\big]^2$ is divisible by
$\prod\Theta$, so  $\hat{f}$ is also divisible by $\prod\Theta$.
Consider any symmetric function $g(x_1,...,x_n)$ in $n$ variables.
One   can  write $g(x_1,...,x_n)$ as a sum
$$\sum_{h,f} h(x_1,...,x_{t})f(x_{t+1},...,x_n)$$ of products
$h(x_1,...,x_{t})f(x_{t+1},...,x_n)$ such that each $h$ is  a
function in $t$ variables $x_1, ..., x_t$ and  each $f$ is always a
symmetric function in $n-t$ variables $x_{t+1}, ..., x_n$ (note that
for the cases $n=4, 5$, see the proofs of
Lemmas~\ref{4-dim}-\ref{5-dim}). Then one has that
\begin{eqnarray*} \hat{g} &=&
{{g(\Theta,\hat{\beta},\hat{\gamma},\hat{\delta})}\over
{\prod\Theta\prod
\hat{\beta}\prod\hat{\gamma}\prod\hat{\delta}}}+{{g(\Theta,\hat{\beta},\hat{\eta},\hat{\varepsilon})}\over
{\prod\Theta\prod
\hat{\beta}\prod\hat{\eta}\prod\hat{\varepsilon}}}+
{{g(\Theta,\hat{\varepsilon},\hat{\gamma},\hat{\lambda})}\over
{\prod\Theta\prod
\hat{\varepsilon}\prod\hat{\gamma}\prod\hat{\lambda}}}\\
&&+ {{g(\Theta,\hat{\lambda},\hat{\eta},\hat{\delta})}\over
{\prod\Theta\prod
\hat{\lambda}\prod\hat{\eta}\prod\hat{\delta}}}\\
&=&\sum{{h(\Theta)}\over{\prod\Theta}}\Bigg\{{{f(\hat{\beta},\hat{\gamma},\hat{\delta})}\over
{\prod
\hat{\beta}\prod\hat{\gamma}\prod\hat{\delta}}}+{{f(\hat{\beta},\hat{\eta},\hat{\varepsilon})}\over
{\prod \hat{\beta}\prod\hat{\eta}\prod\hat{\varepsilon}}}+
{{f(\hat{\varepsilon},\hat{\gamma},\hat{\lambda})}\over {\prod
\hat{\varepsilon}\prod\hat{\gamma}\prod\hat{\lambda}}}+
{{f(\hat{\lambda},\hat{\eta},\hat{\delta})}\over {\prod
\hat{\lambda}\prod\hat{\eta}\prod\hat{\delta}}}\Bigg\}\\
&=& \sum{{h(\Theta)\hat{f}}\over{\prod\Theta}}
\end{eqnarray*}
 so $\hat{g}$ belongs to ${\Bbb Z}_2[\rho_1,...,\rho_{k}]$.
Thus by Theorem~\ref{dks} and Proposition~\ref{lm10}, there is a
$({\Bbb Z}_2)^k$-action such that its changeable part
$\underline{\hat{\omega}}$ is just $\Theta$.
\end{proof}

Together with Corollaries~\ref{lm11}-\ref{dim} and
Proposition~\ref{lm12}, we have completely determined the existence
of all actions in $\mathcal{A}_n^k(4)$. The result is stated as follows.

\begin{thm} \label{exist}
$\mathcal{A}_n^k(4)$ is nonempty if and only if $k\geq 3$ and $n$ is
in the range $\bigcup_{\ell \geq k-3}[3\cdot 2^{\ell}, 5\cdot
2^{\ell}]$.
\end{thm}

We see from the proof of Proposition~\ref{lm12} that
$\underline{\hat{\omega}}$ can always happen as long as $0\leq
|\underline{\hat{\omega}}|\leq 2^{\ell+1}$. However, when we  fix an
integer $t$ with $0\leq t\leq 2^{\ell+1}$, generally there may be
different choices of $\underline{\hat{\omega}}$ with
$|\underline{\hat{\omega}}|=t$ if $t>0$. This can be seen from the
following example.

\begin{example}
Consider the $({\Bbb Z}_2)^4$-manifold $\Omega(\phi_3, {\Bbb R}P^3)$.
This is a $\underline{\hat{\omega}}$-empty $({\Bbb Z}_2)^4$-action,
whose colored graph is shown in Figure~\ref{chang}.
\begin{figure}[h]
    \input{ld4.pstex_t}\centering
    \caption[a]{The  colored graph of $\Omega(\phi_3, {\Bbb R}P^3)$}\label{chang}
\end{figure}
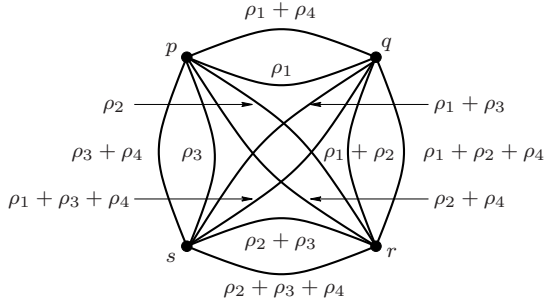
It is easy to see that there are two different choices for
$\underline{\hat{\omega}}$ with $|\underline{\hat{\omega}}|=1$,
which are $\{\rho_1+\rho_2+\rho_3\}$ and
$\{\rho_1+\rho_2+\rho_3+\rho_4\}$ respectively, and there are three different
choices for $\underline{\hat{\omega}}$ with
$|\underline{\hat{\omega}}|=2$, which are $\{\rho_1+\rho_2+\rho_3,
\rho_1+\rho_2+\rho_3\}$, $\{\rho_1+\rho_2+\rho_3,
\rho_1+\rho_2+\rho_3+\rho_4\}$ and $\{\rho_1+\rho_2+\rho_3+\rho_4,
\rho_1+\rho_2+\rho_3+\rho_4\}$ respectively. Also, these
different choices of $\underline{\hat{\omega}}$  with
$|\underline{\hat{\omega}}|=1$ or $2$ determine different actions up to
equivariant cobordism.
\end{example}

\vskip .2cm

 Now fix an $\hat{\omega}$-empty
$2^{\ell-k+3}$-multi-structure
$\mathcal{S}_B^{k,\ell}=\{\hat{\beta}, \hat{\gamma}, \hat{\delta},
\hat{\varepsilon}, \hat{\eta}, \hat{\lambda}\}$ over
$B=\{\beta_1,...,\beta_{k-2},$ $\gamma, \delta\}$ and then  let us
look at how many there are choices of $\underline{\hat{\omega}}$
with $|\underline{\hat{\omega}}|=t$.
 \vskip .2cm

First, write  $2^{k-3}$ different elements in
$\{\beta+\gamma+\delta\ \vert\beta\in\hat{\beta}\}$ as $w_1, ...,
w_{2^{k-3}}$. Next, choose an integer vector ${\bf v}=(v_1,...,
v_{2^{k-3}})$ in ${\Bbb Z}_{\geq 0}^{2^{k-3}}$ such that $$|{\bf
v}|=v_1+\cdots+v_{2^{k-3}}=t$$ and each $v_i\leq 2^{\ell-k+4}$. This
vector ${\bf v}$ determines a unique choice of
$\underline{\hat{\omega}}$ with $|\underline{\hat{\omega}}|=t$,
denoted by $\underline{\hat{\omega}}^{\bf v}_B$,  in such a way that
each $w_i$ has multiplicity  just $v_i$  in
$\underline{\hat{\omega}}$. Thus, it is not difficult to see that
the number of different choices for $\underline{\hat{\omega}}$ with
$|\underline{\hat{\omega}}|=t$ is exactly the number of those
lattices satisfying the equation $x_1+\cdots+x_{2^{k-3}}=t$ with
$0\leq x_i\leq 2^{\ell-k+4}$ in ${\Bbb R}^{2^{k-3}}$, i.e., the
number of all integer vectors in the set
$$\mathcal{I}^k(t):=\Big\{(v_1,..., v_{2^{k-3}})\in {\Bbb Z}_{\geq 0}^{2^{k-3}}\Big|\sum^{2^{k-3}}_{i=1}v_i=t \text{ with each } \
0\leq v_i\leq 2^{\ell-k+4}\Big\}.$$ Obviously, if $k=3$, then
$\underline{\hat{\omega}}$ with $|\underline{\hat{\omega}}|=t$ has a
unique choice.

\vskip .2cm We see from the proof of Proposition~\ref{lm12} that
each $\underline{\hat{\omega}}^{\bf v}_B$ determines a $({\Bbb
Z}_2)^k$-action  with the tangent representation set
$$\{\underline{\hat{\omega}}^{\bf
v}_B\cup\hat{\beta}\cup\hat{\gamma}\cup\hat{\delta},
\underline{\hat{\omega}}^{\bf
v}_B\cup\hat{\beta}\cup\hat{\eta}\cup\hat{\varepsilon},
\underline{\hat{\omega}}^{\bf
v}_B\cup\hat{\varepsilon}\cup\hat{\gamma}\cup\hat{\lambda},
\underline{\hat{\omega}}^{\bf
v}_B\cup\hat{\lambda}\cup\hat{\eta}\cup\hat{\delta}\}, $$ which is
called {\em the $({\Bbb Z}_2)^k$-action of type
$\underline{\hat{\omega}}^{\bf v}_B$} in $\mathcal{A}^k_{3\cdot
2^\ell+t}(4)$.
\begin{lem} \label{changeable}
Let $k>3$ and let ${\bf v}_1$ and ${\bf v}_2$ be two different
vectors in $\mathcal{I}^k(t)$. Then the $({\Bbb Z}_2)^k$-action of
type $\underline{\hat{\omega}}^{{\bf v}_1}_B$ is not equivariantly
cobordant to the $({\Bbb Z}_2)^k$-action of
$\underline{\hat{\omega}}^{{\bf v}_2}_B$.
\end{lem}

\begin{proof}
Obviously, both ${\bf v}_1$ and ${\bf v}_2$ with ${\bf v}_1\not={\bf
v}_2$ give two different $\underline{\hat{\omega}}^{{\bf v}_1}_B$
and $\underline{\hat{\omega}}^{{\bf v}_2}_B$. Then two $({\Bbb
Z}_2)^k$-actions of types $\underline{\hat{\omega}}^{{\bf v}_1}_B$
and $\underline{\hat{\omega}}^{{\bf v}_2}_B$ determine  two
different tangent representation sets. Furthermore, the lemma
follows from this by Theorem~\ref{ce}.
\end{proof}

\begin{rem}
We see by Lemma~\ref{changeable} that up to equivariant cobordism,
the number of all $({\Bbb Z}_2)^k$-actions of types
$\underline{\hat{\omega}}^{\bf v}_B$ with $|{\bf v}|=t$ in
$\mathcal{A}^k_{3\cdot 2^\ell+t}(4)$ is the same as that of elements
in $\mathcal{I}^k(t)$.
\end{rem}

\subsection{The equivariant cobordism classification of all actions in
$\mathcal{A}_n^k(4)$}\label{cl}

Suppose that $\mathcal{A}_n^k(4)$ is nonempty. Then by
Theorem~\ref{exist} there are integers $\ell$ and $t$ with $\ell\geq
k-3$ and $t\leq 2^{\ell+1}$ such that $n=3\cdot 2^\ell+t$.

\vskip .2cm

Now let us consider the equivariant cobordism classification of all
actions in $\mathcal{A}_{3\cdot 2^\ell+t}^k(4)$.

\vskip .2cm

For $t=0$, one knows from Remark~\ref{classify} that for any two
actions $(\Phi_1, M_1)$ and $(\Phi_2, M_2)$ in
$\mathcal{A}^k_{3\cdot 2^\ell}(4)$, there is an automorphism
$\sigma\in \text{\rm GL}(k, {\Bbb Z}_2)$ such that $(\Phi_1, M_1)$
is equivariantly cobordant to $(\sigma\Phi_2, M_2)$. On the other
hand, applying $\Omega$-operation $k-3$ times and $\Delta$-operation
$2^{\ell-k+3}$ times to $(\phi_3, {\Bbb R}P^3)$ gives a $({\Bbb
Z}_2)^k$-action $\Delta^{2^{\ell-k+3}}\Omega^{k-3}(\phi_3, {\Bbb
R}P^3)$ in $\mathcal{A}^k_{3\cdot 2^\ell}(4)$. Thus one has

 \begin{prop}
\label{p5}
 Each of $\mathcal{A}^k_{3\cdot 2^{\ell}}(4)$ with $\ell \geq k-3$
is equivariantly cobordant to one of
$$\sigma\Delta^{2^{\ell-k+3}}\Omega^{k-3}(\phi_3, {\Bbb R}P^3),\
\sigma\in \text{\rm GL}(k, {\Bbb Z}_2).$$
\end{prop}

For $t>0$, as shown in Subsection~\ref{existence}, generally the
changeable part $\underline{\hat{\omega}}$ with
$|\underline{\hat{\omega}}|=t$ may have many possible choices, but
we  can give a description for those possible choices. In
particular, this description also works for the case $t=0$. Thus,
throughout the following assume that $0\leq t\leq 2^{\ell+1}$.

\vskip .2cm

 Beginning with the $({\Bbb Z}_2)^k$-action
$\Delta^{2^{\ell-k+3}}\Omega^{k-3}(\phi_3, {\Bbb R}P^3)$ in
$\mathcal{A}_{3\cdot 2^\ell}^k(4)$, which is equivariantly cobordant
to an action $(\Phi, M)$ fixing exactly four isolated points. Then,
by Lemma~\ref{lm} we can obtain a unique $\hat{\omega}$-empty
$2^{\ell-k+3}$-multi-structure
$\mathcal{S}_B^{k,\ell}=\{\hat{\beta}, \hat{\gamma}, \hat{\delta},
\hat{\varepsilon}, \hat{\eta}, \hat{\lambda}\}$ over
$B=\{\beta_1,...,\beta_{k-2},$ $\gamma, \delta\}$ such that the
tangent representation set of $(\Phi, M)$ consists of
\begin{eqnarray*}
\hat{\beta}\cup\hat{\gamma}\cup\hat{\delta}, \ \
\hat{\beta}\cup\hat{\eta}\cup \hat{\varepsilon}, \ \
 \hat{\gamma}\cup\hat{\varepsilon}\cup\hat{\lambda},
 \ \
 \hat{\delta}\cup\hat{\eta}\cup\hat{\lambda}. \end{eqnarray*}
  Now,  given an integer vector ${\bf v}\in \mathcal{I}^k(t)$ with $|{\bf v}|=t$, in the way as  in Subsection~\ref{exist}, one can obtain a
  changeable part
 $\underline{\hat{\omega}}^{\bf v}_B$. Furthermore, by Lemma~\ref{changeable}, up to equivariant cobordism there is a unique  action denoted by
  $\Lambda^{\bf v}\Delta^{2^{\ell-k+3}}\Omega^{k-3}(\phi_3, {\Bbb R}P^3)$ such that its tangent representation set
  is
$$\big\{\underline{\hat{\omega}}^{\bf
v}_B\cup\hat{\beta}\cup\hat{\gamma}\cup\hat{\delta},
\underline{\hat{\omega}}^{\bf
v}_B\cup\hat{\beta}\cup\hat{\eta}\cup\hat{\varepsilon},
\underline{\hat{\omega}}^{\bf
v}_B\cup\hat{\varepsilon}\cup\hat{\gamma}\cup\hat{\lambda},
\underline{\hat{\omega}}^{\bf
v}_B\cup\hat{\lambda}\cup\hat{\eta}\cup\hat{\delta}\big\}.$$ Let
$\sigma\Lambda^{\bf v}\Delta^{2^{\ell-k+3}}\Omega^{k-3}(\phi_3,
{\Bbb R}P^3)$ denote the action produced by applying an automorphism
$\sigma\in \text{\rm GL}(k,{\Bbb Z}_2)$ to $\Lambda^{\bf
v}\Delta^{2^{\ell-k+3}}\Omega^{k-3}(\phi_3, {\Bbb R}P^3)$.

\begin{thm} \label{classification}
Up to equivariant cobordism,
$$\Big\{\sigma\Lambda^{\bf v}\Delta^{2^{\ell-k+3}}\Omega^{k-3}(\phi_3, {\Bbb
R}P^3)\big|{\bf v}\in \mathcal{I}^k(t), \sigma\in \text{\rm
GL}(k,{\Bbb Z}_2)\Big\}$$ gives all possible nonbounding actions in
$\mathcal{A}_{3\cdot 2^\ell+t}^k(4)$ where $\ell\geq k-3\geq 0$ and
$0\leq t\leq 2^{\ell+1}$.
\end{thm}

\begin{proof}
Let $(\Psi, N)$ be an action in $\mathcal{A}_{3\cdot
2^\ell+t}^k(4)$. Without loss of generality assume that $(\Psi, N)$ has exactly
four fixed points. It is not difficult  to see from
Proposition~\ref{lm10} that $(\Psi, N)$ determines a unique
$\hat{\omega}$-empty $2^{\ell-k+3}$-multi-structure over some basis
$B'$ of $\Hom(({\Bbb Z}_2)^k, {\Bbb Z}_2)$. In addition, $(\Psi, N)$
also determines a unique changeable part
$\underline{\hat{\omega}}^{\bf v}_{B'}$. Thus, the tangent
representation set $\mathcal{N}_{(\Psi, N)}$ is uniquely constructed
by $B'$. Since any two bases in $\Hom(({\Bbb Z}_2)^k, {\Bbb Z}_2)$
can always be translated to each other by automorphisms of
$\Hom(({\Bbb Z}_2)^k, {\Bbb Z}_2)$, there must be one $ \sigma\in
\text{\rm GL}(k,{\Bbb Z}_2)$ such that $(\Psi, N)$ and $\sigma
\Lambda^{\bf v} \Delta^{2^{\ell-k+3}} \Omega^{k-3} (\phi_3, {\Bbb
R}P^3)$ have the same tangent representation set. The theorem then
follows from  Theorem~\ref{ce}.
\end{proof}

\begin{rem}
It should be pointed out that $\Lambda^{\bf
v}\Delta^{2^{\ell-k+3}}\Omega^{k-3}(\phi_3, {\Bbb R}P^3)$ is not a
concrete action, but its tangent representation set is concrete and
can be constructed. We try to construct a concrete action, but fail.
\end{rem}

Now, with Theorems~\ref{exist} and \ref{classification} together, we
complete the proof of Theorem~\ref{thm4}.

\subsection{The characterization of the colored graphs of actions in
$\mathcal{A}_n^k(4)$}\label{char-co}

Suppose that $\mathcal{A}_n^k(4)$ is nonempty. Then one knows from
Theorem~\ref{exist} there are integers $\ell$ and $t$ with $\ell\geq
k-3$ and $t\leq 2^{\ell+1}$ such that $n=3\cdot 2^\ell+t$.

 \vskip .2cm

 Now let $(\Gamma, \alpha)$ be an abstract 1-skeleton of type
$(3\cdot 2^\ell+t,k)$ such that $\Gamma$ contains exactly four
vertices $p,q,r,s$. Then we consider the following question.
\begin{enumerate}
\item[({\bf P})]
{\em When can $(\Gamma, \alpha)$  become a colored graph of some
action in $\mathcal{A}_{3\cdot 2^\ell+t}^k(4)$?}
\end{enumerate}

If $t=0$,  we can characterize the colored graphs of actions in
$\mathcal{A}_{3\cdot 2^\ell+t}^k(4)$, and our result is stated as
follows.

\begin{thm}\label{abstr-graph}
If $t=0$, then $(\Gamma, \alpha)$ is realizable as a colored graph
of some action in $\mathcal{A}_{3\cdot 2^\ell}^k(4)$ if and only if
there is a basis $\{\beta_1,...,\beta_{k-2}, \gamma, \delta\}$ of
$\Hom({\Bbb Z}_2)^k,{\Bbb Z}_2)$  such that

\begin{eqnarray*}
&\alpha(E_p)=\hat{\beta}\cup\hat{\gamma}\cup\hat{\delta}, &
\alpha(E_q)=\hat{\beta}\cup\hat{\eta}\cup\hat{\varepsilon},\\
&\alpha(E_r)=\hat{\varepsilon}\cup\hat{\gamma}\cup\hat{\lambda}, &
\alpha(E_s)=\hat{\lambda}\cup\hat{\eta}\cup\hat{\delta}\end{eqnarray*}
where $\hat{\beta}$ is a multiset consisting of all odd sums with
same multiplicity $2^{\ell-k+3}$ formed by
$\beta_1,...,\beta_{k-2}$,  and

$$\begin{cases}
\hat{\gamma}=\{\gamma+\beta_1+\beta|\beta\in\hat{\beta}\}\\
\hat{\delta}=\{\delta+\beta_1+\beta|\beta\in\hat{\beta}\}\\
\hat{\varepsilon}=\{\gamma+\beta|\beta\in\hat{\beta}\}\\
\hat{\eta}=\{\delta+\beta|\beta\in\hat{\beta}\}\\
\hat{\lambda}=\{\gamma+\delta+\beta_1+\beta|\beta\in\hat{\beta}\}.\\
\end{cases}$$
\end{thm}

\begin{proof}
Clearly, all actions of $\mathcal{A}_{3\cdot 2^\ell}^k(4)$ determine
a unique connected regular graph with four vertices. Then
Theorem~\ref{abstr-graph} immediately follows from
Proposition~\ref{lm10} and Corollary~\ref{data3}.
\end{proof}

If $t>0$, then  generally $\Gamma_{(\Phi, M)}$ of an action $(\Phi,
M^n)$ in $\mathcal{A}_{3\cdot 2^\ell+t}^k(4)$ may not be uniquely
determined, so this leads to the difficulty of determining an
abstract 1-skeleton to be a colored graph of some action. The
problem seems to be quite complicated. Actually, even if $n=5$, as
seen in the 5-dimensional example $(\Lambda, M^5)$ of
Lemma~\ref{5-dim},  we still cannot determine which of six possible
abstract 1-skeleta in Figure~\ref{ab6}  is the colored graphs of
$(\Lambda, M^5)$.
However, we  can  characterize the  tangent representation set of
each action in $\mathcal{A}_{3\cdot 2^\ell+t}^k(4)$. By
Theorem~\ref{tangent} and  Proposition~\ref{lm12}, one has that

\begin{thm} \label{abstr-graph2}
If $t>0$, then the vertex-coloring set $\{\alpha(E_p), \alpha(E_q),
\alpha(E_r), \alpha(E_s)\}$ of $(\Gamma, \alpha)$ is the tangent
representation set of some $G$-action $(\Phi, M)$ in
$\mathcal{A}_{3\cdot 2^\ell+t}^k(4)$ if and only if  there is a
basis $\{\beta_1,...,\beta_{k-2},$ $\gamma, \delta\}$ of
$\Hom(({\Bbb Z}_2)^k,{\Bbb Z}_2)$ such that
\begin{eqnarray*}
\alpha(E_p)=\hat{\beta}\cup\hat{\gamma}\cup\hat{\delta}\cup\hat{\omega},
&& \alpha(E_q)=\hat{\beta}\cup\hat{\eta}\cup
\hat{\varepsilon}\cup\hat{\omega},\\\alpha(E_r)=
 \hat{\gamma}\cup\hat{\varepsilon}\cup\hat{\lambda}\cup\hat{\omega}, &&\alpha(E_s)=
 \hat{\delta}\cup\hat{\eta}\cup\hat{\lambda}\cup\hat{\omega}\end{eqnarray*}
 where $\hat{\beta}$ is a multiset consisting of all sums with  same multiplicity
$2^{\ell-k+3}$ formed by the odd number
 of elements of
$\beta_1,...,\beta_{k-2}$, and
$$\begin{cases}
\hat{\gamma}=\{\gamma+\beta_1+\beta|\beta\in\hat{\beta}\}\\
\hat{\delta}=\{\delta+\beta_1+\beta|\beta\in\hat{\beta}\}\\
\hat{\varepsilon}=\{\gamma+\beta|\beta\in\hat{\beta}\}\\
\hat{\eta}=\{\delta+\beta|\beta\in\hat{\beta}\}\\
\hat{\lambda}=\{\gamma+\delta+\beta_1+\beta|\beta\in\hat{\beta}\}\\
|\hat{\omega}|=t
\end{cases}$$
and  every element in $\hat{\omega}$ is chosen  in the $2^{k-3}$
different elements of $\{\gamma+\delta+\beta|\beta\in\hat{\beta}\}$
and has multiplicity at most $2^{\ell-k+4}$.
\end{thm}

Finally, combining Theorems~\ref{tangent},
\ref{abstr-graph}-\ref{abstr-graph2} and Proposition~\ref{lm12}, we
complete the proof of Theorem~\ref{thm6}.

\section{An observation on the minimum number of fixed points of
actions}\label{s9}

Define
$$m(n, k):=\min \big\{|M^G|\big| (\Phi, M^n)\in \mathcal{A}_n^k\big\}.$$
Theorem~\ref{thm2} has given an  estimation of the lower bound for
the number of fixed points of each action in  $\mathcal{A}_n^k$, so
one has that
$$m(n,k)\geq 1+\lceil{n\over{n-k+1}}\rceil.$$
Then it easily follows that
\begin{lem}\label{min bound}
If there is a $G$-action $(\Phi, M^n)$ in $\mathcal{A}_n^k$ such
that $|M^G|=1+\lceil{n\over{n-k+1}}\rceil$, then $m(n,k)=
1+\lceil{n\over{n-k+1}}\rceil$.
\end{lem}

By Lemma~\ref{min bound}, we know from Examples~\ref{5-1}-\ref{5-2}
that
$$m(n,k)=\begin{cases}
n+1 & \text{if } n=k\geq 2\\
3 & \text{if } n=2^\ell\geq 2^{k-1}\geq 2.
\end{cases}$$
On the other hand,  Theorems~\ref{thm3}(a) tells us that $$m(n,k)=3
\Longleftrightarrow   n=2^\ell\geq 2^{k-1}\geq 2.$$ This gives all
possible values of $n$ and $k$ for $m(n,k)=3$. Set
$$\mathcal{X}_1=\Big\{(n,k)\in{\Bbb N}^2\Big| n=2^\ell\geq 2^{k-1}\geq 2\Big\}$$
and
$$\mathcal{X}_2=\Big\{(n,k)\in{\Bbb N}^2\Big| k\geq 3, n\in \bigcup_{\ell \geq
k-3}[3\cdot 2^\ell, 5\cdot 2^\ell]\Big\}.$$ Since 3 is the minimum
possible value of $m(n,k)$, by Theorem~\ref{thm4}(a) we have that
$$m(n,k)=4\Longleftrightarrow(n,k)\in\mathcal{X}_2\setminus\mathcal{X}_1.$$ Combining the above
arguments, we have that
\begin{prop}\hskip .1cm\label{bo}
\begin{enumerate}
\item[$(a)$] $m(n,k)=3$
if and only if $(n,k)\in\mathcal{X}_1$.
\item[$(b)$] $m(n,k)=4$ if and only if
$(n,k)\in\mathcal{X}_2\setminus\mathcal{X}_1$.
\item[$(c)$] If $n=k\geq 2$, then $m(n,k)=n+1$.
\end{enumerate}
\end{prop}

By applying the $\Omega$-operation and the $\Delta$-operation to
$(\phi_4, {\Bbb R}P^4)$, one can obtain the $({\Bbb Z}_2)^k$-action
$(\Phi, M^n)$ fixing five isolated points with
$$(n, k)\in \mathcal{X}_3=\Big\{(n,k)\in {\Bbb N}^2\Big| n=2^\ell \geq 2^{k-2}\geq 4\Big\}.$$
It is easy to see that $\mathcal{X}_3\setminus(\mathcal{X}_1\cup
\mathcal{X}_2)$ is nonempty. For example, $(4,4)\in \mathcal{X}_3$
but $(4,4)\not\in \mathcal{X}_1\cup \mathcal{X}_2$. By
Proposition~\ref{bo} one has
\begin{cor}
If $(n,k)\in\mathcal{X}_3\setminus(\mathcal{X}_1\cup
\mathcal{X}_2)$, then $m(n,k)=5$.
\end{cor}

By Lemma~\ref{k=2}, it easily follows that

\begin{cor} Let $n\geq 2$ be even and let $n=2^{p_1}+\cdots+2^{p_r}$
with $0<p_1<\cdots <p_r$ be the 2-adic expansion of $n$. Then $m(n,
2)=3^r$.
\end{cor}

\end{document}

%% file: 1.pstex_t
\begin{picture}(0,0)%
\includegraphics{1.pstex}%
\end{picture}%
\setlength{\unitlength}{1500sp}%
\begingroup\makeatletter\ifx\SetFigFont\undefined%
\gdef\SetFigFont#1#2#3#4#5{%
  \reset@font\fontsize{#1}{#2pt}%
  \fontfamily{#3}\fontseries{#4}\fontshape{#5}%
  \selectfont}%
\fi\endgroup%
\begin{picture}(14700,2010)(526,-3961)
\end{picture}%

%% file: n2.pstex_t
\begin{picture}(0,0)%
\includegraphics{n2.pstex}%
\end{picture}%
\setlength{\unitlength}{1579sp}%
\begingroup\makeatletter\ifx\SetFigFont\undefined%
\gdef\SetFigFont#1#2#3#4#5{%
  \reset@font\fontsize{#1}{#2pt}%
  \fontfamily{#3}\fontseries{#4}\fontshape{#5}%
  \selectfont}%
\fi\endgroup%
\begin{picture}(9450,3447)(2026,-9733)
\put(4501,-7486){\makebox(0,0)[lb]{\smash{{\SetFigFont{7}{8.4}{\rmdefault}{\mddefault}{\updefault}$\rho_2$}}}}
\put(9376,-9136){\makebox(0,0)[lb]{\smash{{\SetFigFont{7}{8.4}{\rmdefault}{\mddefault}{\updefault}$\rho_2+\rho_3$}}}}
\put(8851,-9661){\makebox(0,0)[lb]{\smash{{\SetFigFont{7}{8.4}{\rmdefault}{\mddefault}{\updefault}The case $n=3$}}}}
\put(2551,-9661){\makebox(0,0)[lb]{\smash{{\SetFigFont{7}{8.4}{\rmdefault}{\mddefault}{\updefault}The case $n=2$}}}}
\put(9526,-7336){\makebox(0,0)[lb]{\smash{{\SetFigFont{7}{8.4}{\rmdefault}{\mddefault}{\updefault}$\rho_1$}}}}
\put(9001,-8086){\makebox(0,0)[lb]{\smash{{\SetFigFont{7}{8.4}{\rmdefault}{\mddefault}{\updefault}$\rho_2$}}}}
\put(3001,-9136){\makebox(0,0)[lb]{\smash{{\SetFigFont{7}{8.4}{\rmdefault}{\mddefault}{\updefault}$\rho_1+\rho_2$}}}}
\put(10351,-8011){\makebox(0,0)[lb]{\smash{{\SetFigFont{7}{8.4}{\rmdefault}{\mddefault}{\updefault}$\rho_3$}}}}
\put(2251,-7486){\makebox(0,0)[lb]{\smash{{\SetFigFont{7}{8.4}{\rmdefault}{\mddefault}{\updefault}$\rho_1$}}}}
\put(7801,-7336){\makebox(0,0)[lb]{\smash{{\SetFigFont{7}{8.4}{\rmdefault}{\mddefault}{\updefault}$\rho_1+\rho_2$}}}}
\put(10726,-7336){\makebox(0,0)[lb]{\smash{{\SetFigFont{7}{8.4}{\rmdefault}{\mddefault}{\updefault}$\rho_1+\rho_3$}}}}
\end{picture}%

%% file: n3.pstex_t
\begin{picture}(0,0)%
\includegraphics{n3.pstex}%
\end{picture}%
\setlength{\unitlength}{1579sp}%
\begingroup\makeatletter\ifx\SetFigFont\undefined%
\gdef\SetFigFont#1#2#3#4#5{%
  \reset@font\fontsize{#1}{#2pt}%
  \fontfamily{#3}\fontseries{#4}\fontshape{#5}%
  \selectfont}%
\fi\endgroup%
\begin{picture}(3975,3381)(8401,-6367)
\put(9976,-6286){\makebox(0,0)[lb]{\smash{{\SetFigFont{8}{9.6}{\rmdefault}{\mddefault}{\updefault}$\rho_1+\rho_2+\rho_3$}}}}
\put(10201,-5536){\makebox(0,0)[lb]{\smash{{\SetFigFont{8}{9.6}{\rmdefault}{\mddefault}{\updefault}$\rho_1+\rho_2$}}}}
\put(12151,-3961){\makebox(0,0)[lb]{\smash{{\SetFigFont{8}{9.6}{\rmdefault}{\mddefault}{\updefault}$\rho_2+\rho_3$}}}}
\put(8401,-3961){\makebox(0,0)[lb]{\smash{{\SetFigFont{8}{9.6}{\rmdefault}{\mddefault}{\updefault}$\rho_1+\rho_3$}}}}
\put(11401,-4411){\makebox(0,0)[lb]{\smash{{\SetFigFont{8}{9.6}{\rmdefault}{\mddefault}{\updefault}$\rho_2$}}}}
\put(9976,-4486){\makebox(0,0)[lb]{\smash{{\SetFigFont{8}{9.6}{\rmdefault}{\mddefault}{\updefault}$\rho_1$}}}}
\end{picture}%

%% file: n1.pstex_t
\begin{picture}(0,0)%
\includegraphics{n1.pstex}%
\end{picture}%
\setlength{\unitlength}{1381sp}%
\begingroup\makeatletter\ifx\SetFigFont\undefined%
\gdef\SetFigFont#1#2#3#4#5{%
  \reset@font\fontsize{#1}{#2pt}%
  \fontfamily{#3}\fontseries{#4}\fontshape{#5}%
  \selectfont}%
\fi\endgroup%
\begin{picture}(3715,3672)(3751,-9892)
\put(5026,-9811){\makebox(0,0)[lb]{\smash{{\SetFigFont{6}{7.2}{\rmdefault}{\mddefault}{\updefault}$\rho_2+\rho_3$}}}}
\put(7051,-9586){\makebox(0,0)[lb]{\smash{{\SetFigFont{6}{7.2}{\rmdefault}{\mddefault}{\updefault}$r$}}}}
\put(7051,-6736){\makebox(0,0)[lb]{\smash{{\SetFigFont{6}{7.2}{\rmdefault}{\mddefault}{\updefault}$q$}}}}
\put(7126,-8011){\makebox(0,0)[lb]{\smash{{\SetFigFont{6}{7.2}{\rmdefault}{\mddefault}{\updefault}$\rho_1+\rho_2$}}}}
\put(4951,-7261){\makebox(0,0)[lb]{\smash{{\SetFigFont{6}{7.2}{\rmdefault}{\mddefault}{\updefault}$\rho_2$}}}}
\put(5251,-6436){\makebox(0,0)[lb]{\smash{{\SetFigFont{6}{7.2}{\rmdefault}{\mddefault}{\updefault}$\rho_1$}}}}
\put(5026,-8836){\makebox(0,0)[lb]{\smash{{\SetFigFont{6}{7.2}{\rmdefault}{\mddefault}{\updefault}$\rho_1+\rho_3$}}}}
\put(3751,-8086){\makebox(0,0)[lb]{\smash{{\SetFigFont{6}{7.2}{\rmdefault}{\mddefault}{\updefault}$\rho_3$}}}}
\put(3826,-9586){\makebox(0,0)[lb]{\smash{{\SetFigFont{6}{7.2}{\rmdefault}{\mddefault}{\updefault}$s$}}}}
\put(3751,-6736){\makebox(0,0)[lb]{\smash{{\SetFigFont{6}{7.2}{\rmdefault}{\mddefault}{\updefault}$p$}}}}
\end{picture}%

%% file: n7.pstex_t
\begin{picture}(0,0)%
\includegraphics{n7.pstex}%
\end{picture}%
\setlength{\unitlength}{1776sp}%
\begingroup\makeatletter\ifx\SetFigFont\undefined%
\gdef\SetFigFont#1#2#3#4#5{%
  \reset@font\fontsize{#1}{#2pt}%
  \fontfamily{#3}\fontseries{#4}\fontshape{#5}%
  \selectfont}%
\fi\endgroup%
\begin{picture}(3933,2405)(-299,-4183)
\put(-299,-4111){\makebox(0,0)[lb]{\smash{{\SetFigFont{7}{8.4}{\rmdefault}{\mddefault}{\updefault}The connected regular graph $\Gamma$ of valence 4}}}}
\end{picture}%

%% file: ld3.pstex_t
\begin{picture}(0,0)%
\includegraphics{ld3.pstex}%
\end{picture}%
\setlength{\unitlength}{1934sp}%
\begingroup\makeatletter\ifx\SetFigFont\undefined%
\gdef\SetFigFont#1#2#3#4#5{%
  \reset@font\fontsize{#1}{#2pt}%
  \fontfamily{#3}\fontseries{#4}\fontshape{#5}%
  \selectfont}%
\fi\endgroup%
\begin{picture}(10383,4356)(751,-5149)
\put(2101,-2236){\makebox(0,0)[lb]{\smash{{\SetFigFont{7}{8.4}{\rmdefault}{\mddefault}{\updefault}$\rho_1+\rho_3$}}}}
\put(2626,-2761){\makebox(0,0)[lb]{\smash{{\SetFigFont{7}{8.4}{\rmdefault}{\mddefault}{\updefault}$\rho_1+\rho_2$}}}}
\put(4051,-2761){\makebox(0,0)[lb]{\smash{{\SetFigFont{7}{8.4}{\rmdefault}{\mddefault}{\updefault}$\rho_3$}}}}
\put(1576,-3886){\makebox(0,0)[lb]{\smash{{\SetFigFont{7}{8.4}{\rmdefault}{\mddefault}{\updefault}$\rho_2+\rho_3$}}}}
\put(1351,-4561){\makebox(0,0)[lb]{\smash{{\SetFigFont{7}{8.4}{\rmdefault}{\mddefault}{\updefault}$\rho_1+\rho_2+\rho_3$}}}}
\put(5701,-961){\makebox(0,0)[lb]{\smash{{\SetFigFont{7}{8.4}{\rmdefault}{\mddefault}{\updefault}$\rho_1+\rho_2+\rho_3$}}}}
\put(5326,-1786){\makebox(0,0)[lb]{\smash{{\SetFigFont{7}{8.4}{\rmdefault}{\mddefault}{\updefault}$\rho_1$}}}}
\put(4651,-2386){\makebox(0,0)[lb]{\smash{{\SetFigFont{7}{8.4}{\rmdefault}{\mddefault}{\updefault}$\rho_2$}}}}
\put(5701,-2836){\makebox(0,0)[lb]{\smash{{\SetFigFont{7}{8.4}{\rmdefault}{\mddefault}{\updefault}$\rho_1+\rho_2$}}}}
\put(6001,-4561){\makebox(0,0)[lb]{\smash{{\SetFigFont{7}{8.4}{\rmdefault}{\mddefault}{\updefault}$\rho_1+\rho_2+\rho_3$}}}}
\put(9376,-1786){\makebox(0,0)[lb]{\smash{{\SetFigFont{7}{8.4}{\rmdefault}{\mddefault}{\updefault}$\rho_1$}}}}
\put(6976,-2611){\makebox(0,0)[lb]{\smash{{\SetFigFont{7}{8.4}{\rmdefault}{\mddefault}{\updefault}$\rho_1+\rho_2+\rho_3$}}}}
\put(8401,-2836){\makebox(0,0)[lb]{\smash{{\SetFigFont{7}{8.4}{\rmdefault}{\mddefault}{\updefault}$\rho_3$}}}}
\put(8851,-2461){\makebox(0,0)[lb]{\smash{{\SetFigFont{7}{8.4}{\rmdefault}{\mddefault}{\updefault}$\rho_2$}}}}
\put(9901,-2836){\makebox(0,0)[lb]{\smash{{\SetFigFont{7}{8.4}{\rmdefault}{\mddefault}{\updefault}$\rho_1+\rho_2$}}}}
\put(10801,-2686){\makebox(0,0)[lb]{\smash{{\SetFigFont{7}{8.4}{\rmdefault}{\mddefault}{\updefault}$\rho_1+\rho_2+\rho_3$}}}}
\put(9226,-3961){\makebox(0,0)[lb]{\smash{{\SetFigFont{7}{8.4}{\rmdefault}{\mddefault}{\updefault}$\rho_2+\rho_3$}}}}
\put(9526,-2236){\makebox(0,0)[lb]{\smash{{\SetFigFont{7}{8.4}{\rmdefault}{\mddefault}{\updefault}$\rho_1+\rho_3$}}}}
\put(5401,-2236){\makebox(0,0)[lb]{\smash{{\SetFigFont{7}{8.4}{\rmdefault}{\mddefault}{\updefault}$\rho_1+\rho_3$}}}}
\put(4951,-3511){\makebox(0,0)[lb]{\smash{{\SetFigFont{7}{8.4}{\rmdefault}{\mddefault}{\updefault}$\rho_2+\rho_3$}}}}
\put(1351,-2461){\makebox(0,0)[lb]{\smash{{\SetFigFont{7}{8.4}{\rmdefault}{\mddefault}{\updefault}$\rho_2$}}}}
\put(1501,-5086){\makebox(0,0)[lb]{\smash{{\SetFigFont{7}{8.4}{\rmdefault}{\mddefault}{\updefault}$(\Gamma, \alpha_1)$}}}}
\put(5026,-5086){\makebox(0,0)[lb]{\smash{{\SetFigFont{7}{8.4}{\rmdefault}{\mddefault}{\updefault}$(\Gamma, \alpha_2)$}}}}
\put(9451,-5086){\makebox(0,0)[lb]{\smash{{\SetFigFont{7}{8.4}{\rmdefault}{\mddefault}{\updefault}$(\Gamma, \alpha_3)$}}}}
\put(10651,-3811){\makebox(0,0)[lb]{\smash{{\SetFigFont{7}{8.4}{\rmdefault}{\mddefault}{\updefault}$r$}}}}
\put(6451,-3811){\makebox(0,0)[lb]{\smash{{\SetFigFont{7}{8.4}{\rmdefault}{\mddefault}{\updefault}$r$}}}}
\put(8401,-3811){\makebox(0,0)[lb]{\smash{{\SetFigFont{7}{8.4}{\rmdefault}{\mddefault}{\updefault}$s$}}}}
\put(3076,-1786){\makebox(0,0)[lb]{\smash{{\SetFigFont{7}{8.4}{\rmdefault}{\mddefault}{\updefault}$q$}}}}
\put(976,-1786){\makebox(0,0)[lb]{\smash{{\SetFigFont{7}{8.4}{\rmdefault}{\mddefault}{\updefault}$p$}}}}
\put(901,-3811){\makebox(0,0)[lb]{\smash{{\SetFigFont{7}{8.4}{\rmdefault}{\mddefault}{\updefault}$s$}}}}
\put(3151,-3811){\makebox(0,0)[lb]{\smash{{\SetFigFont{7}{8.4}{\rmdefault}{\mddefault}{\updefault}$r$}}}}
\put(4201,-1786){\makebox(0,0)[lb]{\smash{{\SetFigFont{7}{8.4}{\rmdefault}{\mddefault}{\updefault}$p$}}}}
\put(4201,-3811){\makebox(0,0)[lb]{\smash{{\SetFigFont{7}{8.4}{\rmdefault}{\mddefault}{\updefault}$s$}}}}
\put(6376,-1786){\makebox(0,0)[lb]{\smash{{\SetFigFont{7}{8.4}{\rmdefault}{\mddefault}{\updefault}$q$}}}}
\put(8401,-1786){\makebox(0,0)[lb]{\smash{{\SetFigFont{7}{8.4}{\rmdefault}{\mddefault}{\updefault}$p$}}}}
\put(10576,-1786){\makebox(0,0)[lb]{\smash{{\SetFigFont{7}{8.4}{\rmdefault}{\mddefault}{\updefault}$q$}}}}
\put(1426,-1111){\makebox(0,0)[lb]{\smash{{\SetFigFont{7}{8.4}{\rmdefault}{\mddefault}{\updefault}$\rho_1+\rho_2+\rho_3$}}}}
\put(1951,-1786){\makebox(0,0)[lb]{\smash{{\SetFigFont{7}{8.4}{\rmdefault}{\mddefault}{\updefault}$\rho_1$}}}}
\put(751,-2761){\makebox(0,0)[lb]{\smash{{\SetFigFont{7}{8.4}{\rmdefault}{\mddefault}{\updefault}$\rho_3$}}}}
\end{picture}%

%% file: n6.pstex_t
\begin{picture}(0,0)%
\includegraphics{n6.pstex}%
\end{picture}%
\setlength{\unitlength}{1776sp}%
\begingroup\makeatletter\ifx\SetFigFont\undefined%
\gdef\SetFigFont#1#2#3#4#5{%
  \reset@font\fontsize{#1}{#2pt}%
  \fontfamily{#3}\fontseries{#4}\fontshape{#5}%
  \selectfont}%
\fi\endgroup%
\begin{picture}(10966,3396)(568,-4624)
\put(9601,-4486){\makebox(0,0)[lb]{\smash{{\SetFigFont{7}{8.4}{\rmdefault}{\mddefault}{\updefault}$\Gamma_2$}}}}
\put(2026,-4561){\makebox(0,0)[lb]{\smash{{\SetFigFont{7}{8.4}{\rmdefault}{\mddefault}{\updefault}$\Gamma_1$}}}}
\end{picture}%

%% file: n4.pstex_t
\begin{picture}(0,0)%
\includegraphics{n4.pstex}%
\end{picture}%
\setlength{\unitlength}{1934sp}%
\begingroup\makeatletter\ifx\SetFigFont\undefined%
\gdef\SetFigFont#1#2#3#4#5{%
  \reset@font\fontsize{#1}{#2pt}%
  \fontfamily{#3}\fontseries{#4}\fontshape{#5}%
  \selectfont}%
\fi\endgroup%
\begin{picture}(11381,9081)(-74,-9499)
\put(1651,-8461){\makebox(0,0)[lb]{\smash{{\SetFigFont{7}{8.4}{\rmdefault}{\mddefault}{\updefault}$\rho_2+\rho_3$}}}}
\put(5251,-8011){\makebox(0,0)[lb]{\smash{{\SetFigFont{7}{8.4}{\rmdefault}{\mddefault}{\updefault}$\rho_2+\rho_3$}}}}
\put(1951,-6211){\makebox(0,0)[lb]{\smash{{\SetFigFont{7}{8.4}{\rmdefault}{\mddefault}{\updefault}$\rho_1$}}}}
\put(5626,-6286){\makebox(0,0)[lb]{\smash{{\SetFigFont{7}{8.4}{\rmdefault}{\mddefault}{\updefault}$\rho_1$}}}}
\put(9451,-6211){\makebox(0,0)[lb]{\smash{{\SetFigFont{7}{8.4}{\rmdefault}{\mddefault}{\updefault}$\rho_1$}}}}
\put(9526,-1711){\makebox(0,0)[lb]{\smash{{\SetFigFont{7}{8.4}{\rmdefault}{\mddefault}{\updefault}$\rho_1$}}}}
\put(9151,-4486){\makebox(0,0)[lb]{\smash{{\SetFigFont{7}{8.4}{\rmdefault}{\mddefault}{\updefault}$(\Gamma_1, \alpha_3)$}}}}
\put(9151,-9436){\makebox(0,0)[lb]{\smash{{\SetFigFont{7}{8.4}{\rmdefault}{\mddefault}{\updefault}$(\Gamma_2, \alpha_6)$}}}}
\put(10576,-7411){\makebox(0,0)[lb]{\smash{{\SetFigFont{7}{8.4}{\rmdefault}{\mddefault}{\updefault}$\rho_1+\rho_2$}}}}
\put(6001,-7336){\makebox(0,0)[lb]{\smash{{\SetFigFont{7}{8.4}{\rmdefault}{\mddefault}{\updefault}$\rho_1+\rho_2$}}}}
\put(5176,-9436){\makebox(0,0)[lb]{\smash{{\SetFigFont{7}{8.4}{\rmdefault}{\mddefault}{\updefault}$(\Gamma_2, \alpha_5)$}}}}
\put(1876,-9436){\makebox(0,0)[lb]{\smash{{\SetFigFont{7}{8.4}{\rmdefault}{\mddefault}{\updefault}$(\Gamma_2, \alpha_4)$}}}}
\put(3076,-1786){\makebox(0,0)[lb]{\smash{{\SetFigFont{7}{8.4}{\rmdefault}{\mddefault}{\updefault}$q$}}}}
\put(976,-1786){\makebox(0,0)[lb]{\smash{{\SetFigFont{7}{8.4}{\rmdefault}{\mddefault}{\updefault}$p$}}}}
\put(901,-3811){\makebox(0,0)[lb]{\smash{{\SetFigFont{7}{8.4}{\rmdefault}{\mddefault}{\updefault}$s$}}}}
\put(3151,-3811){\makebox(0,0)[lb]{\smash{{\SetFigFont{7}{8.4}{\rmdefault}{\mddefault}{\updefault}$r$}}}}
\put(4201,-1786){\makebox(0,0)[lb]{\smash{{\SetFigFont{7}{8.4}{\rmdefault}{\mddefault}{\updefault}$p$}}}}
\put(4201,-3811){\makebox(0,0)[lb]{\smash{{\SetFigFont{7}{8.4}{\rmdefault}{\mddefault}{\updefault}$s$}}}}
\put(6376,-1786){\makebox(0,0)[lb]{\smash{{\SetFigFont{7}{8.4}{\rmdefault}{\mddefault}{\updefault}$q$}}}}
\put(8401,-1786){\makebox(0,0)[lb]{\smash{{\SetFigFont{7}{8.4}{\rmdefault}{\mddefault}{\updefault}$p$}}}}
\put(10576,-1786){\makebox(0,0)[lb]{\smash{{\SetFigFont{7}{8.4}{\rmdefault}{\mddefault}{\updefault}$q$}}}}
\put(8476,-3886){\makebox(0,0)[lb]{\smash{{\SetFigFont{7}{8.4}{\rmdefault}{\mddefault}{\updefault}$s$}}}}
\put(1951,-1786){\makebox(0,0)[lb]{\smash{{\SetFigFont{7}{8.4}{\rmdefault}{\mddefault}{\updefault}$\rho_1$}}}}
\put(751,-2761){\makebox(0,0)[lb]{\smash{{\SetFigFont{7}{8.4}{\rmdefault}{\mddefault}{\updefault}$\rho_3$}}}}
\put(2101,-2236){\makebox(0,0)[lb]{\smash{{\SetFigFont{7}{8.4}{\rmdefault}{\mddefault}{\updefault}$\rho_1+\rho_3$}}}}
\put(2626,-2761){\makebox(0,0)[lb]{\smash{{\SetFigFont{7}{8.4}{\rmdefault}{\mddefault}{\updefault}$\rho_1+\rho_2$}}}}
\put(4051,-2761){\makebox(0,0)[lb]{\smash{{\SetFigFont{7}{8.4}{\rmdefault}{\mddefault}{\updefault}$\rho_3$}}}}
\put(1576,-3886){\makebox(0,0)[lb]{\smash{{\SetFigFont{7}{8.4}{\rmdefault}{\mddefault}{\updefault}$\rho_2+\rho_3$}}}}
\put(5326,-1786){\makebox(0,0)[lb]{\smash{{\SetFigFont{7}{8.4}{\rmdefault}{\mddefault}{\updefault}$\rho_1$}}}}
\put(4651,-2386){\makebox(0,0)[lb]{\smash{{\SetFigFont{7}{8.4}{\rmdefault}{\mddefault}{\updefault}$\rho_2$}}}}
\put(5701,-2836){\makebox(0,0)[lb]{\smash{{\SetFigFont{7}{8.4}{\rmdefault}{\mddefault}{\updefault}$\rho_1+\rho_2$}}}}
\put(8401,-2836){\makebox(0,0)[lb]{\smash{{\SetFigFont{7}{8.4}{\rmdefault}{\mddefault}{\updefault}$\rho_3$}}}}
\put(8851,-2461){\makebox(0,0)[lb]{\smash{{\SetFigFont{7}{8.4}{\rmdefault}{\mddefault}{\updefault}$\rho_2$}}}}
\put(9901,-2836){\makebox(0,0)[lb]{\smash{{\SetFigFont{7}{8.4}{\rmdefault}{\mddefault}{\updefault}$\rho_1+\rho_2$}}}}
\put(9526,-2236){\makebox(0,0)[lb]{\smash{{\SetFigFont{7}{8.4}{\rmdefault}{\mddefault}{\updefault}$\rho_1+\rho_3$}}}}
\put(5401,-2236){\makebox(0,0)[lb]{\smash{{\SetFigFont{7}{8.4}{\rmdefault}{\mddefault}{\updefault}$\rho_1+\rho_3$}}}}
\put(4951,-3511){\makebox(0,0)[lb]{\smash{{\SetFigFont{7}{8.4}{\rmdefault}{\mddefault}{\updefault}$\rho_2+\rho_3$}}}}
\put(1351,-2461){\makebox(0,0)[lb]{\smash{{\SetFigFont{7}{8.4}{\rmdefault}{\mddefault}{\updefault}$\rho_2$}}}}
\put(751,-7261){\makebox(0,0)[lb]{\smash{{\SetFigFont{7}{8.4}{\rmdefault}{\mddefault}{\updefault}$\rho_3$}}}}
\put(1276,-6886){\makebox(0,0)[lb]{\smash{{\SetFigFont{7}{8.4}{\rmdefault}{\mddefault}{\updefault}$\rho_2$}}}}
\put(2251,-6811){\makebox(0,0)[lb]{\smash{{\SetFigFont{7}{8.4}{\rmdefault}{\mddefault}{\updefault}$\rho_1+\rho_3$}}}}
\put(751,-6436){\makebox(0,0)[lb]{\smash{{\SetFigFont{7}{8.4}{\rmdefault}{\mddefault}{\updefault}$p$}}}}
\put(751,-8311){\makebox(0,0)[lb]{\smash{{\SetFigFont{7}{8.4}{\rmdefault}{\mddefault}{\updefault}$s$}}}}
\put(4426,-6286){\makebox(0,0)[lb]{\smash{{\SetFigFont{7}{8.4}{\rmdefault}{\mddefault}{\updefault}$p$}}}}
\put(8476,-8461){\makebox(0,0)[lb]{\smash{{\SetFigFont{7}{8.4}{\rmdefault}{\mddefault}{\updefault}$s$}}}}
\put(4426,-7261){\makebox(0,0)[lb]{\smash{{\SetFigFont{7}{8.4}{\rmdefault}{\mddefault}{\updefault}$\rho_3$}}}}
\put(8326,-7261){\makebox(0,0)[lb]{\smash{{\SetFigFont{7}{8.4}{\rmdefault}{\mddefault}{\updefault}$\rho_3$}}}}
\put(8851,-6886){\makebox(0,0)[lb]{\smash{{\SetFigFont{7}{8.4}{\rmdefault}{\mddefault}{\updefault}$\rho_2$}}}}
\put(9751,-6811){\makebox(0,0)[lb]{\smash{{\SetFigFont{7}{8.4}{\rmdefault}{\mddefault}{\updefault}$\rho_1+\rho_3$}}}}
\put(4951,-6886){\makebox(0,0)[lb]{\smash{{\SetFigFont{7}{8.4}{\rmdefault}{\mddefault}{\updefault}$\rho_2$}}}}
\put(5776,-6736){\makebox(0,0)[lb]{\smash{{\SetFigFont{7}{8.4}{\rmdefault}{\mddefault}{\updefault}$\rho_1+\rho_3$}}}}
\put(6751,-6361){\makebox(0,0)[lb]{\smash{{\SetFigFont{7}{8.4}{\rmdefault}{\mddefault}{\updefault}$q$}}}}
\put(10651,-6361){\makebox(0,0)[lb]{\smash{{\SetFigFont{7}{8.4}{\rmdefault}{\mddefault}{\updefault}$q$}}}}
\put(4351,-8311){\makebox(0,0)[lb]{\smash{{\SetFigFont{7}{8.4}{\rmdefault}{\mddefault}{\updefault}$s$}}}}
\put(8326,-6286){\makebox(0,0)[lb]{\smash{{\SetFigFont{7}{8.4}{\rmdefault}{\mddefault}{\updefault}$p$}}}}
\put(3076,-6286){\makebox(0,0)[lb]{\smash{{\SetFigFont{7}{8.4}{\rmdefault}{\mddefault}{\updefault}$q$}}}}
\put(7651,-586){\makebox(0,0)[lb]{\smash{{\SetFigFont{7}{8.4}{\rmdefault}{\mddefault}{\updefault}$\rho_1+\rho_2+\rho_3$}}}}
\put(301,-661){\makebox(0,0)[lb]{\smash{{\SetFigFont{7}{8.4}{\rmdefault}{\mddefault}{\updefault}$\rho_1+\rho_2+\rho_3$}}}}
\put(9376,-4936){\makebox(0,0)[lb]{\smash{{\SetFigFont{7}{8.4}{\rmdefault}{\mddefault}{\updefault}$\rho_1+\rho_2+\rho_3$}}}}
\put(1201,-4711){\makebox(0,0)[lb]{\smash{{\SetFigFont{7}{8.4}{\rmdefault}{\mddefault}{\updefault}$(\Gamma_1, \alpha_1)$}}}}
\put(3601,-4936){\makebox(0,0)[lb]{\smash{{\SetFigFont{7}{8.4}{\rmdefault}{\mddefault}{\updefault}$\rho_1+\rho_2+\rho_3$}}}}
\put(5251,-4636){\makebox(0,0)[lb]{\smash{{\SetFigFont{7}{8.4}{\rmdefault}{\mddefault}{\updefault}$(\Gamma_1, \alpha_2)$}}}}
\put(6901,-9361){\makebox(0,0)[lb]{\smash{{\SetFigFont{7}{8.4}{\rmdefault}{\mddefault}{\updefault}$\rho_1+\rho_2+\rho_3$}}}}
\put(-74,-5686){\makebox(0,0)[lb]{\smash{{\SetFigFont{7}{8.4}{\rmdefault}{\mddefault}{\updefault}$\rho_1+\rho_2+\rho_3$}}}}
\put(2476,-7336){\makebox(0,0)[lb]{\smash{{\SetFigFont{7}{8.4}{\rmdefault}{\mddefault}{\updefault}$\rho_1+\rho_2$}}}}
\put(9151,-3586){\makebox(0,0)[lb]{\smash{{\SetFigFont{7}{8.4}{\rmdefault}{\mddefault}{\updefault}$\rho_2+\rho_3$}}}}
\put(9151,-8011){\makebox(0,0)[lb]{\smash{{\SetFigFont{7}{8.4}{\rmdefault}{\mddefault}{\updefault}$\rho_2+\rho_3$}}}}
\put(3151,-8311){\makebox(0,0)[lb]{\smash{{\SetFigFont{7}{8.4}{\rmdefault}{\mddefault}{\updefault}$r$}}}}
\put(10651,-3736){\makebox(0,0)[lb]{\smash{{\SetFigFont{7}{8.4}{\rmdefault}{\mddefault}{\updefault}$r$}}}}
\put(6751,-8236){\makebox(0,0)[lb]{\smash{{\SetFigFont{7}{8.4}{\rmdefault}{\mddefault}{\updefault}$r$}}}}
\put(10651,-8311){\makebox(0,0)[lb]{\smash{{\SetFigFont{7}{8.4}{\rmdefault}{\mddefault}{\updefault}$r$}}}}
\put(6451,-3811){\makebox(0,0)[lb]{\smash{{\SetFigFont{7}{8.4}{\rmdefault}{\mddefault}{\updefault}$r$}}}}
\end{picture}%

%% file: n8.pstex_t
\begin{picture}(0,0)%
\includegraphics{n8.pstex}%
\end{picture}%
\setlength{\unitlength}{1579sp}%
\begingroup\makeatletter\ifx\SetFigFont\undefined%
\gdef\SetFigFont#1#2#3#4#5{%
  \reset@font\fontsize{#1}{#2pt}%
  \fontfamily{#3}\fontseries{#4}\fontshape{#5}%
  \selectfont}%
\fi\endgroup%
\begin{picture}(3594,3397)(8926,-5228)
\put(12151,-2236){\makebox(0,0)[lb]{\smash{{\SetFigFont{6}{7.2}{\rmdefault}{\mddefault}{\updefault}$q$}}}}
\put(12151,-5011){\makebox(0,0)[lb]{\smash{{\SetFigFont{6}{7.2}{\rmdefault}{\mddefault}{\updefault}$r$}}}}
\put(10651,-2011){\makebox(0,0)[lb]{\smash{{\SetFigFont{6}{7.2}{\rmdefault}{\mddefault}{\updefault}$\hat{\beta}$}}}}
\put(9976,-2611){\makebox(0,0)[lb]{\smash{{\SetFigFont{6}{7.2}{\rmdefault}{\mddefault}{\updefault}$\hat{\gamma}$}}}}
\put(11326,-3136){\makebox(0,0)[lb]{\smash{{\SetFigFont{6}{7.2}{\rmdefault}{\mddefault}{\updefault}$\hat{\eta}$}}}}
\put(9001,-5086){\makebox(0,0)[lb]{\smash{{\SetFigFont{6}{7.2}{\rmdefault}{\mddefault}{\updefault}$s$}}}}
\put(8926,-2236){\makebox(0,0)[lb]{\smash{{\SetFigFont{6}{7.2}{\rmdefault}{\mddefault}{\updefault}$p$}}}}
\put(8926,-3586){\makebox(0,0)[lb]{\smash{{\SetFigFont{6}{7.2}{\rmdefault}{\mddefault}{\updefault}$\hat{\delta}$}}}}
\put(10651,-5161){\makebox(0,0)[lb]{\smash{{\SetFigFont{6}{7.2}{\rmdefault}{\mddefault}{\updefault}$\hat{\lambda}$}}}}
\put(12151,-3586){\makebox(0,0)[lb]{\smash{{\SetFigFont{6}{7.2}{\rmdefault}{\mddefault}{\updefault}$\hat{\varepsilon}$}}}}
\end{picture}%

%% file: ld4.pstex_t
\begin{picture}(0,0)%
\includegraphics{ld4.pstex}%
\end{picture}%
\setlength{\unitlength}{1737sp}%
\begingroup\makeatletter\ifx\SetFigFont\undefined%
\gdef\SetFigFont#1#2#3#4#5{%
  \reset@font\fontsize{#1}{#2pt}%
  \fontfamily{#3}\fontseries{#4}\fontshape{#5}%
  \selectfont}%
\fi\endgroup%
\begin{picture}(6075,4272)(1651,-5017)
\put(3901,-1486){\makebox(0,0)[lb]{\smash{{\SetFigFont{8}{9.6}{\rmdefault}{\mddefault}{\updefault}$p$}}}}
\put(3901,-4486){\makebox(0,0)[lb]{\smash{{\SetFigFont{8}{9.6}{\rmdefault}{\mddefault}{\updefault}$s$}}}}
\put(7051,-4411){\makebox(0,0)[lb]{\smash{{\SetFigFont{8}{9.6}{\rmdefault}{\mddefault}{\updefault}$r$}}}}
\put(5026,-4261){\makebox(0,0)[lb]{\smash{{\SetFigFont{8}{9.6}{\rmdefault}{\mddefault}{\updefault}$\rho_2+\rho_3$}}}}
\put(7726,-2311){\makebox(0,0)[lb]{\smash{{\SetFigFont{8}{9.6}{\rmdefault}{\mddefault}{\updefault}$\rho_1+\rho_3$}}}}
\put(6151,-2986){\makebox(0,0)[lb]{\smash{{\SetFigFont{8}{9.6}{\rmdefault}{\mddefault}{\updefault}$\rho_1+\rho_2$}}}}
\put(2551,-2986){\makebox(0,0)[lb]{\smash{{\SetFigFont{8}{9.6}{\rmdefault}{\mddefault}{\updefault}$\rho_3+\rho_4$}}}}
\put(4126,-2986){\makebox(0,0)[lb]{\smash{{\SetFigFont{8}{9.6}{\rmdefault}{\mddefault}{\updefault}$\rho_3$}}}}
\put(7726,-3661){\makebox(0,0)[lb]{\smash{{\SetFigFont{8}{9.6}{\rmdefault}{\mddefault}{\updefault}$\rho_2+\rho_4$}}}}
\put(3001,-2311){\makebox(0,0)[lb]{\smash{{\SetFigFont{8}{9.6}{\rmdefault}{\mddefault}{\updefault}$\rho_2$}}}}
\put(1651,-3661){\makebox(0,0)[lb]{\smash{{\SetFigFont{8}{9.6}{\rmdefault}{\mddefault}{\updefault}$\rho_1+\rho_3+\rho_4$}}}}
\put(7576,-2986){\makebox(0,0)[lb]{\smash{{\SetFigFont{8}{9.6}{\rmdefault}{\mddefault}{\updefault}$\rho_1+\rho_2+\rho_4$}}}}
\put(6976,-1411){\makebox(0,0)[lb]{\smash{{\SetFigFont{8}{9.6}{\rmdefault}{\mddefault}{\updefault}$q$}}}}
\put(5026,-961){\makebox(0,0)[lb]{\smash{{\SetFigFont{8}{9.6}{\rmdefault}{\mddefault}{\updefault}$\rho_1+\rho_4$}}}}
\put(4726,-4936){\makebox(0,0)[lb]{\smash{{\SetFigFont{8}{9.6}{\rmdefault}{\mddefault}{\updefault}$\rho_2+\rho_3+\rho_4$}}}}
\put(5401,-1786){\makebox(0,0)[lb]{\smash{{\SetFigFont{8}{9.6}{\rmdefault}{\mddefault}{\updefault}$\rho_1$}}}}
\end{picture}%